\documentclass[11pt, letterpaper, reqno]{article}
\usepackage{times}
\usepackage[letterpaper,margin=1in]{geometry}
\usepackage{amsthm}
\usepackage{amssymb}
\usepackage{amsmath,amsfonts}
\usepackage{mathrsfs}
\usepackage{cases}
\usepackage{latexsym,bm}
\usepackage{indentfirst}
\usepackage{color}
\usepackage{ifpdf}
\usepackage{graphicx}
\usepackage{psfrag}
\usepackage{dsfont}
\usepackage[ruled,vlined]{algorithm2e}
\usepackage{enumerate}
\usepackage[
pdfauthor={},
pdftitle={},
pdfstartview=XYZ,
bookmarks=true,
colorlinks=true,
linkcolor=blue,
urlcolor=blue,
citecolor=blue,
bookmarks=true,
linktocpage=true,
hyperindex=true
]{hyperref}

\usepackage{pgf,tikz,pgfplots}
\usepackage{tikz-3dplot}
\usepackage{pifont}
\usepackage{refcount}
\newtheorem{THM}{\textbf{Theorem}}[section]

\newtheorem{LEM}[THM]{\textbf{Lemma}}
\newtheorem{CLA}[THM]{\textbf{Claim}}
\newtheorem{CON}[THM]{\textbf{Conjecture}}

\newtheorem{OBV}{\textbf{Observation}}

\newcommand{\pf}{\noindent\textbf{Proof}.\quad}

\newtheorem*{THM1}{\textbf{Theorem 1.3}}

\newtheorem{Pro}{\textbf{Procedure}}

\newcommand{\ve}{\varepsilon}

\newcommand{\CC}{\mathcal{C}}

\DeclareMathOperator{\df}{def}
\newcommand{\pbar}{\overline{\varphi}}

\newcommand{\hl}[1]{{\bf \noindent \underline{#1}}}

\begin{document}
	\title{Towards the Overfull  Conjecture}
	\author	{Songling Shan \thanks{Partially supported by NSF  grant 
			DMS-2345869.}\\ 
		\medskip  Auburn University, Auburn, AL 36849\\
		\medskip 
		{\tt szs0398@auburn.edu}
	}
	
	\date{September 5, 2024}
	\maketitle

	\emph{\textbf{Abstract}.}
	Let $G$ be a simple graph with  maximum degree  denoted as $\Delta(G)$. An overfull subgraph $H$ of $G$ is a subgraph satisfying the condition $|E(H)| > \Delta(G)\lfloor \frac{1}{2}|V(H)| \rfloor$. In 1986, Chetwynd and Hilton proposed the Overfull Conjecture, stating that a  graph $G$ with  maximum degree $\Delta(G)> \frac{1}{3}|V(G)|$  has  chromatic index equal to  $\Delta(G)$ if and only if it does not contain any overfull subgraph.
	The Overfull Conjecture has many implications.   For example,  it implies  a  polynomial-time algorithm 
	for determining the chromatic index of   graphs $G$ with  $\Delta(G) > \frac{1}{3}|V(G)|$, and implies several longstanding 
	conjectures   in the area of graph edge colorings. 
	In this paper, we make  the first breakthrough  towards the conjecture when not imposing a minimum degree condition on the graph:  for any  $0<\ve  \le \frac{1}{14}$, there exists  a positive  integer $n_0$ 
	such that if $G$ is a  graph on $n\ge n_0$ vertices with $\Delta(G) \ge (1-\ve)n$, then the 
	Overfull Conjecture holds for $G$.  The previous best result in this direction,   due to Chetwynd and Hilton from 1989, asserts the conjecture 
	for graphs $G$ with $\Delta(G) \ge |V(G)|-3$.   Our result also 
	implies the Average Degree Conjecture of Vizing from 1968 
	for the same class of graphs $G$.

	
	\vspace{2mm}
	
	\section{Introduction}

	Graphs in this paper have neither multiple edge nor loop, whereas  multigraphs may have multiple edges but no loop.  
	However, for simpler terminology, for a multigraph $G$ and $H\subseteq G$, we still call $H$ 
	a subgraph of $G$ instead of a sub-multigraph. 
	As we will need multigraphs in the proof of the main result, we define terminologies with respect to multigraphs for generality. 
	For two integers $p,q$, let $[p,q]=\{ i\in \mathbb{Z} \,:\, p \le i \le q\}$.
	For a nonnegative integer $k$,  
	an  \emph{edge $k$-coloring} of a  multigraph $G$ is a mapping $\varphi$ from $E(G)$ to the set of integers
	$[1,k]$, called  \emph{colors}, such that  no two adjacent edges receive the same color with respect to $\varphi$.   
	Each set of edges colored by the same color under $\varphi$ is a \emph{color class} of $G$ with respect to $\varphi$. 
	The \emph{chromatic index} of $G$, denoted $\chi'(G)$, is  the smallest integer $k$ so that $G$ has an edge $k$-coloring.  
	A  multigraph $G$ with $|E(G)|>\Delta(G) \lfloor \frac{1}{2} |V(G)| \rfloor$ 
	is \emph{overfull}  and a  subgraph $H$ of $G$ with  $|E(H)|>\Delta(G) \lfloor \frac{1}{2} |V(H)| \rfloor$ 
	is called a \emph{$\Delta(G)$-overfull subgraph}.  Since each color class of $G$ is a matching and each matching in $G$
	has size at most $\lfloor  \frac{1}{2}|V(G)|\rfloor$, it is clear that if $G$ contains a $\Delta(G)$-overfull subgraph $H$, then $\chi'(G) \ge \chi'(H) \ge  |E(H)|/ \lfloor \frac{1}{2} |V(H)|\rfloor>\Delta(H)=\Delta(G)$.

	In 1960's, Gupta~\cite{Gupta-67}  and, independently, Vizing~\cite{Vizing-2-classes}  showed
	that for all multigraphs $G$,  $\Delta(G) \le \chi'(G) \le \Delta(G)+\mu(G)$, where $\mu(G)$, the \emph{maximum multiplicity} of $G$,  
	is the maximum number of  edges joining two vertices in $G$.  Restricting this result to graphs, it leads 
	to a natural classification of  them. Following Fiorini and Wilson~\cite{fw},   a  graph $G$ is of \emph{class 1} if $\chi'(G) = \Delta(G)$ and of \emph{class 2} if $\chi'(G) = \Delta(G)+1$.  Holyer~\cite{Holyer} showed that it is NP-complete to determine whether an arbitrary  graph is of class 1.  However, 
	by a conjecture of  Chetwynd and  Hilton~\cite{MR848854,MR975994} from  1986, there is a polynomial-time algorithm to 
	determine the chromatic index for  graphs $G$ with $\Delta(G) >\frac{1}{3} |V(G)|$. The conjecture involves  the 
	overfullness of graphs; in fact,    
	a number of longstanding conjectures listed in {\it Twenty Pretty Edge Coloring Conjectures} in~\cite{StiebSTF-Book} lie in deciding when a 
	graph is overfull.   The conjecture by Chetwynd and  Hilton~\cite{MR848854,MR975994} is stated as follows. 
	\begin{CON}[Overfull Conjecture]\label{overfull-con}
		Let $G$ be a  graph  with $\Delta(G)>\frac{1}{3}|V(G)|$. Then $\chi'(G)=\Delta(G)$  if and only if $G$ contains no $\Delta(G)$-overfull subgraph.  
	\end{CON}

	The  graph $P^*$, obtained from the Petersen graph by deleting one vertex, has $\chi'(P^*)=4$, 
	satisfies $\Delta(P^*)=\frac{1}{3}|V(P^*)|$ but contains no  3-overfull subgraph. 
	Thus the degree condition  $\Delta(G)>\frac{1}{3}|V(G)|$ in the conjecture above is best possible.  
	
	Applying Edmonds' matching polytope theorem, Seymour~\cite{seymour79}  showed  that whether a  graph  $G$ contains an overfull subgraph of maximum degree $\Delta(G)$ can be determined in polynomial time in $|V(G)|$.   Independently,  Niessen~\cite{MR1814514} in 2001 showed that 
	for graphs  $G$ with 
	$ \Delta(G) > \frac{1}{3}|V(G)|$, there are at most three induced  $\Delta(G)$-overfull subgraphs, and it is possible to find one  in polynomial time in $|V(G)|$. When 
	$\Delta(G) \geq \frac{1}{2}|V(G)|$, there is at most one induced $\Delta(G)$-overfull subgraph, and it is possible to find it in linear time in  $(|V(G)|+|E(G)|)$.
	Thus if the Overfull Conjecture is true, then the NP-complete problem of 
	determining the chromatic index becomes  polynomial-time solvable 
	for graphs $G$ with $\Delta(G)>\frac{1}{3}|V(G)|$. The Overfull Conjecture also
	implies
	several other longstanding conjectures in  edge colorings. 
	For example, the Just-overfull Conjecture~\cite[Conjecture 4.23]{StiebSTF-Book},
	the Vertex-splitting Conjecture~\cite[Conjecture 1]{MR1460574}, and the Independence Number Conjecture, 2-factor Conjecture
	and Average Degree Conjecture of Vizing~\cite{vizing-2factor,vizing-ind} when restricted to graphs $G$ 
	with $\Delta(G)>\frac{1}{3}|V(G)|$ (these three conjectures are also listed as  Conjectures 9-11 in ~\cite[Chapter 9]{StiebSTF-Book}).

	There have been some fairly strong results supporting the Overfull Conjecture when imposing a lower bound on the 
	minimum degree of $G$.  In case when $G$ is regular with even order, it is easy to verify that $G$ has no $\Delta(G)$-overfull subgraph if its vertex degrees are at least $\frac{1}{2}|V(G)|$.  Thus the well-known 1-Factorization Conjecture,   first stated in~\cite{MR772711} but may go back to Dirac in the early 1950s,  is a special case of the Overfull Conjecture.
	
	\begin{CON}[1-Factorization Conjecture]\label{con:1-factorization}
		Let  $G$ be a graph  of even order $n$. If $G$ is $k$-regular for some  $k\ge 2\lceil  \frac{1}{4}n\rceil-1$,  then  $G$ is 1-factorable; equivalently, $ \chi'(G) = \Delta(G)$.
	\end{CON}

	Hilton and Chetwynd~\cite{MR1001390} verified the 1-Factorization Conjecture if the vertex degree is at least  $0.823|V(G)|$.  Perkovi\'c and Reed~\cite{MR1439301} showed in 1997 that the 1-Factorization Conjecture is true for large regular graphs with vertex degree at least $ \frac{1}{2-\ve}|V(G)|$ for any given $0<\ve <1$. This was generalized by Vaughan~\cite{MR2993074} to multigraphs of bounded multiplicity.    In 2016,  Csaba, K\"uhn, Lo, Osthus and Treglown~\cite{MR3545109} verified the conjecture for sufficiently large $|V(G)|$. 
	Results were also obtained  for large even order  $G$
	when requiring the minimum degree to be  at least $\frac{2}{3}|V(G)|$ by Plantholt~\cite{MR4394718}. 
	Recently, Plantholt and the author~\cite{MR4563210}  showed that  for large even order  $G$, if its minimum degree is 
	arbitrarily  close to $\frac{1}{2}|V(G)|$ from above, then the Overfull Conjecture holds for $G$;
	this result has been extended to the same class of graphs $G$ with odd order by the 
	author~\cite{ 2205.08564}.

	Most results listed above, along with much of the existing research on the Overfull Conjecture, have
	focused on graphs $G$ that satisfy a minimum degree condition or with its subgraph induced by the maximum
	degree vertices restricted in some way. Consequently, the truth of the conjecture when these constraints are
	lifted remains largely unexplored. In this direction, the best result was due to Chetwynd and Hilton from
	1989~\cite{MR975994}, 
	which states that any  graph $G$ with $\Delta(G) \ge |V(G)| -3$ satisfies the Overfull Conjecture. After
	a significant research gap of over 30 years,
	this paper presents the following result.

	\begin{THM}\label{thm:1}
		For all $0<\ve  \le  \frac{1}{14}$, there exists $n_0 \in \mathbb{N}$
		for which the following statement holds:
		if $G$ is a graph on $n\ge n_0$ vertices with $\Delta(G) \ge (1-\ve )n$,   then $\chi'(G)=\Delta(G)$ if and only if $G$  contains no $\Delta(G)$-overfull 
		subgraph. 
	\end{THM}

	We say a graph $G$ is \emph{edge-chromatic critical} if  $\chi'(G)= \Delta(G)+1$ but any proper subgraph $H$ of $G$ satisfies $\chi'(H) \le \Delta(G)$. 
	The  Average Degree Conjecture of Vizing~\cite{vizing-ind}
	states that any $n$-vertex edge-chromatic critical graph  $G$ has average degree at least $\Delta(G)-1+\frac{3}{n}$. 
	The conjecture was completely confirmed  for graphs with maximum at most 6. 
	For more details and partial results on the conjecture, see~\cite{MR4596511}
	and the references therein.  
	As an application of Theorem~\ref{thm:1},  
	we show that the Average Degree Conjecture  holds for 
	large graphs $G$ with maximum degree at least $\frac{13}{14}|V(G)|$.
	
	\begin{THM}\label{thm2}
		There exists $n_0 \in \mathbb{N}$
		for which the following statement holds: if $G$ is an edge-chromatic critical graph with on $n\ge n_0$
		vertices with $\Delta(G)\ge \frac{13}{14}n$, then the average degree of $G$ is at least $\Delta(G)-1+\frac{3}{n}$. 
	\end{THM}
	
	\proof Let $n_0$ be the same integer as specified in Theorem~\ref{thm:1},
	and $G$ be an  edge-chromatic critical graph  on $n\ge n_0$
	vertices with $\Delta(G)\ge \frac{13}{14}n$. By Theorem~\ref{thm:1}, $G$ 
	contains a $\Delta(G)$-overfull subgraph $H$. Since $G$ is edge-chromatic critical
	and $\chi'(H)=\Delta(G)+1$, it follows that $H=G$. This, in particular, implies 
	$2|E(G)|\ge \Delta(G)(n-1)+2$ (in fact, $2|E(G)|= \Delta(G)(n-1)+2$ by $G$ being edge-chromatic critical). 
	Hence,  $2|E(G)|/n \ge \Delta(G)-1+\frac{3}{n}$, as 
	$\Delta(G) \le n-1$. 
	\qed

	We provide an overview of the proof  of Theorem~\ref{thm:1} before finishing this section.  The ``only if'' direction is trivial. To prove the ``if'' direction, we proceed  by a contradiction: suppose that $G$ contains no 
	$\Delta(G)$-overfull 
	subgraph but $\chi'(G)= \Delta(G)+1$.  By deleting vertices or edges if necessary, it is then natural to assume that $G$ is edge-chromatic critical.  Applying the main result 
	from~\cite{MR4480636} that  any	edge-chromatic critical $R$ with 
	$\Delta(R) -\frac{7\delta(R)}{4} \ge \frac{3|V(R)|-17}4$ is overfull, we get  $\delta(G)> \frac{1}{7}(1-4\ve) n$. 
	The rest proof is divided into two cases according to ``how far away'' that $G$ is from being $\Delta(G)$-regular.  To have a measure 
	for this property, for $v\in V(G)$, let $\df_G(v)=\Delta(G)-d_G(v)$
	be the \emph{deficiency} of $v$.   Then the deficiency of $G$,  denoted  \emph{$\df(G)$}, 
	is the sum of the deficiencies of the vertices of $G$.  If $\df(G)$ is at least  $\alpha n^2$
	for some constant  real $0<\alpha <1$, by adding new vertices and edges to $G$, we can  form a supergraph $G'$  of $G$ on some $n' >n$ vertices, avoiding the overfullness,   with $\Delta(G')=\Delta(G)$ 
	such that  either $\delta(G') \ge \frac{1}{2}(1+\ve')n'$ for some constant $0<\ve' \le \ve$ or $G'$ is an expander graph. In both cases, existing results imply that $G'$ is class 1 and so $G$
	is class 1. This gives a contradiction. 
	
	The most difficult case is when $\df(G)$ is small such as being only linear or sublinear in $n$.
	We will add edges to $G$ and also an extra vertex when $n$ is odd to create a  $\Delta(G)$-regular multigraph $G'$ such that $G'$
	still contains no $\Delta(G)$-overfull subgraph (there may be some vertex identifications in the process of forming $G'$). This step is very intricate, as the edges need to be added in a way such that no overfullness is created. A new lemma (Lemma~\ref{lem:graphical-biparite}) on degree sequence  is  developed 
	for this purpose, and several procedures (inducing identifying some vertices) are established to form $G'$.  
	The remaining  proof procedure relies on a mini-regularity type framework, typically applicable to graphs with even order and minimum degree greater than half the order (e.g., \cite{MR4563210, 2205.08564, MR2993074}). This work overcomes this limitation by extending the framework to $n$-vertex graphs where the minimum degree could be 
	lower than $\frac{1}{7}n$. 
	Another challenge in the proof arises from the potential high multiplicity of  the multigraph $G'$, 
	which could reach approximately  $\frac{1}{3} \Delta(G)$.  Precisely, 
	the edge decomposition procedure for $G'$  can be broken down into five steps:

	\begin{enumerate}[Step 1.]
		\item  Partition $V(G')$  into two subsets $A$ and $B$ such that  $|A|=|B|$,   and $||N_{G'}(v)\cap A|-|N_{G'}(v)\cap B||  \le c\eta n$ for any $v\in V(G')$, where $c>1$
		is some constant and $\eta>0$ is far less than $\ve$. We then refine the partition to ensure that (almost) all pairs of vertices in $G'$
		with a large number of edges joining them are placed in separate subsets, i.e., one vertex in $A$ and the other in $B$.

		\item  Let $\Delta'= \lceil \frac{\Delta(G)}{2}+5.3\eta n \rceil $ and 
		$
		k=\Delta'+ \lceil\sqrt{\Delta'} \rceil.
		$  
		Define edge sets $E_1 \subseteq  E(G'[A]), F_1\subseteq E(G'[B])$, and vertex-disjoint $F_{21}, F_{22}\subseteq E(G'[B])\setminus F_1$ 
		such that $|E_1|=|F_1|$ and $|F_{21}|=|F_{22}|$, where  edges of $E_1$ are incident with a fixed vertex in $G'$, $F_1$ is a set of multiple edges between two fixed vertices called $v_p$ and $v_{p+1}$, and all edges of $F_{21}$ are incident in $G'$ with $v_p$ and all edges of $F_{22}$ are incident in $G'$ with $v_{p+1}$.  We carefully define a multigraph $G_{A,B}$ as the union of $G'[A]-E_1$, $G'[B]-(F_1\cup F_{21}\cup F_{22})$, and a subset of edges between $A$
		and $B$. The challenge  lies in ensuring that $G_{A,B}$ does not contain any $k$-overfull  subgraph.  After the construction, we apply  a known edge coloring result to  equitably edge color $G_{A,B}$ using $k$ colors. 
		We will ``ignore''  the edges in $F_1$ for now, as they were all added to $G$ when forming $G'$. 
		Ignoring these edges also simplifies $G'$ (reduces the multiplicity of $G'[B]$), which is helpful for the 1-factor extension in Step 4.
		We will temporarily disregard the edges in $E_1$ until Step 5, where we will reintroduce them to form a nearly-bipartite multigraph and apply Theorem~\ref{thm:chromatic-index-nearly-bipartite-graph}. 
		The other reason for ignoring edges  in  $F_1$ during Steps 2-4 is to ensure that sets $A$ and $B$ have the same number of remaining edges.  Edges in $F_{21}\cup F_{22}$ will be removed in Steps 2-3 as to guarantee that 
		the degrees of $v_p$ and $v_{p+1}$ in the multigraph $G_{A,B}$  do not 
		exceed $\Delta'$. 
		These edges  of  $F_{21}\cup F_{22}$ will be colored  in Step 4.

		\item 
		As $G'$ has even order, we can pair up vertices  $u,v$ of $G'$ which are missing the same color, say $i$,  under the given edge $k$-coloring. We will use the maximum degree condition and the condition that $\df(G)$ is small 
		in finding an alternating path $P$ that starts at $u$, ends at $v$, and with edges alternating  between uncolored edges in-between $A$ and $B$ and edges with color $i$ within $A$ or $B$.  By swapping the color $i$
		and ``uncolor'' on edges of $P$, we increase the size of the color class $i$ by one. Repeating this procedure for all pairs of vertices missing the same color and all $k$ color classes will eventually extend each color class into a 1-factor of $G'$. 
		\item  In Step 3, some edges of $A$ and $B$ which were colored in Step 2 got  uncolored in Step 3.  We 
		will color these edges together with edges of $F_{21}\cup F_{22}$ using another $\ell$ colors and extend each of the new color class into a 1-factor. This step is divided into two stages to guarantee that 
		the multigraph  $R$ obtained at the end of Step 4  does not contain any $\Delta(R)$-overfull subgraph. 
		\item The remaining uncolored multigraph $R$   now have all uncolored edges between $A$ and $B$ and all edges of $E_1\cup F_1$, and   is  $(\Delta(G)-k-\ell)$-regular.  As all edges of $E_1$ are incident in $R$ with a fixed vertex, $R-F_1$
		is nearly-bipartite. We can show that $R-F_1$ contains no $\Delta(R-F_1)$ subgraph. 
		By a result of  Eggan and Plantholt~\cite{MR0830594},  $\chi'(R-F_1)=\Delta(R-F_1)$. 
		This implies $\chi'(G'-E_1) \le k+\ell+(\Delta(G')-k+\ell) =\Delta(G')=\Delta(G)$. As $G\subseteq G'-E_1$, 
		we get a contradiction to the assumption that $\chi'(G)=\Delta(G)+1$. 
	\end{enumerate}
	This work introduces several novel techniques for constructing the auxiliary multigraph $G'$,   defining the multigraph $G_{A,B}$ in Step 2, extending 1-factors in Step  4, and using a nearly-bipartite structure in Step 5.  These innovations are crucial to address the significantly weaker minimum degree lower bound of  $\frac{1}{7}(1-4\ve) n$ 
	compared to the previously used framework (e.g., \cite{MR4563210, 2205.08564, MR2993074}) which relied on a minimum degree of at least $\frac{1}{2}(1+\ve)n$.  We anticipate that these novel techniques will serve as stepping stones for further progress on the Overfull Conjecture.

	We are not able to deduce an algorithm that finds an optimal edge coloring for $G$ because  the proof is  by contradiction and 
	the properties of $G$, when assumed to be edge-chromatic critical,  are  used in the proof. 
	The remainder of this paper is organized as follows.
	In the next section,  we 
	introduce some notation and   preliminary results. 
	In   Section 3, we prove 
	Theorem~\ref{thm:1}.

\section{Notation and preliminaries}

Let $G$ be a multigraph. 
We use $V(G)$ and  $E(G)$ to denote the vertex set and the edge set of $G$,
respectively, and  let  $e(G)=|E(G)|$. 
For $v\in V(G)$, $N_G(v)$ is the set of neighbors of $v$ 
in $G$, and 
$d_G(v)$, the degree of $v$
in $G$, is the number of edges of $G$ that are incident with $v$.
We let $d_G^s(v)=|N_G(v)|$, and call it the \emph{simple degree}
of $v$ in $G$.   We let $G^s$ be the underlying simple graph of $G$. 
Thus we have 
 $d_{G^s}(v)=d_G^s(v)$.    
We let $\delta^s(G)$ be the minimum simple degree among all the simple degrees of 
vertices of $G$.  
Let $V_1,
V_2\subseteq V(G)$ be two disjoint  vertex sets. Then $E_G(V_1,V_2)$ is the set
of edges in $G$  with one end in $V_1$ and the other end in $V_2$, and  $e_G(V_1,V_2):=|E_G(V_1,V_2)|$.  We write $E_G(v,V_2)$ and $e_G(v,V_2)$
if $V_1=\{v\}$ is a singleton.   
We also use $G[V_1,V_2]$
to denote the bipartite subgraph of $G$ with vertex set $V_1\cup V_2$
and edge set $E_G(V_1, V_2)$.  
Let 
$S\subseteq V(G)$ and $v\in V(G)$. Then  $N_G(v, S):=N_G(v)\cap S$ and $d_G(v, S):=e_G(v,S\setminus\{v\})$,  and $d_G^s(v, S):=|N_G(v, S)|$.  
The subgraph of $G$ induced by  $S$ is  $G[S]$, and  $G-S:=G[V(G)\setminus S]$. 
If $F\subseteq E(G)$, then $G-F$ is obtained from $G$ by deleting all
the edges of $F$. 
Let $\mu(G)=\max\{e_G(u,v)\,:\, u,v\in V(G)\}$ be 
the multiplicity of $G$, and  for $v\in V(G)$, let $\mu_G(v)=\max\{e_G(u,v): u\in N_G(v)\}$
be the multiplicity of $v$ in $G$.

We will use the following notation: $0<a \ll b \le 1$. 
Precisely, if we say a claim is true provided that $0<a \ll b \le 1$, 
then this means that there exists a non-decreasing function $f:(0,1]\rightarrow (0,1]$ such that the statement holds for all $0<a,b\le 1$ satisfying $a \le f(b)$. 

We present in the rest of this section some preliminaries that will be needed in the proof of Theorem~\ref{thm:1}.  

\subsection{Results on degree sequences}

The results in this subsection will be used to construct  a regular supergraph  $G'$ based on a given graph $G$ by adding edges.

\begin{LEM}\label{lem:graphical-regular}
	Let  $m\ge 3$ and  $d\ge 2$ be  positive integers with $d$ even and  $d<m$. Then there is a  $d$-regular   graph  $L$
	on $m$ vertices. 
	Furthermore $L$ can be constructed in  $O(dm)$ steps. 
\end{LEM}

\pf  Let $\{v_1, v_2, \ldots, v_m\}$ be a set of  $m$ vertices.   The indices  of the vertices will be taken modular $m$ in the arguments below. 
For each  $i\in [1,d/2]$,   we add an  edge  joining the vertices $v_j$ and $v_{j+i}$ for each  $j\in [1,m]$ (adding a cycle containing all the vertices for each $i$).   Denote the resulting multigraph by $L$.  Note that $$N_L(v_j)=\{v_{j+1}, v_{j+2}, \ldots, v_{j+d/2},v_{j-1}, v_{j-2}, \ldots, v_{j-d/2}\}$$ for each  $j\in [1,m]$. 
Since $m>d$, we have $d_L(v_j)=|N_L(v_j)|=d$.    Thus $L$ is a desired $d$-regular graph. 
By the process,  $L$ can be constructed  in $O(dm)$ steps. 
\qed

Let  $m\ge 2$ be an integer. A sequence of nonnegative and non-increasing integers $(d_1, \ldots, d_m)$ is  \emph{graphic}
if there is a  graph $L$ such that the degree sequence of $L$ is $(d_1, \ldots, d_m)$. In this case, 
we say that $L$ \emph{realizes} $(d_1, \ldots, d_m)$. The Havel–Hakimi Algorithm 
determines efficiently if a given nonnegative integer sequence is graphic. 

\begin{LEM}[Havel–Hakimi Algorithm, \cite{MR0089165,MR148049}]\label{lem:Havel–Hakimi Algorithm}
	Let $s\ge 1$ and $m\ge 1$ be positive integers with $m\ge s$. 
	A nonnegative  and non-increasing integer sequence $(s, d_1, \ldots, d_m)$ 
	is graphic if and only if $(d_1-1,  \ldots, d_s-1,  d_{s+1}, \ldots, d_m)$ (after reordering the entries so the sequence is  non-increasing) is graphic. Furthermore, 
	if $(s, d_1, \ldots, d_m)$  is graphic, then  the sequence can be realized in polynomial time in $m$. 
\end{LEM}

Let  $m\ge 2$ be an integer. A sequence of nonnegative and non-increasing integers $(d_1, \ldots, d_m)$ is 
\emph{admissible}  if $\sum_{i=1}^m d_i$ is even and $d_1\le \sum_{i=2}^m d_i$. 

\begin{LEM}\label{lem:degree-sequence-d-d-1}
	Let $m$ and $d$ be positive integers such that $m\ge d+1 \ge 3$.  
	Suppose $(d_1, \ldots, d_m)$ is  a  sequence of positive integers with $d_1= \ldots =d_{t}=d$
	and $d_{t+1} =\ldots = d_m=d-1$ and $\sum_{i=1}^m d_i$ even, where $t\in [1,m]$. Then $(d_1, \ldots, d_m)$  
	is graphic.  Furthermore,  the sequence can be realized in polynomial time in $m$.  
\end{LEM}

\pf  Applying the Havel–Hakimi Algorithm repeatedly,  the sequence  $(d,\ldots, d, d-1, \ldots, d-1)$
can be reduced to a sequence of the form $(2, \ldots, 2,  1, \ldots, 1)$ or $(1, \ldots, 1)$. 
Since  $\sum_{i=1}^m d_i$  is even  and we reduce the sum of the sequence by an even amount every time when applying  the Havel–Hakimi Algorithm, 
the sum of the  integers in the reduced sequence is even. As  the sum of the  integers in the reduced sequence is even, if the sequence  contains $2$, then we know that it contains another $2$ or at least two  $1$s; 
if the sequence does not contain $2$, then it contains an even number of $1$s.  Thus the reduced sequence is 
admissible. It is clear that any  admissible sequence of the form 
$(2, \ldots, 2,  1, \ldots, 1)$ or $(1, \ldots, 1)$ is graphic. 
Therefore, $(d_1, \ldots, d_m)$  is graphic by the Havel–Hakimi Algorithm.  The sequence can be realized in polynomial time in $m$ also by Lemma~\ref{lem:Havel–Hakimi Algorithm}. 
\qed

\begin{LEM}\label{lem:graphical-biparite}
	Let  $m\ge 2$ and $(d_1, d_2, \ldots, d_m)$ be an  admissible sequence.  Then there is a bipartite  multigraph $L$ on  $\{v_1, \ldots, v_m\}$ 
	and an even index   $p\in [2, m]$ satisfying the following properties, where $d_i$ is taken to be $0$ if  $i\not\in [1,m]$. 
	\begin{enumerate}[(a)]
		\item $d_L(v_i)=d_i-d_{i+1}$ and  $d_L(v_{i+1})=0$ for any odd  $i\in [1,p]$.
		\item $d_L(v_{p+1})  \le d_{p+1} $ and   $d_L(v_i) =d_i$ for any  $i\in [p+2,m]$.
		\item $\{v_1, \ldots, v_{p}\}$ and $\{v_{p+1}, \ldots, v_m\}$ form a bipartition of $L$. 
		\item the underlying simple graph $L^s$ of $L$ is a forest such that any vertex has at most two non-leaf neighbors in $L^s$.
		\item  Let distinct $i, j\in [1, \frac{p}{2}]$ such that both $N_L(v_{2i-1}) \ne \emptyset$ and $N_L(v_{2j-1}) \ne \emptyset$. 
		Let $v_{p_1}, \ldots, v_{p_s}$ be all 
		neighbors of $v_{2i-1}$ in $L$ and $v_{q_1}, \ldots, v_{q_t}$ be all 
		neighbors of $v_{2j-1}$ in $L$, where $s, t\ge 1$ are integers  and $p_1, \ldots, p_s$  and 
		$q_1, \ldots, q_t$ are integers satisfying  $p_1, \ldots, p_s, q_1, \ldots, q_t\in [p+1, m]$, $p_1<  \ldots <p_s$, and $q_1<\ldots <q_t$. 
		Then  $p_1, \ldots, p_s$ are consecutive integers;  and if $i<j$, then  $p_1 \ge q_t$. 
\item  Let distinct $i, j\in [p+1, m]$ such that both $N_L(v_i) \ne \emptyset$ and $N_L(v_j) \ne \emptyset$. 
Let $v_{i_1}, \ldots, v_{i_s}$ be all 
neighbors of $v_{i}$ in $L$ and $v_{j_1}, \ldots, v_{j_t}$ be all 
neighbors of $v_{j}$ in $L$, where $s, t\ge 1$ are integers  and $i_1, \ldots, i_s$  and 
$j_1, \ldots, j_t$ are integers satisfying  $i_1, \ldots, i_s, j_1, \ldots, j_t\in \{ 2i-1: \in [1, \frac{p}{2}]\}$, $i_1<  \ldots <i_s$, and $j_1<\ldots <j_t$. 
If $i<j$, then  $i_1 \ge j_t$. 
	\end{enumerate}
	Furthermore, there is a polynomial-time algorithm
	that finds  $L$ and $p$.  See Figure~\ref{f3} for a construction of $L$ with $m=14$ and $p=8$. 
\end{LEM}

\pf  
We apply induction on $m$. 
If $m=2$, then $d_1  \le d_2$ follows from  the condition that $d_1\le \sum_{i=2}^m d_i$. Thus $d_1=d_2$. 
Letting  $L=\overline{K_2}$  and $p=2$ gives a desired graph and index. 

Thus we assume $m\ge 3$ and 
start by letting $L_0$ be the empty graph on the vertex set  $\{v_1, \ldots, v_m\}$. 
Let $i_0\in [1,m]$ be the smallest odd index such that $d_{i_0}>d_{i_0+1}$.  If $i_0$ does not exist, then we let $L=L_0$
and $p=m$. 
If $i_0=m$, then we simply let $L=L_0$ and $p=m-1$.

Therefore we assume  $i_0<m$.  If $d_{i_0}-d_{i_0+1}  \ge d_m$, then 
we add $d_m$ edges joining $v_{i_0}$ and $v_m$ in $L_0$ and still call the resulting multigraph $L_0$. 
Let $f_{i_0}=d_{i_0}-d_m$ and $f_j=d_j$ for each $j\in [1,m-1]$ with $j\ne i_0$. 
Note that $f_{i_0}\ge f_{i_0+1}$ by $d_{i_0}-d_{i_0+1}  \ge d_m$.  Thus  $(f_1, \ldots, f_{m-1})$ 
is still a non-increasing and non-negative integer sequence. 
As $\sum_{i=1}^m d_i$ is even and $\sum_{i=1}^{m-1} f_i=\sum_{i=1}^m d_i-2d_m$, we know that 
$\sum_{i=1}^{m-1} f_i$ is even.  Furthermore,  if $i_0\ge 3$,  then we have $f_1=f_2$ and so 
$f_1\le \sum_{i=2}^{m-1}f_i$; if $i_0=1$, then $f_1=d_1-d_m \le \sum_{i=2}^m d_i-d_m=\sum_{i=2}^{m-1}f_i$. 
Thus $(f_1, \ldots, f_{m-1})$ is admissible. 
Applying induction,  there is a bipartite multigraph $J$ on $\{v_1, \ldots, v_{m-1}\}$ and an even index $q\in [1,m-1]$ 
satisfying all the properties (a)-(f).  
Let $L=L_0\cup J$  and $p=q$.   As  $d_L(v_i)=d_{J}(v_i)$ for any $i\in [1, m-1]$ with $i\ne i_0$,  $d_L(v_{i_0})=d_{J}(v_{i_0})+d_m$,  $d_L(v_m)=d_m$,  and all the edges of $L_0$ are joining $v_{i_0}$ and $v_m$, it follows that $L$ and $p$ satisfies property (a)-(d). 
Note that $d_m>0$ implies  $d_{m-1}>0$. Thus if $e_{L_0}(v_{i_0}, v_m)>0$ and   if $v_{i_0}$ has a neighbor in $J$, 
then $v_{m-1}$ is a neighbor of $v_{i_0}$ in $J$ by (e) and (f).   Again since all the edges of $L_0$ are joining $v_{i_0}$ and $v_m$, 
it follows that $L$ also satisfies properties (e) and (f).

If $d_{i_0}-d_{i_0+1}  <d_m$, we add $d_{i_0}-d_{i_0+1}$ edges joining $v_{i_0}$ and $v_m$ in $L_0$ and still call the resulting multigraph $L_0$. 
Let $f_{i_0}=d_{i_0+1}$, $f_m=d_m-(d_{i_0}-d_{i_0+1})$ and $f_j=d_j$ for each $j\in [1,m-1]$ with $j\ne i_0$. 
Then we let  $j_0\in [i_0+2, m]$ be the smallest  odd index such that  $f_{j_0}>f_{j_0+1}$. Note that $j_0$ exists as we have $f_{m-1}=d_{m-1} \ge d_m>d_m-(d_{i_0}-d_{i_0+1})=f_m$. 
If $j_0=m$,   we simply let $L=L_0$ and $p=m-1$,  which are desired bipartite multigraph and index for $(d_1, \ldots, d_m)$, respectively. 
Thus we assume $j_0<m$.  

If $f_{j_0}-f_{{j_0}+1}  \ge f_m$, then 
we add $f_m$ edges joining $v_{j_0}$ and $v_m$ in $L_0$ and still call the resulting multigraph $L_0$. 
Let $g_{j_0}=f_{j_0}-f_m$ and $g_j=f_j$ for each $j\in [1,m-1]$ with $j\ne j_0$. 
Note that $g_{j_0}\ge g_{j_0+1}$ by $f_{j_0}-f_{j_0+1}  \ge f_m$.   It is routine to check that $(g_1, \ldots, g_{m-1})$ 
is an admissible sequence. 
Applying induction,  there is a bipartite multigraph $J$ on $\{v_1, \ldots, v_{m-1}\}$ and an even index $q\in [1,m-1]$ 
satisfying properties (a)-(f). If  $e_{L_0}(v_{j_0}, v_m)>0$ and the vertex $v_{j_0}$ 
has a neighbor in $J$, then $v_{m-1}$ is  a neighbor of $v_{j_0}$ in $J$. Also, in $J$, by property (a), we know that 
vertices from $\{v_1, \ldots, v_{j_0-1}\}$ are isolated vertices in $J$. 
Thus  $L:=L_0\cup J$ and $p:=q$ are respectively  desired bipartite multigraph and index for the sequence $(d_1,\ldots, d_m)$. 

Thus we assume  $f_{j_0}-f_{j_0+1}  <f_m$.  We add $f_{j_0}-f_{j_0+1}$ edges joining $v_{j_0}$ and $v_m$ in $L_0$ and still call the resulting multigraph $L_0$. 
Let $g_{j_0}=g_{j_0+1}$, $g_m=f_m-(f_{j_0}-f_{j_0+1})$ and $g_j=f_j$ for each $j\in [1,m-1]$ with $j\ne j_0$. 
We now consider the sequence $(g_{j_0}, g_{j_0+1}, \ldots, g_m)$.  As $g_i=g_{i+1}$ for each odd index $i\in [1,j_0-2]$, it follows that 
$\sum_{i=j_0}^mg_i$ is even. Furthermore,  $g_{j_0} \le \sum_{i=j_0+1}^m g_i$ by $g_{j_0}=g_{j_0+1}$. Thus by applying the induction hypothesis on 
$(g_{j_0}, g_{j_0+1}, \ldots, g_m)$, we find a bipartite multigraph $J$ on $\{v_{j_0}, v_{j_0+1}, \ldots, v_m\}$  and an index $q\in [1,m-j_0+1]$ (translating the indices $j_0, \ldots, m$ to $1, \ldots, m-j_0+1$) satisfying the 
properties (a)-(f).  As all edges are incident with the vertex $v_m$ in $L_0$, 
letting $L=L_0\cup J$ and $p=q+j_0-1$ gives  respectively  a desired bipartite multigraph and index for the sequence $(d_1, \ldots, d_m)$. 

The inductive process above indicates that $L$ and $p$ can be found in  polynomial time in $m$. 
\qed 

\begin{figure}[!htb]
	\begin{center}

		\begin{tikzpicture}[line cap=round,line join=round,>=triangle 45,x=1cm,y=1cm]
			\clip(-2.422441573688839,-6.963959497952107) rectangle (11.397150945208852,4.3406227501517956);
			\draw (-2.1890901052646754,-5.045291868686766-0.5) node[anchor=north west] {$d_2=8$};
			\draw (0.2222084017850152,-3.33404776690956960-0.5) node[anchor=north west] {$d_3=8$};
			\draw (2.400155440410542,-3.3340477669095696-0.5) node[anchor=north west] {$d_5=7$};
			\draw (-0.11485483038322115,-5.045291868686766-0.5) node[anchor=north west] {$d_4=7$};
			\draw (2.218659853858415,-5.097147750558801-0.5) node[anchor=north west] {$d_6=5$};
			\draw (4.318823069675888,-5.071219809622784-0.5) node[anchor=north west] {$d_8=4$};
			\draw (1.3111819210977786,3.381288935519125) node[anchor=north west] {$d_{10}=3$};
			\draw (3.4372730778512697,3.381288935519125) node[anchor=north west] {$d_{11}=3$};
			\draw (5.433724529924669,3.3553609945831067) node[anchor=north west] {$d_{12}=2$};
			\draw (7.352392159190015,3.381288935519125) node[anchor=north west] {$d_{13}=1$};
			\draw [line width=2pt] (8,2)-- (-1.4104275186587365,-3.489361334775554);
			\draw [line width=2pt] (6,2)-- (-1.4104275186587365,-3.489361334775554);
			\draw [line width=2pt] (6,2)-- (-1.4104275186587365,-3.489361334775554);
			\draw [line width=2pt] (4,2)-- (-1.4104275186587365,-3.489361334775554);
			\draw (-1.9038827549684756,-3.3599757078455874-0.5) node[anchor=north west] {$d_1=14$};
			\draw [line width=2pt] (2,2)-- (0.593630891292589,-3.489361334775554);
			\draw [line width=2pt] (2,2)-- (2.989286921579231,-3.4663261806381827);
			\draw (4.500318656228015,-3.3340477669095696-0.5) node[anchor=north west] {$d_7=5$};
			\draw (-0.6593415900396029,3.2775771717750524) node[anchor=north west] {$d_9=3$};
			\draw [line width=2pt] (0,2)-- (5.016380485667928,-3.4663261806381827);
			\draw (9.530339197815541,3.3553609945831067) node[anchor=north west] {$d_{14}=0$};
			\draw [shift={(-5.057720560650102,-2.0999777856952764)},line width=2pt]  plot[domain=-0.16819172495497625:0.5262726512094615,variable=\t]({1*8.162183369382111*cos(\t r)+0*8.162183369382111*sin(\t r)},{0*8.162183369382111*cos(\t r)+1*8.162183369382111*sin(\t r)});
			\draw [shift={(16.11352267373852,-15.350332308736593)},line width=2pt]  plot[domain=2.1802946421631786:2.54657763461773,variable=\t]({1*21.160611115720783*cos(\t r)+0*21.160611115720783*sin(\t r)},{0*21.160611115720783*cos(\t r)+1*21.160611115720783*sin(\t r)});
			\draw [shift={(-12.053630592908812,12.411793713229521)},line width=2pt]  plot[domain=5.302236767578541:5.707820816381954,variable=\t]({1*19.134380145184306*cos(\t r)+0*19.134380145184306*sin(\t r)},{0*19.134380145184306*cos(\t r)+1*19.134380145184306*sin(\t r)});
			\draw [shift={(-25.681946034987465,37.022836188225)},line width=2pt]  plot[domain=5.2521584593785455:5.447744242486777,variable=\t]({1*47.22652601273039*cos(\t r)+0*47.22652601273039*sin(\t r)},{0*47.22652601273039*cos(\t r)+1*47.22652601273039*sin(\t r)});
			\begin{scriptsize}
				\draw [fill=black] (-1.4334626727961082,-5.101822124391561) circle (4.5pt);
				\draw [fill=black] (-1.4104275186587365,-3.489361334775554) circle (4.5pt);
				\draw [fill=black] (0.593630891292589,-3.489361334775554) circle (4.5pt);
				\draw [fill=black] (5.016380485667928,-3.4663261806381827) circle (4.5pt);
				\draw [fill=black] (0.5705957371552174,-5.07878697025419) circle (4.5pt);
				\draw [fill=black] (4,2) circle (4.5pt);
				\draw [fill=black] (6,2) circle (4.5pt);
				\draw [fill=black] (8,2) circle (4.5pt);
				\draw [fill=black] (2,2) circle (4.5pt);
				\draw [fill=black] (2.989286921579231,-3.4663261806381827) circle (4.5pt);
				\draw [fill=black] (0,2) circle (4.5pt);
				\draw [fill=black] (2.966251767441859,-5.101822124391561) circle (4.5pt);
				\draw [fill=black] (4.970310177393185,-5.0557518161168185) circle (4.5pt);
				\draw [fill=black] (10,2) circle (4.5pt);
			\end{scriptsize}
		\end{tikzpicture}
	\end{center}
	\vspace{-1cm}
	\label{f3}
	\caption{The construction of a bipartite multigraph $L$ described in  Lemma~\ref{lem:graphical-biparite}, where $p=8$.}
\end{figure}

\subsection{Results on edge colorings}

Let  $k\ge 0$ be an integer and $\CC^k(G)$ be the set of all edge $k$-colorings of a multigraph $G$, and $\varphi\in \CC^k(G)$. 
For any $v\in V(G)$, the set of colors \emph{present} at $v$ is 
$\varphi(v)=\{\varphi(e)\,:\, \text{$e \in E(G)$ is incident to $v$}\}$, and the set of colors \emph{missing} at $v$ is $\pbar(v)=[1,k]\setminus\varphi(v)$.  For a subset $X$ of $V(G)$ and a color $i\in [1,k]$, define 
$\pbar_X^{-1}(i)= \{v\in X: i\in \pbar(v)\}$. We simply write $\pbar^{-1}(i)$ for $\pbar^{-1}_{V(G)}(i)$.

Recall that a graph $G$ is \emph{edge-chromatic critical} if $\chi'(G)= \Delta(G)+1$ but any proper subgraph $H$ of $G$ satisfies $\chi'(H) \le \Delta(G)$.  
An edge $e\in E(G)$  is \emph{critical} if $\chi'(G)= \Delta(G)+1$  but $\chi'(G-e) \le \Delta(G)$. It is clear that a graph is edge-chromatic critical 
if and only if it is connected and all of its edges are critical. 

\begin{LEM}[Vizing's Adjacency Lemma (VAL)~\cite{Vizing-2-classes}] \label{lem:val}Let $G$ be a class 2 graph with maximum degree $\Delta$. If $e=xy$ is a critical edge of $G$, then $x$ has at least $\Delta-d_G(y)+1$ neighbors  of degree $\Delta$ from  $V(G)\setminus \{y\}$. As a consequence, $d_G(x)+d_G(y)\ge \Delta+2$. 
	\label{thm:val}
\end{LEM}

In 1960's, Gupta~\cite{Gupta-67} and, independently, Vizing~\cite{Vizing-2-classes}   provided an upper bound on the chromatic index
of multigraphs, and K\"onig~\cite{MR1511872} obtained  an exact value 
of the chromatic index for bipartite multigraphs. 

\begin{THM}[\cite{Gupta-67, Vizing-2-classes}]\label{thm:chromatic-index}
Every multigraph  $G$ satisfies $\chi'(G) \le \Delta(G)+\mu(G)$. 
\end{THM}

\begin{THM}[\cite{MR1511872}]\label{konig}
	Every bipartite multigraph $G$ satisfies $\chi'(G)=\Delta(G)$. 
\end{THM}

A multigraph $G$  is \emph{nearly-bipartite} if deleting at most one vertex in $G$ results in a bipartite multigraph.  The following result was due to 
Eggan and Plantholt. 
\begin{THM}[{\cite[Theorem 2]{MR0830594}}]\label{thm:chromatic-index-nearly-bipartite-graph}
	Let $G$ be a nearly-bipartite multigraph. Then $\chi'(G)=\Delta(G)$ if and only if $G$ contains no $\Delta(G)$-overfull subgraph. 
\end{THM}

For any multigraph $G$ on at least three vertices, let $\rho(G)$,   the \emph{density} of $G$, be defined as 
$$\rho(G) =\max \left \{ \frac{e(H)}{\lfloor  \frac{1}{2}|V(H)|\rfloor }: H\subseteq G,  \quad  3\le  |V(H)| \equiv 1 \pmod{2} \right\}.$$
As $\chi'(G) \ge \chi'(H)$ for any subgraph $H$ of $G$, and $\chi'(H) \ge \frac{e(H)}{\lfloor  \frac{1}{2}|V(H)|\rfloor }$, it follows that 
$\chi'(G) \ge \rho(G)$.   We will use the following upper bound  on $\chi'(G)$ involving $\Delta(G)$ and $\rho(G)$~\cite[Corollary 5.20 and Theorem 5.24]{StiebSTF-Book}.

\begin{THM}\label{thm:chromatic-index-bound}
Every multigraph $G$ satisfies $\chi'(G) \le \max\{\Delta(G)+\sqrt{(\Delta(G)-1)/2}, \rho(G)\}$. Furthermore, 
an edge $\max\{\Delta(G)+\sqrt{(\Delta(G)-1)/2}, \rho(G)\}$-coloring can be realized in polynomial time in $|V(G)|$. 
\end{THM}

An edge $k$-coloring of a multigraph $G$ is said to be \emph{equalized} if each color
class contains either $\lfloor |E(G)|/k \rfloor$ or $\lceil |E(G)|/k \rceil$ edges.  McDiarmid~\cite{MR300623} observed the following result. 
\begin{THM}\label{lem:equa-edge-coloring}
	Let $G$ be a multigraph with chromatic index $\chi'(G)$. Then for all $k\ge \chi'(G)$, there is an equalized edge-coloring of $G$ with $k$ colors. 
\end{THM}

Given an edge coloring of $G$ and a given color $i$, 
since  the color class $i$ 
is a matching, we have the Parity Lemma below, see~\cite[Lemma 2.1]{MR2028248}.

\begin{LEM}[Parity Lemma]
	Let $G$ be a multigraph and $\varphi\in \CC^k(G)$ for some integer $k\ge \Delta(G)$. 
	Then 
	$|\pbar^{-1}(i)| \equiv |V(G)| \pmod{2}$ for every color $i\in [1,k]$. 
\end{LEM}

\subsection{Results  on the existence of overfull subgraphs and matchings}

\begin{LEM}[{\cite[Lemma 2.9]{MR4563210}}]\label{lem:overfull2}
	Let  $G$  be a  graph of  even order  $n$ with  minimum degree greater than $\frac{n}{2}$.  Then $G$
	contains no $\Delta(G)$-overfull subgraph  provided that $G$ has at least two vertices of minimum degree. 
\end{LEM}

\begin{LEM}\label{lem:matching-in-bipartite}
	Let $G[X,Y]$ be bipartite multigraph with $|X|=|Y|=n$. Suppose that $\delta^s(G) \ge t$ 
	 where 
	$t\in [1,n-1]$,  
	and except at most $t$ vertices all other vertices of $G$ have simple degree 
	at least $\frac{n}{2}$ in $G$. Then 
	$G$ has a perfect matching. 
\end{LEM}
\pf  
We show that $G$ satisfies Hall's Condition. If not, 
we let $S\subseteq X$ with smallest cardinality such that 
$|S|>|N_{G}(S)|$.  Since $\delta(G)\ge 1$, it follows  that $|S|\ge 1$. 
Then by the choice of $S$,  we have $|S|=|N_{G}(S)|+1$ and $|N_{G}(S)|<|Y|$.  
As  $|S|>|N_{G}(S)|$, it follows that 
$|S| \ge \delta^s(G)+1 \ge t+1$. As $G$
has at most $t$ vertices of  simple degree less than $\frac{n}{2}$, $S$ contains a vertex 
of simple degree at least $\frac{n}{2}$ in $G$ and so $|N_{G}(S)|\ge \frac{n}{2}$. 
Thus $|S|>\frac{n}{2}$ and consequently  
$|X\setminus S|  <\frac{n}{2}$. Since $|N_{G}(S)|<|Y|$, there exists $z\in Y\setminus N_{G}(S)$ such that $N_{G}(z)\subseteq X\setminus S$. As $\delta^s(G) \ge t$, we have $|X\setminus S| \ge t$. 
As $|Y\setminus N_{G}(S)|=|Y|-|S|+1 =|X|-|S|+1\ge t+1$ and $G$
has at most $t$ vertices of simple degree less than $\frac{n}{2}$, $Y\setminus N_{G}(S)$
contains a vertex of simple degree at least $\frac{n}{2}$ in $G$.  As a result, $|X\setminus S|\ge \frac{n}{2}$,
contradicting the earlier assertion that 
$|X\setminus S| <\frac{n}{2}$.
Hence $G$ has a perfect matching.
\qed

\section{Proof of Theorem~\ref{thm:1}}

We will need the following results. 

\begin{THM}[\cite{MR4480636}]\label{thm:overfull-min-degree}
	Let $G$ be an  edge-chromatic critical graph of order $n$.  If 
	$\Delta(G) -\frac{7\delta(G)}{4} \ge \frac{3n-17}4$,  then $G$ is overfull. 
\end{THM}

\begin{THM}[\cite{MR4563210}]\label{thm:PS}
	For all $0<\ve <1$, there exists $n_1 \in \mathbb{N}$
	such that the following statement holds:
	if $G$ is a graph on $2n\ge n_1$ vertices with $\delta(G) \ge (1+\ve)n$,  then $\chi'(G)=\Delta(G)$ if and only if $G$ contains no $\Delta(G)$-overfull subgraph. Furthermore, there is a polynomial-time algorithm that finds an optimal coloring. 
\end{THM}

\begin{THM}[\cite{2205.08564}]\label{thm:S}
	For all $0<\ve <1$, there exists $n_2  \in \mathbb{N}$
	such that the following statement holds:
	if $G$ is a graph on $2n-1\ge n_2$ vertices with $\delta(G) \ge (1+\ve)n$,  then $\chi'(G)=\Delta(G)$ if and only if $G$  is not overfull. Furthermore, there is a polynomial-time algorithm that finds an optimal coloring. 
\end{THM}

\begin{THM}[\cite{MR4072961}]\label{thm:H-cycle}
	Let $G$ be an  edge-chromatic critical graph of order $n$.  If $\Delta(G)\ge \frac{2n}{3}+12$, then $G$ is hamiltonian. 
\end{THM}

Given $0<\nu \le \tau <1$,  we say that  a graph $G$  on $n$ vertices is a \emph{robust $(\nu,\tau)$-expander},  if for all $S\subseteq V(G)$ with 
$  \tau n\le  |S| \le (1-\tau)n $ the number of vertices that have at least $\nu n$ neighbors in $S$ is at least $|S|+\nu n$. 
The  \emph{$\nu$-robust neighbourhood} $RN_G(S)$
is the set of all those
vertices of $G$ which have at least  $\nu n$ neighbours in $S$.

\begin{THM}[\cite{MR3028574}]\label{thm:robust-expander}
	For every $\alpha >0$ there exists $\tau >0$ such that for every $\nu >0$ 
	there exists $n_4(\alpha,\nu,\tau) \in \mathbb{N}$ for which the following holds.  Suppose that 
	\begin{enumerate}[(i)]
		\item  $G$ is an $r$-regular graph on $n\ge n_4$ vertices, where $r \ge \alpha n$ is even; 
		\item $G$ is a robust $(\nu,\tau)$-expander. 
	\end{enumerate}
	Then $G$ has a Hamilton decomposition. Moreover, this decomposition can be found in time polynomial in $n$. 
\end{THM}

We will also need the following result, which was proved using 
Chernoff bound. 

\begin{LEM}[\cite{SHAN2022429}, Lemma 3.2]\label{lem:partition}
	There exists  $n_3  \in \mathbb{N}$ such that for all $n\ge n_3$ the
	following holds. Let $G$ be a graph on $2n$ vertices, and $N=\{x_1,y_1,\ldots, x_t,y_t\}\subseteq V(G)$, where $t\in [1,n]$. 
	Then $V(G)$ can be partitioned into two  parts 
	$A$ and $B$ satisfying the properties below:
	\begin{enumerate}[(i)]
		\item  $|A|=|B|$;
		\item $|A\cap \{x_i,y_i\}|=1$ for each $i\in [1,t]$;
		\item $| d_G(v, A)-d_G(v, B)| \le n^{\frac{2}{3}}$ for each $v\in V(G)$.
	\end{enumerate}
	Furthermore, one such partition can be constructed  in time polynomial in $n$.
\end{LEM}

We are now ready to prove the main theorem. 

\begin{THM1}
For all $0<\ve  \le  \frac{1}{14}$, there exists $n_0 \in \mathbb{N}$
for which the following statement holds:
if $G$ is a  graph on $n\ge n_0$ vertices with $\Delta(G) \ge (1-\ve)n$,   then $\chi'(G)=\Delta(G)$ if and only if $G$  contains no $\Delta(G)$-overfull 
subgraph. 
\end{THM1}

\pf  Throughout the proof,  we let $\Delta=\Delta(G)$, $\delta=\delta(G)$,  $V_\Delta$ and $V_\delta$ be the set of 
maximum degree and minimum degree vertices of $G$, respectively.

Let $\tau^*=\tau(\frac{3}{7})$  and $n^*_4=n_4(\frac{3}{7},  \frac{1}{2}\eta^2, \tau^*)$ be as defined in Theorem~\ref{thm:robust-expander}.  
Choose $n_0 \in \mathbb{N}$ and $\eta$ such that 
\begin{equation}\label{eqn:parameters}
	0< \frac{2}{n_0 } \le  \min \left \{\frac{1}{n^*_1}, \frac{1}{n_2}, \frac{1}{n_3},   \frac{1}{n^*_4} \right \}  \ll  \eta  \ll   \ll \tau^*, \ve,   
\end{equation}
where  $n^*_1=n_1(\eta^2/2)$ is defined in Theorem~\ref{thm:PS}, $n_2 =n_2(\ve)$  is specified in Theorem~\ref{thm:S},  and $n_3$ is defined in Lemma~\ref{lem:partition}.

If $\chi'(G)=\Delta$, then clearly $G$ contains no $\Delta$-overfull subgraph. 
Thus we assume that $G$  contains no $\Delta$-overfull subgraph  and show $\chi'(G)=\Delta$. 
Suppose to the contrary that $\chi'(G)=\Delta+1$. Then $G$ is class 2 and so   $G$ contains 
an edge-chromatic critical class 2 subgraph with the same maximum degree. The graph  $G$ containing  no $\Delta$-overfull subgraph
implies that any subgraph of $G$ contains no $\Delta$-overfull subgraph.  Furthermore, as  $\Delta \ge (1-\ve) n$,  the subgraph contains 
more than $(1-\ve) n$ vertices.  By our setting up of the lower bound   $n_0$ that  $\frac{2}{n_0 } \le  \min \left \{\frac{1}{n^*_1}, \frac{1}{n_2}, \frac{1}{n_3},   \frac{1}{n^*_4} \right \} $ in ~\eqref{eqn:parameters} and so $\frac{1}{(1- \ve)n_0 } \le   \min \left \{\frac{1}{n^*_1}, \frac{1}{n_2}, \frac{1}{n_3},   \frac{1}{n^*_4} \right \} $, 
we may simply assume  that $G$ is edge-chromatic critical. 
 Since $G$ is not overfull, 
by Theorem~\ref{thm:overfull-min-degree}, we have 
$\Delta< \frac{7 \delta}{4}+ \frac{3n-17}{4}$. Thus 
\begin{eqnarray}
\delta &>&\frac{1}{7}  \left(4\Delta - (3n-17) \right)  
>  \frac{1}{7}  \left(1-4\ve  \right)n+2. 
\label{eqn: G-min-degree}
\end{eqnarray}
By Theorems~\ref{thm:PS} and~\ref{thm:S}, we may assume that 
\begin{equation}\label{eqn:upper-bound-delta}
\delta <\frac{1}{2}(1+\ve) n. 
\end{equation}
Then by  Vizing's Adjacency  Lemma, Lemma~\ref{lem:val}, we know that 
\begin{equation}\label{eqn:number-maximum-degree-vertex}
	|V_\Delta| \ge \Delta -\delta +1 >\left (\frac{1}{2}-\frac{3\ve}{2} \right) n.  
\end{equation}

Let 
$$U=\{v\in V(G)\,:\, \Delta-d_G(v) \ge    \eta n \}.$$  
We consider two cases regarding the size of $U$.  

\medskip 

{\bf \noindent Case A.  $|U| \ge 2\eta n$. } 

\medskip 

Recall that   the deficiency of  a vertex $u$
in $G$ is $\df_{G}(u)=\Delta-d_G(u)$ and $\df(G) =\sum_{u\in V(G)} \df_G(u)$. 
We split  Case A  into two subcases . 

\medskip 

{\bf \noindent Subcase A.1: $\df(G) \ge 1.01\ve n^2$.} 

\medskip 

In this case,  we add vertices to $G$ and construct a simple graph $H$ on an even number of vertices in three steps such that $G\subseteq H$, 
$\Delta(H)=\Delta$,     $\delta(H) \ge \frac{1}{2}(1+ \frac{1}{2}\eta^2) |V(H)|$, and $H$ contains no $\Delta$-overfull subgraph.  Then  we get $\chi'(G) \le \chi'(H)=\Delta(H)=\Delta$ by applying Theorem~\ref{thm:PS}, showing a contradiction to the assumption that  $\chi'(G)=\Delta+1$. 

Since $|U|\ge 2\eta n$, we have $\df(G)\ge 2\eta^2 n^2$.  We construct $H$ in three steps as follows. 

{\bf Step 1}:  We add a set  $W$ of  new vertices  to $G$  such that $ |W| \in \{ \lfloor  2\Delta -n-\eta^2 n \rfloor, \lfloor 2\Delta -n-\eta^2 n \rfloor-1\}$  
and  $|W| \equiv  n \pmod{2} $. 

The condition on the size of $W$ gives $\Delta \le \frac{1}{2}(|W|+n+\eta^2 n)+1 $.  
As $  (1-\ve) n\le \Delta \le n-1$, the condition on the size of $W$  also gives $ (1-2\ve-\eta^2) n-2\le |W| <(1-\eta^2) n$. Thus   
\begin{eqnarray*}
	\Delta -\df(G)/ |W|-2 & \le &  \frac{1}{2}(|W|+n+\eta^2 n) -1-\df(G)/|W| \\  
	&<& \frac{1}{2}(|W|+n+\eta^2 n-2) -1.01 \ve  n < \frac{1}{2}(|W|+n-2\ve n-\eta^2 n-2)  \le |W|. 
\end{eqnarray*}
Let $d$  be the even number in $\{ \lfloor \Delta -\df(G)/ |W| \rfloor-3, \lfloor \Delta -\df(G)/ |W| \rfloor-2\}$, 
which is less than $|W|$ by the inequality above.  We claim also that $d\ge 2$. 
Note that $\df(G)  \le (n-|V_\Delta|) (\Delta-\delta) <\frac{1}{2} n^2$ by~\eqref{eqn:upper-bound-delta} and~\eqref{eqn:number-maximum-degree-vertex}. As  $|W| \ge (1-2\ve-\eta^2 )n -2$, it then follows that $ \df(G)/|W| <\Delta -5$. 
Thus $d\ge 2$.

{\bf Step 2}:  We obtain a $d$-regular graph on $W$ by 
applying Lemma~\ref{lem:graphical-regular}.

{\bf Step 3}: We add $\df(G)$ edges joining vertices of $V(G)\setminus V_\Delta$ and $W$. 

Let 
$W=\{w_1, w_2, \ldots, w_{|W|}\}$.   
In particular, we list these  $\df(G)$ edges  as  $e_1, e_2, \ldots, e_{\df(G)}$   and add them one by one 
such that 
\begin{enumerate}[(1)]
	\item For $v\in V(G)\setminus V_\Delta$, the $\df_G(v)$ edges joining $v$ with vertices  of $W$ are listed consecutively in the ordering above; 
	\item  For each $i\in [1, \df(G)-1]$, the edge $e_i$ is incident with the vertex $w_i$ from $W$, where the index $i$ in $w_i$ is taken  modular $|W|$. 
	If $\df(G) \not\equiv (|W|-1) \pmod{|W|}$, the edge $e_{\df(G)}$ is also incident with $w_{\df(G)}$. 
	If $\df(G) \equiv (|W|-1) \pmod{|W|}$, we add the edge $e_{\df(G)}$ such that it is incident with the vertex $w_1$. 
\end{enumerate}

We denote by $H$ the resulting multigraph.   As $\delta > \frac{1}{7}(1-4\ve ) n+2$ by~\eqref{eqn: G-min-degree}
and so 
\begin{eqnarray*}
	(\delta+|W|-2 ) -\Delta   &>& \frac{1}{7}(1-4\ve) n +2+ (2\Delta -n-\eta^2 n-2-2) - \Delta\\
	&= & \frac{1}{7}(1-4\ve) n  +\Delta -n-\eta^2 n-2\\
	&\ge &  \frac{1}{7}(1-4\ve) n+(1-\ve)n -n-\eta^2 n-2\\
	& \ge & \frac{1}{7} n-\frac{11.5}{7} \ve n \ge 0, 
\end{eqnarray*}
as $\ve  \le \frac{1}{12}$. 
Thus $\df_G(v) \le \Delta-\delta<|W|-2$
for any $v\in V(G)$.  Then by  the constraints (1) and (2) in Step 3 of constructing $H$, we know that $H$
is simple. 
Note that for any $v\in V(G)$, we have $d_H(v)=\Delta =\Delta(H)$;
and for any $v\in W$, we have 
\begin{eqnarray*}
	\Delta-4\le d+ \lfloor \df(G)/|W| \rfloor -1 \le &d_H(v) &\le d+ \lceil \df(G)/|W| \rceil +1 \le  \Delta. 
\end{eqnarray*}
As $|W| \equiv  n \pmod{2} $,  we know that the order of $H$ is even and 
\begin{eqnarray*}
	(2-2\ve-\eta^2)n-2 \le  2\Delta -\eta^2 n-2 \le & |V(H)| &=n+|W| \le 2\Delta -\eta^2 n <(2-\eta^2 ) n. 
\end{eqnarray*}
Thus we have $2\Delta \ge n+|W|+\eta^2 n$ and so 
$$
\delta(H) \ge \Delta-4 \ge  \frac{1}{2}(n+|W|+\eta^2 n)-4 >\frac{1}{2}\left(1+ \frac{1}{2}\eta^2 \right) |V(H)|. 
$$

Furthermore, by the constraint (2) in Step 3, we know that $H$ has at least two vertices of minimum
degree.   Thus $H$ contains no $\Delta$-overfull subgraph by Lemma~\ref{lem:overfull2}. 
Therefore $\chi'(H)=\Delta(H)=\Delta$ by applying Theorem~\ref{thm:PS}.   Since $G \subseteq H$, 
it follows that $\chi'(G)=\Delta$, a contradiction to $G$ being class 2.

\medskip 
{\bf \noindent Subcase A.2: $\df(G) < 1.01\ve n^2$.}  

\medskip 
If $\Delta$ is odd,  by Theorem~\ref{thm:H-cycle}, we take  off a   matching  $M$ saturating $V(G)$  if $n$ is even, and saturating all but one vertex of minimum degree   if $n$ is odd. 
Thus we may assume that $\Delta$ is even, and  still use $ (1-\ve) n$ as a lower bound for  $\Delta$ for notation simplicity. (To be rigorous, for example, we can let $ \ve':=1.001\ve$ and have $\ve'$ in the place of $\ve$.  However, as  all calculations involving $\ve$ also work for $1.01 \ve$,   we just use $\ve$ for notation simplicity.)
When $n$ is odd and $\Delta$ is odd,  the matching  $M$ saturates all but one vertex of $G$. 
Thus,   $\df(G) \ge |U| \times \eta n-1 \ge 2\eta^2 n^2 -1$ regardless the parity of $\Delta$ and $n$. 
We add vertices to $G$ and construct a simple graph $H$ on an even number of vertices in three steps such that $G\subseteq H$, 
$\Delta(H)=\Delta$,   and $H$ is a $\Delta$-regular robust $(\eta^2, \eta)$-expander. Then  we get $\chi'(G) =\Delta$ by applying Theorem~\ref{thm:robust-expander}, showing
a contradiction to the assumption $\chi'(G)=\Delta+1$.

We add a set  $W$ of  new vertices  to $G$  such that
$|W|$ is the smallest integer satisfying the following properties:
\begin{enumerate}[(1)]
	\item $|W| \equiv  n \pmod{2} $ and $|W| \ge \frac{1}{2}(\Delta+1)$; 
	\item $ |W| \ge  \Delta - \left\lceil\frac{\df(G)}{ |W|}  \right\rceil +2$. 
\end{enumerate}

Since  $\df(G) \ge 2\eta^2 n^2-1> \Delta+1$, we know that when $|W| =\Delta$ or 
$|W|=\Delta+1$, the conditions above are satisfied.  Thus a smallest 
integer $|W|$ satisfying the conditions above exists. 
As $ \lceil \df(G)/\Delta  \rceil  \ge  \df(G)/\Delta >2\eta n>2\eta \Delta$, we have $\Delta- \lceil \df(G)/\Delta  \rceil  <(1-2\eta) \Delta$. Thus 
\begin{eqnarray*}
	&& 	\df(G)/(\Delta- \lceil \df(G)/\Delta  \rceil ) - \lceil \df(G)/\Delta  \rceil    \\
	&>& \df(G)/ \left((1-2\eta ) \Delta  \right)-  \lceil \df(G)/\Delta  \rceil  \\ 
	& >&  \frac{\df(G)}{ \Delta}\left (\frac{1}{1-2\eta} -1 \right ) -1\\  
	&>& 2\eta n \left (\frac{1}{1-2\eta} -1 \right ) -1 \ge 2, 
\end{eqnarray*}
and so $- \lceil \df(G)/\Delta  \rceil  > - \df(G)/(\Delta- \lceil \df(G)/\Delta  \rceil ) +2$ and thus 
$\Delta - \lceil \df(G)/\Delta  \rceil  \ge \Delta- \lceil \df(G)/(\Delta - \lceil \df(G)/\Delta  \rceil ) \rceil  +2$. 
This implies 
that $|W| =\Delta -\lceil \df(G)/\Delta \rceil  $ satisfies the inequality in (2) above.  Furthermore, since $\df(G)  \le (n-|V_\Delta|) (\Delta-\delta) <\frac{1}{2}(\Delta^2-\Delta)$ by~\eqref{eqn:upper-bound-delta} and~\eqref{eqn:number-maximum-degree-vertex}
and $|W|  \le \Delta+1$, it follows that $\Delta - \lceil \df(G)/\Delta  \rceil \ge \frac{1}{2}(\Delta+1)$. 
Therefore  we have 
$|W| \le  \Delta- \lceil \df(G) /\Delta  \rceil+1$.  

Suppose $\df(G) \equiv \ell \pmod{|W|}$, where $\ell \in [0, |W|-1]$.  
Let $d_1=\ldots =d_{\ell }= \Delta- \lceil\frac{\df(G)}{ |W|} \rceil =d-1$, and $d_{\ell+1}=\ldots = d_{|W|} =\Delta- (\lceil\frac{\df(G)}{ |W|} \rceil -1)=d$. 
Since $|W| \equiv n \pmod{2}$, we know that 
\begin{eqnarray*}
	\sum\limits_{i=1}^{|W|} d_i= \Delta |W| -\df(G) \equiv \Delta n - \df(G) =2e(G)+\df(G)-\df(G)  \equiv 0 \pmod{2}. 
\end{eqnarray*}
Since  $\df(G)<1.01 \ve n^2$ by the assumption of Subcase A.2 and  $|W| \ge \frac{1}{2}(\Delta+1)$, we have  $d=\Delta- (\lceil\frac{\df(G)}{ |W|} \rceil -1) \ge 2$. 
By Constraint (2) of the choice of $W$, we have $|W| \ge d+1$. 
Now by Lemma~\ref{lem:degree-sequence-d-d-1}, there is a graph $R$
on $W$ that realizes $(d_1, \ldots, d_{|W|})$.

Let 
$W=\{w_1, w_2, \ldots, w_{|W|}\}$, where $d_R(w_i) =d-1$ for each $i\in [1, \ell]$.   We now add $\df(G)$ edges  joining vertices from $V(G)\setminus V_\Delta$ and $V(R)$. 
We list these  $\df(G)$ edges  as  $e_1, e_2, \ldots, e_{\df(G)}$   and add them one by one 
such that 
\begin{enumerate}[(1)]
	\item For $v\in V(G)\setminus V_\Delta$, the $\df_G(v)$ edges joining $v$ with vertices  of $W$ are listed consecutively in the ordering above; 
	\item  For each $i\in [1, \df(G)]$, the edge $e_i$ is incident with the vertex $w_i$ from $W$, where the index $i$ in $w_i$ is taken  modular $|W|$. 
\end{enumerate}

We denote by $H$ the resulting multigraph.   We claim that $H$ is simple. By the construction of $H$, it suffice to show that 
$\delta+|W| \ge \Delta$.    We first prove the following lower bound on the size of $W$. 

\begin{CLA}\label{claim:lower-bound-W}
	It holds that	$|W|  \ge   \frac{5}{6} \Delta$. 
\end{CLA}

\pf  By using the condition that  $|W| \ge \frac{1}{2}(\Delta+1)$ 
and $|W|$ is smallest satisfying the inequality 
$|W|  \ge \Delta - \left\lceil\frac{\df(G)}{ |W|}  \right\rceil +2 \ge  \Delta - \frac{\df(G)}{ |W|}  +1$, we have 
$$
|W|  \ge  \frac{(\Delta +1) +\sqrt{(\Delta+1)^2-4\df(G)}}{2}. 
$$
Since  $\df(G)<1.01\ve n^2$ and $\ve <\frac{1}{10}$, we have 
\begin{eqnarray*}
	&& (\Delta+1)^2-4.04\ve n^2 -  (\frac{2}{3} \Delta)^2 \\  
	&> &  \frac{5}{9} \Delta^2 -4.04\ve n^2   \ge \frac{5}{9}(1-\ve)^2 n^2 -4.04\ve n^2 \\
	&>& \frac{5}{9} \left (1-2\ve -\frac{36.36}{5} \ve \right )n^2 >0. 
\end{eqnarray*}
Thus we get 
\begin{eqnarray*}
	&& \frac{(\Delta +1) +\sqrt{(\Delta+1)^2-4\df(G)}}{2} - \frac{5}{6} \Delta \\
	&>&  \frac{\sqrt{(\Delta+1)^2-4\df(G)}}{2} -\frac{1}{3} \Delta \\ 
	&>&   \frac{\sqrt{(\Delta+1)^2-4.04\ve n^2}-\frac{2}{3} \Delta}{2} >0, 
\end{eqnarray*} 
showing that $|W|  \ge   \frac{5}{6} \Delta$. 
\qed 

Then using the condition on the size of $W$ again, we get 
\begin{eqnarray*}
	|W| &\ge & \Delta - \left\lceil\frac{\df(G)}{ |W|}  \right\rceil +2 
	 > \Delta - \frac{1.01\ve n^2}{(\frac{5}{6} \Delta)} +1. 
\end{eqnarray*}
Hence 
\begin{eqnarray*}
	(\delta+|W| )-  \Delta  &\ge & \delta- \frac{1.01\ve n^2}{(\frac{5}{6} \Delta)} +1 \\
	&>& \frac{1}{7}(1-4\ve )n- \frac{1.01\ve n^2}{5(1-\ve)n/6} + 1 \\
	&>&  \left(  \frac{1}{7} -\frac{4\ve}{7} -\frac{6.06\ve}{ 5(1-\ve)}\right) n \\
	&=& \frac{1}{7}\left((1-4\ve) -\frac{8.484 \ve}{ 1-\ve}\right) n >0, 
\end{eqnarray*}
since $(1-4\ve)(1-\ve) > 8.484 \ve$ when $\ve \le \frac{1}{14}$.  Thus $\df_G(v) \le \Delta-\delta\le  |W|$
for any $v\in V(G)$.  Then by  the constraints (1) and (2) in  constructing $H$ above, we know that $H$
is simple.  

Note that  $H$ is  $\Delta$-regular,  $\Delta$ is even, and $ n_H:=|V(H)| =n+|W| \le n+\Delta+1$.  Thus 
$\Delta \ge (1-\ve) n \ge (1-\ve)(n_H-(\Delta+1))$ or $(2-\ve) \Delta \ge (1-\ve) (n_H-1)$, implying 
\begin{equation}\label{eqn:Delta-with-n_H}
	\Delta \ge \frac{1-\ve}{2-\ve} (n_H-1)> \frac{3}{7}  n_H. 
\end{equation}

\begin{CLA}\label{claim:robust-H}
	The graph $H$ is a robust $(\eta^2/2, \eta)$-expander.  
\end{CLA}
\proof  Let $S\subseteq V(H)$ with $  \eta n_H\le |S| \le (1-\eta)n_H$. 
If $n_H-|S| \le  \Delta-\frac{\eta^2}{2} n_H$, then as  $H$ is $\Delta$-regular, 
every vertex  of $H$ has at least $\frac{\eta^2}{2} n_H$ neighbors from  $S$. 
Thus $RN_H(S)=V(H)$ and so $|RN_H(S)| \ge |S|+\frac{\eta^2}{2} n_H$. 
Thus we assume $n_H-|S| >  \Delta-\frac{\eta^2}{2} n_H$ or $|S|<n_H-\Delta+\frac{\eta^2}{2} n_H$.  We consider two cases as follows. 

{\bf \noindent  Case 1: $|S| < \Delta-0.6\eta n_H$. }

Then we have 
\begin{eqnarray*}
	|S| \Delta = & \sum_{v\in S} d_H(v)  &=\sum_{v\in S} d_H(v, RN_H(S))  +\sum_{v\in S} d_H(v, V(H) \setminus RN_H(S))\\ 
	&& \le |S| |RN_H(S)|+  \frac{\eta^2}{2} n^2_H.
\end{eqnarray*}
As $|S| \ge \eta n_H$, we have 
\begin{eqnarray*}
	|RN_H(S)| &\ge & \Delta -\eta^2 n^2_H/(2|S|) \ge \Delta- \frac{1}{2}\eta n_H\\
	&=&  \Delta -0.6\eta n_H+ 0.6\eta n_H-\frac{1}{2}\eta n_H \\
	&>& |S|+ 0.1\eta n_H \ge  |S|+\frac{\eta^2}{2} n_H.  
\end{eqnarray*}

{\bf \noindent  Case 2: $ \Delta-0.6\eta n_H \le |S| <n_H-\Delta+\frac{\eta^2}{2} n_H$. }

By the choice of the size of $|W|$ and Claim~\ref{claim:lower-bound-W}, we must have $|W|-2<d+1$ and so $\delta(R)=d-1 > |W|-4$. 
We consider two subcases to finish the proof.

{\bf \noindent  Subcase 2.1:  $|S\cap W| < \frac{\eta^2}{2} n_H+3$. } 

Since $|S|> \Delta-0.6\eta n_H$, we know that $ |S\cap V(G)| > \Delta-0.6\eta n_H - (\frac{\eta^2}{2} n_H+3) > \Delta-0.7\eta n_H$.  
Thus all vertices of $G$ with degree at least $$n- (\Delta-0.7\eta n_H) +\frac{\eta^2}{2} n_H <\ve n+\eta n_H$$ 
have at least 
$\frac{\eta^2}{2} n_H$ neighbors from $S\cap V(G)$, and so  are contained in $RN_H(S)$. 
As $\delta(G) >\frac{1}{7}(1-4\ve) n>\ve n+ \eta n_H$ implied by $\ve \le \frac{1}{14}$, it follows that $V(G) \subseteq RN_H(S)$. 
If $|U\setminus  S| \ge |S\cap W|+\frac{\eta^2}{2} n_H$,  then we get 
$$ |RN_H(S)| \ge  |S\cap V(G)|+|V(G)\setminus S|\ge |S\cap V(G)|+|U\setminus S|\ge   |S\cap V(G)|+ |S\cap W|+\frac{\eta^2}{2} n_H= |S|+\frac{\eta^2}{2} n_H $$
already. 
Thus we assume that $|U\setminus  S| < |S\cap W|+\frac{\eta^2}{2} n_H$. As $|S\cap W| < \frac{\eta^2}{2} n_H+3$,
we have $|U\setminus  S|  <\eta^2 n_H+3$. 
By the construction of $H$,  every vertex of  $W$ has in $H$ at least $\lfloor \frac{\df(G)}{|W|} \rfloor \ge 2\eta n$  neighbors 
from $U$. So for any $v\in W$, we have 
$$|N_H(v)\cap (U\cap S)| \ge  2\eta n-(U\setminus S) \ge 2\eta n-(\eta^2 n_H+3) >\frac{\eta^2}{2} n_H.$$ This implies that $W \subseteq RN_H(S)$. 
Combining the arguments above, we get 
$|N_H(S)|=|V(G)|+|W| =n_H >|S|+\frac{\eta^2}{2} n_H$. 

	{\bf \noindent  Subcase 2.2:  $|S\cap W| \ge \frac{\eta^2}{2} n_H+3$. } 
	
	Then as $\delta(R)=d-1 > |W|-4$, we get $W\subseteq RN_H(S)$.  
	Consider first that $|S\cap V(G)| \ge  n-\Delta+\frac{\eta^2}{2} n_H$.  This implies that every vertex from $V_\Delta$
	has at least $\frac{\eta^2}{2} n_H$ neighbors in $S\cap V(G)$ and so $V_\Delta \subseteq RN_H(S)$. 
	As $|V_\Delta| \ge \Delta-\delta+1$ by~\eqref{eqn:number-maximum-degree-vertex} and 
	$2\Delta > n+\delta +\eta^2 n_H$ (this is because $\Delta\ge (1-\ve)n$, $\delta< \frac{1}{2}(1+\ve) n$ by~\eqref{eqn:upper-bound-delta},
	$\ve \le \frac{1}{24}$, and $n_H<2n$), 
	we get 
	\begin{eqnarray*}
		|RN_H(S)|  &\ge &|W|+|V_\Delta|  \ge |W|+\Delta-\delta+1 \\
		&=&|W|+2\Delta -\Delta-\delta+1  \\  
		& >& n+|W|-\Delta +\frac{\eta^2}{2} n_H+\frac{\eta^2}{2} n_H \\
		&=& n_H-\Delta +\frac{\eta^2}{2} n_H+\frac{\eta^2}{2} n_H \\
		&>&|S| +\frac{\eta^2}{2} n_H. 
	\end{eqnarray*}
	Thus we assume 	$|S\cap V(G)| <  n-\Delta+\frac{\eta^2}{2} n_H$. 
	Let $S_1=(S\cap V(G)) \setminus RN_H(S)$. 
	If $|W\setminus S| \ge |S_1|+ \frac{\eta^2}{2} n_H$, then we again get $|RN_H(S)| \ge |S| +\frac{\eta^2}{2} n_H.$ 
	Thus we also assume $|W\setminus S| < |S_1|+ \frac{\eta^2}{2} n_H$.

	If $|S_1| \le  \frac{1}{2}\eta n $,  then $|W\setminus S| < |S_1|+ \frac{\eta^2}{2} n_H < \frac{1}{2}\eta n + \frac{\eta^2}{2} n_H$. 
	Since every vertex from $U$ has in $H$ at least $\eta n$ neighbors   from $W$, 
	we then know that   $U\subseteq RN_H(S)$ and so 
	\begin{eqnarray*}
		|RN_H(S)| &\ge& |W|+|(S\cap V(G)) \setminus S_1|+|U|  \ge |S\cap W| +|(S\cap V(G)) \setminus S_1|+|S_1|+(|U| -|S_1|)  \\ 
		&=&|S| +(|U| -|S_1|)   \ge |S| +(2\eta n-\frac{1}{2} \eta n) \\
		&>& |S| +\frac{\eta^2}{2} n_H.  
	\end{eqnarray*}
	
	Thus we assume $|S_1| >   \frac{1}{2}\eta n $.  As every vertex of $S_1$ has in  $H $ less than $  \frac{\eta^2}{2} n_H$ neighbors from $W\cap S$, 
	we know for any $v\in S_1$, 
	\begin{eqnarray*}
		d_G(v) &\ge&  \Delta-(|W\setminus S| + \frac{\eta^2}{2} n_H)  \\
		&\ge & \Delta-(|S_1| + \frac{\eta^2}{2} n_H+\frac{\eta^2}{2} n_H) \\
		&\ge& \Delta-|S\cap V(G)|-\eta^2n_H  \\
		&\ge&  \Delta -(n-\Delta+\frac{\eta^2}{2} n_H)-\eta^2n_H   \\
		&=& 2\Delta -n - \frac{3}{2}\eta^2 n_H \ge  \frac{1}{2}n. 
	\end{eqnarray*}
	Therefore 
	\begin{eqnarray*}
		\frac{1}{2} n |S_1|   \le & \sum_{v\in S_1}d_G(v)   &=\sum_{v\in S_1}d_G(v, RN_G(S_1))+ \sum_{v\in S_1}d_G(v, V(G) \setminus RN_G(S_1))\\ 
		&& \le |S_1| |RN_G(S_1)|+  \frac{1}{2}\eta^2 n^2_H.
	\end{eqnarray*}
	As $|S_1| >  \frac{1}{2}\eta n > \frac{1}{4}\eta n_H$, the inequalities 
	above  imply that 
	\begin{eqnarray*}
		|RN_G(S_1)|  &\ge & \frac{1}{2}n-\eta^2 n^2_H/(2|S_1|)    \\  
		&\ge & \frac{1}{2}n - 2\eta n_H  > \ve n+ \eta^2 n_H  \\
		& \ge  &  n-\Delta+\frac{\eta^2}{2} n_H +\frac{\eta^2}{2} n_H  \\
		& >& |S\cap V(G)|  + \frac{\eta^2}{2} n_H. 
	\end{eqnarray*}
	Therefore, $ |RN_H(S)| \ge |W| + |RN_G(S_1)|  > |S|+\frac{\eta^2}{2} n_H. $
	\qed

	By Claim~\ref{claim:robust-H}, the graph $H$ is a robust $(\eta^2/2, \eta)$-expander.   
	Furthermore, $H$ is $\Delta$-regular, $\Delta$ is even, and $\Delta \ge  \frac{3}{7} n_H$ by~\eqref{eqn:Delta-with-n_H}. 
	As $\eta \ll \tau^*$ by~\eqref{eqn:parameters}, $H$ is a robust $(\eta^2/2, \eta)$-expander implies that $H$ is  a robust $(\eta^2/2, \tau^*)$-expander. 
	As $n$ is taken to be at least  $n_4$  as defined in Theorem~\ref{thm:robust-expander}, Theorem~\ref{thm:robust-expander} implies that $H$ 
	has a Hamilton decomposition.  As $n_H$ is even, we conclude that $\chi'(H)= \Delta$.  As $G\subseteq H$ and $\Delta(G)=\Delta(H)$, it follows 
	that $\chi'(G)=\Delta$, a contradiction.

\medskip 

{\bf \noindent Case B.  $|U|<2\eta n$. } 

\medskip

The most obstacle comes in  scenarios like the following:   suppose $n$ is even, $d_G(v_1)=d_G(v_2)=d_G(v_3) =\frac{1}{3} \Delta$, and all the other vertices of $G$ have degree $\Delta$ in $G$. The only way to construct a $\Delta$-regular multigraph  $G'$ based on $G$ is to add $\frac{1}{3} \Delta$ edges between each pair of vertices of $\{v_1, v_2, v_3\}$.  However, the multiplicity of 
the resulting multigraph $G'$  is  at least $\frac{1}{3} \Delta$ and no matter how  we partition  $V(G')$ into subsets $A$ and $B$, the multiplicity of one of $G'[A]$ and $G'[B]$ is going to be at least $\frac{1}{3} \Delta$.   This makes the maximum degree of one of $G'[A]$ and $G'[B]$ 
to be about  $\frac{1}{2} \Delta +\frac{1}{6} \Delta $. This presents a significant challenge for the second step of the procedure, as will be detailed later,  where we require the difference between the maximum degree of $G'[A]\cup G'[B]$
and $\frac{1}{2} \Delta$ to be much smaller  than $\frac{1}{2} \delta$. 
This challenge, compounded by the significantly weaker minimum degree lower bound of  
$\frac{1}{7}(1-4\ve) n$  also hinders the fourth step. There, the typical application 
is for $n$-vertex graphs with  minimum degree   at least $\frac{1}{2}(1+\ve) n$. 
 Case B is largely devoted to dealing with the two challenges above. We will construct 
 a $\Delta$-regular multigraph   $G_3$ such that $G\subseteq G_3$ and   $G_3$ contains no $\Delta$-overfull subgraph. 
Disregard some multiple edges,  we will decompose the edges of a spanning subgraph of $G_3$ into  edge-disjoint matchings by following  the edge-decomposition framework described in Section 1. 

\subsection{Construct a $\Delta$-regular multigraph $G_3$: $G \rightarrow G_0\rightarrow G_1\rightarrow G_2\rightarrow G_3$}

We will add edges to $G$ (when $n$ is even) or the union of $G$ and a new vertex (when $n$ is odd)
to obtain a $\Delta$-regular multigraph called  $G_3$ such that $G_3$ contains no $\Delta$-overfull subgraph.

\hl{Construction of $G_0$}. \quad 
If $n$ is even, we let $G_0=G$ and $m=n$. 
If $n$ is odd, we let $G_0$ be the union of $G$ and a new single vertex, and  let $m=n+1$.   

\medskip 

\hl{Construction of $G_1$}.\quad 
The graph $G_0$ contains no  $\Delta$-overfull subgraph
as $G$ contains no $\Delta$-overfull subgraph. 
 Let $V(G_0) =\{v_1,  \ldots, v_m\}$ and suppose $d_{G_0}(v_1) \le \ldots \le d_{G_0}(v_m)$. 
Let $d_i=\df_{G_0}(v_i):= \Delta-d_{G_0}(v_i)$ for each $i\in [1, m]$.  As $m$ is even, it follows that $\sum_{i=1}^n d_i =m \Delta -2e(G_0)$ is even. 
Since $G_0$ contains no  $\Delta$-overfull subgraph, we have 
\begin{eqnarray*}
2e(G_0-v_1)&=&\sum_{i=2}^m d_{G_0}(v_i)- d_{G_0}(v_1)
=\sum_{i=2}^m(\Delta-d_i)-(\Delta-d_1) 
 \le   \Delta(m-2). 
\end{eqnarray*}
This gives  $d_1\le \sum_{i=2}^m d_i$. Thus $(d_1, \ldots, d_m)$ is a feasible sequence. 
Applying  Lemma~\ref{lem:graphical-biparite}, we find a multigraph $L$ on $\{v_1, \ldots, v_m\}$ 
and an even index $p\in [2,m]$
such that $L$ 
satisfies all the properties listed in~Lemma~\ref{lem:graphical-biparite}.

Let $G_1$ be  the multigraph obtained from $G_0$ by adding all the edges of $L$ (multiple edges may be created).  
It is clear that $\Delta(G_1)=\Delta$
and $\{v_1, \ldots, v_p, v_{p+1}\} $ are all the possible vertices of $G_1$ with degree less than $\Delta$.  In particular, 
we have 
\begin{eqnarray}
d_{G_1}(v_{2i-1}) &=&d_{G_1}(v_{2i}) =d_G(v_{2i})=\Delta-d_{2i} \quad \text{for each $i\in [1,\frac{p}{2}]$}; \label{eqn:deg-in-G1} \\
\df_{G_1}(v_{p+1}) &\equiv& 0 \pmod{2} \quad \text{as $m=|V(G_1)|$ is even and so $\sum_{v\in V(G_1)} \df_{G_1}(v)$ is even.}\label{eqn:def(v-p+1)-even}
\end{eqnarray}

\hl{Construction of $G_2$}.\quad 
We define a multigraph $G_2$ based on $G_1$ such that $$e_{G_2}(\{v_{p-1}, v_p, v_{p+1}\}, V(G_2) \setminus \{v_{p-1}, v_p, v_{p+1}\})  \ge \Delta.$$  This will ensure that the number of edges within $v_{p-1}, v_p, v_{p+1}$
is at most $\Delta$ when we completing a $\Delta$-regular multigraph based on $G_2$. 
The multigraph $G_2$ is obtained through the following algorithm. 
 \begin{Pro}[Vertex Identification Procedure]\label{opr:vertex-contraction} 
 	We apply the following algorithm on $G_1$: 
 	\begin{enumerate}[]
 		\item  Step 1: If  $e_{G_1}(\{v_{p-1}, v_p, v_{p+1}\}, V(G_1) \setminus \{v_{p-1}, v_p, v_{p+1}\}) <\Delta$, go to Step 2; Otherwise, go to Step 3. 
 		\item  Step 2: Identify  $v_{p}$ and $v_{p+1}$ with   $v_{p-1}$ (remove any resulting loops).  Update $G_1$, set $p:=p-2$, and go to Step 1.  (The identified vertex plays the role of $v_{p+1}$.)
 		\item Step 3: Name the current multigraph as $G_2$ and stop.  
 	\end{enumerate}
 \end{Pro}
 Since $G_1-\{v_1, \ldots, v_{p+1}\} \subseteq G$ is not $\Delta$-overfull and all vertices from $G_1-\{v_1, \ldots, v_{p+1}\}$ have degree $\Delta$ in $G_1$, it follows that 
 \begin{eqnarray*}
 	2e(G_1-\{v_1, \ldots, v_{p+1}\} )& =& \Delta (m-(p+1))-e_{G_1}(\{v_1, \ldots, v_{p+1}\}, V(G_1) \setminus \{v_1, \ldots, v_{p+1}\})  \\
 	& \le& \Delta(m-(p+1) -1).
 \end{eqnarray*}
 This implies that $e_{G_1}(\{v_1, \ldots, v_{p+1}\}, V(G_1) \setminus \{v_1, \ldots, v_{p+1}\}) \ge \Delta$. 
 Thus the procedure above  will stop.  In the multigraph $G_2$, we have 
 \begin{equation}\label{claim:degree-of-v-p+1}
 e_{G_2}(\{v_{p-1}, v_p, v_{p+1}\}, V(G_2) \setminus \{v_{p-1}, v_p, v_{p+1}\})  \ge \Delta.
 \end{equation}

To clarify the construction of $G_3$ in the next step, 
we examine some properties of $G_2$ below.

\begin{CLA}\label{claim:G2-properties}
	The following statements hold. 
	\begin{enumerate}[(i)]
		\item $\chi'(G_2) =\Delta$ implies that  $\chi'(G) =\Delta$; 
		\item  For any $v\in \{v_1, \ldots, v_p\}$, we have $N_{G_2}(v) = N_{G_1}(v)$ and 
		  $d_{G_2}(v_{2i-1}) =d_{G_2}(v_{2i}) =d_G(v_{2i})=\Delta-d_{2i}$  for each $i\in [1,\frac{p}{2}]$. 
		Moreover, $N_G(v_{p+1})\subseteq  N_{G_2}(v_{p+1})$ and  $d_{G_2}(v_{p+1}) 
		\ge d_G(v_{p+1})$; 
		\item For any $v\in V(G_2)\setminus \{v_1, \ldots, v_{p+1}\}$, we have $d_{G_2}(v) =d_{G_1}(v)$ and  $d_{G_2}^s(v) \ge  d^s_{G_1}(v)-|U| \ge d^s_{G_1}(v) -2\eta n$; 
		\item For any $v\in V(G)\cap V(G_2)$, we have $d_{G_2}^s(v) \ge \delta$; 
		\item $d_{G_2-\{v_{p-1}, v_p\}}(v_{p+1}) \ge  \frac{\Delta}{3}$; 
		\item $L\cap G_2$ still satisfies all the properties listed in Lemma~\ref{lem:graphical-biparite}
		with bipartition as $\{v_1,\ldots, v_p\}$ and $V(L\cap G_2)\setminus \{v_1,\ldots, v_p\}$. 
	\end{enumerate}
\end{CLA}

\pf   Note that the index $p$ in $G_2$ might be different from the index $p$ in $G_1$. 
In the proof of this claim, we use $q$ for the old index $p$.  Thus $q\ge p$. 
If $G_2=G_1$, then  Claims(i)-(iv) are obviously true. Thus we assume $G_2\ne G_1$ in the   proof of Claims(i)-(iv).  
The assumption $G_2\ne G_1$ implies that  $G_2$ was obtained from $G_1$ by identifying vertices through Procedure~\ref{opr:vertex-contraction},
and so  $e_{G_1}(\{v_{q-1}, v_{q}, v_{q+1}\}, V(G_1) \setminus \{v_{q-1}, v_q, v_{q+1}\}) <\Delta$. We claim 
below that 
\begin{equation}\label{eqn:ind-set}
\text{$e_{G_1}(\{v_{q-1}, v_{q}, v_{q+1}\}, V(G_1) \setminus \{v_{q-1}, v_q, v_{q+1}\}) <\Delta$ implies that $\{v_1, \ldots, v_{q+1}\}$ is  independent  in $G_0$}.
\end{equation}
As $G$ is edge-chromatic critical,  Vizing's Adjacency Lemma implies that the degree sum of any two adjacent vertices in $G$ is at least $\Delta+2$. 
Since  $d_{G_0}(v_1) \le d_{G_0}(v_2) \le \ldots \le d_{G_0}(v_{q-1})$, it follows that $d_{G_0}(v_q)+d_{G_0}(v_{q+1}) \ge d_{G_0}(u)+d_{G_0}(v)$ for any two distinct $u,v\in \{v_1, \ldots, v_{q+1}\}$. Thus if $uv\in E(G)$ for some $u,v\in \{v_1, \ldots, v_{q+1}\}$, then 
\begin{eqnarray}
	&& e_{G_1}(\{v_{q-1}, v_q, v_{q+1}\}, V(G_1) \setminus \{v_{q-1}, v_q, v_{q+1}\})   \nonumber \\
	& \ge & d_{G_0-\{v_{q}, v_{q+1}\}}(v_{q-1}) +d_{G_0-\{v_{q-1}, v_{q+1}\}}(v_{q}) +d_{G_0-\{v_{q-1}, v_q\}}(v_{q+1}) \nonumber  \\
	&\ge&  \Delta+2+d_{G_0}(v_{q-1})-2\left(e_{G_0}(v_{q-1}, v_q)+e_{G_0}(v_{q-1}, v_{q+1})+e_{G_0}(v_{q}, v_{q+1}) \right). \label{eqn:degree-sum}
\end{eqnarray}
If $d_{G_0}(v_{q-1}) =0$, then we have $e_{G_0}(v_{q-1}, v_q)=e_{G_0}(v_{q-1}, v_{q+1})=0$ and obtain a contradiction 
from~\eqref{eqn:degree-sum}. Thus $d_{G_0}(v_{q-1}) \ge \delta >\frac{1}{7}(1-4\ve)n>6$. 
Since  $2(e_{G_0}(v_{q-1}, v_q)+e_{G_0}(v_{q-1}, v_{q+1})+e_{G_0}(v_{q}, v_{q+1}) )\le 6$, we obtain a contradiction 
from~\eqref{eqn:degree-sum} again. Thus $\{v_1, \ldots, v_{q+1}\}$ is an independent set in $G_0$.

For (i),   as  $G_2\ne G_1$, 
we know that  $\{v_1, \ldots, v_{q+1}\}$ is an independent set in $G_0$  by~\eqref{eqn:ind-set}.
Since only vertices from  $\{v_{p+2}, \ldots, v_{q+1}\}$ were identified together, we have 
$e_{G_2}(\{v_1, \ldots, v_{p+1}\}, V(G_2) \setminus \{v_1, \ldots, v_{p+1}\}) =e_{G_0}(\{v_1, \ldots, v_{q+1}\}, V(G_0) \setminus \{v_1, \ldots, v_{q+1}\})$. 
As $\{v_1, \ldots, v_{q+1}\}$ is an independent set in $G_0$, it follows that any  edge $\Delta$-coloring  $\varphi$ of $G_2-E(G_2[\{v_1, \ldots, v_{p+1}\}])$
is an  edge $\Delta$-coloring of $G_0$ as all edges incident to the identified vertex  $v_{p+1}$ of  $G_2$
are colored with distinct colors  under $\varphi$.  Thus $\chi'(G_2) =\Delta$ implies that $\chi'(G)  \le \chi'(G_0)=\Delta$ and so $\chi'(G) =\Delta$.

For (ii), 
since $\{v_1, \ldots, v_{q}\}$ is an independent set in both $G_0$ and $G_1$  and 
 only vertices of  $\{v_{p+2}, \ldots, v_{q+1}\}$ were identified into $v_{p+1}$ when we get $G_2$ from $G_1$, 
it follows that for any 
$v\in \{v_1, \ldots, v_{p}\}$, $e$ is an edge incident with $v$ in $G_1$ if and only if  $e$ is an edge  incident 
with $v$ in $G_2$.  This 
together with~\eqref{eqn:deg-in-G1}, proves the first part of (ii). 
The second part of (ii) follows by the following fact:  by the constructions of $G_1$ and $G_2$,   
if $e$ is an edge incident with $v_{p+1}$ in $G$, then $e$ is   incident 
with $v$ in $G_1$ and $G_2$. (Note that it is possible to have an edge $e$ that is incident with $v_{p+1}$ in $G_1$ but  was removed in $G_2$ after applying  the Vertex Identification Procedure.)

For (iii), 
since  $d_{G_0}(v_q) \le d_{G_0}(v_{q+1})$ and $e_{G_1}(\{v_{q-1}, v_q, v_{q+1}\}, V(G_1) \setminus \{v_{q-1}, v_q, v_{q+1}\}) <\Delta$, it follows that $v_{1}, \ldots, v_{q}\in U$.  Thus   $q-p \le |U|$, and the conclusion follows.

For (iv), for any $v\in \{v_1, \ldots, v_{p+1}\}$ or $v\in V(G_2)\setminus U$ with $v\in V(G)$, we have $d_{G_2}^s(v) \ge \delta$ by (ii) and (iii). 
Suppose $v\in  U\setminus \{v_1, \ldots, v_{p+1}\}$. If $v$ is not adjacent in $G$ to any vertex from $\{v_1, \ldots, v_{q+1}\}$, then we have 
$d_{G_2}^s(v) \ge d_G(v) \ge \delta$. Thus there exists $i\in [1,q+1]$ such that $vv_{i} \in E(G)$. 
Then as $d_G(v) \ge d_G(v_i)$, VAL (Lemma~\ref{lem:val}) implies that $d_G(v) \ge \frac{2\Delta}{3}$ if $i \le q-1$ (as $d_G(v_i) \le d_{G_0}(v_{q-1}) \le \frac{\Delta}{3}+2$
and $d_G(v) \ge d_G(v_{q+1})$ by   $d_G(v_1) \le \ldots \le d_G(v_m)$),
and so  $d_{G_2}^s(v) \ge d_G(v)-|U| \ge \frac{2\Delta}{3}-2\eta n  \ge \delta$. 
Thus  we assume that $v$ is   adjacent in $G$ to a vertex  from $\{v_q, v_{q+1}\}$ but is not 
adjacent to any other vertex $v_i$ with $i\in [1,q-1]$.   If $d_G(v) \ge \delta+1$
or $v$ is adjacent in $G$ to only one of $v_q$ and $v_{q+1}$,  then we again get 
$d_{G_2}^s(v) \ge  \delta$. Thus we assume that $d_G(v) =\delta$ and $vv_q, vv_{q+1} \in E(G)$. 
Then as $d_G(v)+d_G(v_q) \ge \Delta+2$ and $d_G(v) =\delta$, we get  $d_G(v_q)+d_G(v_{q+1}) =\delta+\delta \ge \Delta+2$, 
a contradiction to the assumption $e_{G_1}(\{v_{q-1}, v_q, v_{q+1}\}, V(G_1) \setminus \{v_{q-1}, v_q, v_{q+1}\}) <\Delta$.

For (v), if we have $uv\in E(G)$ for any distinct  $u, v\in \{v_1, \ldots, v_{p+1}\}$, then we 
have $d_G(v_{p+1}) \ge \max\{d_G(u), d_G(v)\} \ge \frac{1}{2} \Delta+1$
by Vizing's Adjacency Lemma. Thus $d_{G_2-\{v_{p-1}, v_p\}}(v_{p+1}) \ge  \frac{1}{2} \Delta+1 -2>\frac{\Delta}{3}$. 
Thus we assume that $\{v_1, \ldots, v_{p+1}\}$ is an independent set in $G_0$.  This implies that $\{v_1, \ldots, v_{p+1}\}$ 
is an independent set in $G_2$ by the construction of $G_1$ and $G_2$. 
Thus we have $d_{G_2-\{v_{p-1}, v_p\}}(v_{p+1})=d_{G_2}(v_{p+1})$,   $d_{G_2-\{v_{p-1}, v_{p+1}\}}(v_{p}) = d_{G_2}(v_p)$, and  $d_{G_2-\{v_{p-1}, v_{p+1}\}}(v_{p-1}) = d_{G_2}(v_{p-1})$. 
 As $d_{G_2}(v_{p+1})  \ge  d_{G}(v_{p+1}) $ by (ii),  $d_{G}(v_{p+1}) \ge d_{G}(v_p)$,  
 and $ d_{G}(v_p)=d_{G_2}(v_p)=d_{G_2}(v_{p-1})$ by  (ii), it follows that  $d_{G_2-\{v_{p-1}, v_p\}}(v_{p+1})  \ge d_{G}(v_{p+1}) \ge  d_{G}(v_p) =d_{G_2}(v_{p}) =d_{G_2-\{v_{p-1}, v_{p+1}\}}(v_{p}) =d_{G_2}(v_{p-1})= d_{G_2-\{v_{p}, v_{p+1}\}}(v_{p-1})$. Thus from~\eqref{claim:degree-of-v-p+1},
 we get   $d_{G_2-\{v_{p-1}, v_p\}}(v_{p+1}) \ge  \frac{\Delta}{3}$.

For (vi), as the identified vertices $v_{p+2}, \ldots, v_{q+1}$ have consecutive indices,
and they were all identified into $v_{p+1}$,  $L$ satisfying  the properties 
listed in Lemma~\ref{lem:graphical-biparite} implies that $L\cap G_2$ satisfies the properties 
listed in Lemma~\ref{lem:graphical-biparite}, and $p$ is the desired index for the bipartition of $L\cap G_2$.  
\qed 

By Claim~\ref{claim:G2-properties}(i), we will build $G_3$ based on $G_2$ for our goal of achieving 
a contradiction to the assumption $\chi'(G) \ge \Delta+1$.   Furthermore, for notation simplicity, 
by Claim~\ref{claim:G2-properties}(vi), 
we will just use $L$ for $L\cap G_2$ in the rest of the paper. 
By Claim~\ref{claim:G2-properties}(ii)-(iii), 
we know that  vertices of $\{v_1, \ldots, v_{p+1}\}$ are the only possible vertices of degree degree less than 
$\Delta$ in $G_2$, and  $\df_{G_2}(v_{2i-1}) =\df_{G_2}(v_{2i}) =d_{2i}$
 for each $i\in [1,\frac{p}{2}]$. As $\df(G_2)$ is even (note that  $|V(G_2)|$ is even, as three vertices were identified into one when each time Step 2 of the Vertex Identification Procedure was applied), it follows that $\df_{G_2}(v_{p+1})$ is even. 
 
 \medskip 
\hl{Construction of $G_3$}.\quad 
 We add edges to $G_2$ to form a $\Delta$-regular multigraph  $G_3$. 
 \begin{Pro}[Edge Addition Procedure]\label{Pro2:Edge-Addition}
 	We add the following edges to $G_2$ to get $G_3$: 
 	\begin{enumerate}[(i)] 
 		\item $\frac{1}{2}\df_{G_2}(v_{p+1})$ edges ($\df_{G_2}(v_{p+1})$ is even) respectively between  $v_{p+1}$ and $v_{p-1}$ and between $v_{p+1}$ and $v_{p}$.  (Thus the degree of $v_{p+1}$ in the resulting multigraph is $\Delta$.)
 		\item  $\df_{G_2}(v_{p-1}) -\frac{1}{2}\df_{G_2}(v_{p+1})$ edges between $v_{p-1}$ and $v_p$.
 		 (Thus  the degrees of $v_{p-1}$ and $ v_p$ in the resulting multigraph are $\Delta$.)
 		\item For each $i\in [1,\frac{p}{2}-1]$, we add $\df_{G_2}(v_{2i-1})$ edges between $v_{2i-1}$ and $v_{2i}$. 
 		(Thus the degrees of $v_{2i-1}$ and $ v_{2i}$ in the resulting multigraph are $\Delta$.)
 	\end{enumerate}
 \end{Pro}
 By the construction, the multigraph $G_3$ is $\Delta$-regular. 
We next  study  some crucial properties of  $G_3$.

\begin{CLA}\label{claim:G1-no-overfullness}
	The multigraph $G_3$ contains no $\Delta$-overfull subgraph. 
\end{CLA}

\pf   Suppose to the contrary that $G_3$ contains a $\Delta$-overfull subgraph $R$.    Then there exists $X \subseteq V(G_3)$ 
such that $R=G_3[X]$ is $\Delta$-overfull.  As $G_3$
is regular,  $R=G_3[X]$ is $\Delta$-overfull if and only if $e_{G_3}(X, V(G_3)\setminus X) \le \Delta-2$. 
This further implies that  $ 3\le |X| \le m-3$ and 
so $G_3-X$ is  $\Delta$-overfull as well. 
Therefore we assume, without loss of generality, that $|X| \le m/2$.

If $|X\cap (V(G_3)\setminus (U\cup \{v_1\})| \ge 3$ (note that when $n$ is odd,  $d_{G_0}(v_1) =0$ and $v_1\not\in U$),  
since $d^s_{G_3} (v) \ge \Delta -\eta n -2\eta n$ for any $v\in V(G_3)\setminus (U\cup \{v_1\})$ by Claim~\ref{claim:G2-properties}(iii), 
we then  get 
\begin{eqnarray*}
e_{G_3}(X, V(G_3) \setminus X)& \ge&  3(\Delta -3\eta n -|X|+1 )\\ 
 &\ge & 3\left( (1-\ve)n-3\eta n -\frac{n+1}{2}\right) >n >\Delta, 
 \end{eqnarray*}
a contradiction.   Thus $|X\cap (V(G_3)\setminus (U\cup \{v_1\})|  \le 2$ and so $|X| \le |U|+1+2 <2\eta n+3$.  
As $\delta^s(G_3) >  \frac{1}{7} (1-4\ve) n  \ge \frac{8}{7} \ve n \gg |X|$ when $m=n$ and $\delta^s(G_3-v_1)  >\frac{1}{7} (1-4\ve) n  \ge \frac{8}{7} \ve n \gg|X|$  if $m=n+1$ (recall $\ve \le 1/12$),  
by counting $e_{G_3}(X, V(G_3) \setminus X)$  the same way as above, we  further conclude that $|X| < \frac{1}{ \ve}$.
Since $R=G_3[X]$ is $\Delta$-overfull, we have $\chi'(R) \ge  \Delta +1$. 
By Theorem~\ref{konig}, $G_3[X]$  must contain an odd cycle. 

Consider first that $|X|=3$. 
If $X$ is an independent set in $G_0$, then $G_3[X] \subseteq G_3-E(G_0)$. 
Since $L$ is bipartite with one bipartition as $\{v_1, \ldots, v_p\}$, by the construction of $G_3$,  
we know that the vertex set of the only possible odd cycle  of $G_3-E(G_0)$ contains $\{v_{p-1},  v_p,  v_{p+1}\}$. 
As $|X|=3$,  it must be  that  $X=\{v_{p-1},  v_p,  v_{p+1}\}$. 
However,  \eqref{claim:degree-of-v-p+1}   implies that $G_3[X] $ is not $\Delta$-overfull, a contradiction. 
Thus $X$ is not an independent set in $G_0$. Let $u, v\in X$ such that 
$uv\in E(G)$. Then  as $G$ is edge-chromatic critical, we have $d_G(u)+d_G(v)  \ge \Delta+2$  
by VAL (Lemma~\ref{lem:val}).   Thus we have $d_{G_3}(u)+d_{G_3}(v)  \ge \Delta+2$  as  $d_{G_3}(u)+d_{G_3}(v) \ge d_G(u)+d_G(v)$ 
by Claim~\ref{claim:G2-properties}(ii) and (iii). 
Let $w\in X\setminus \{u,v\}$.  If $d_{G}(w) \ge 4$, then we get 
$e_{G_3}(X, V(G_3) \setminus X) \ge e_{G}(X, V(G) \setminus X) \ge \Delta+2+4-6 =\Delta$, a contradiction to $G_3[X] $ being  $\Delta$-overfull. Thus $d_G(w) \le 3$ 
and this implies that $w=v_1$ and $n$ is odd. 
Then we get $e_{G_3}(X, V(G_3) \setminus X) \ge e_{G_0}(X, V(G_0) \setminus X)  \ge  \Delta+2-2 =\Delta$, a contradiction again.

Consider next that $|X| \ge 5$. If there are $u,v\in X$ with $uv\in E(G)$, then we again have  
$d_G(u)+d_G(v)  \ge \Delta+2$  
by VAL (Lemma~\ref{lem:val}).  Then  by counting the 
edges  from $u,v$ and another vertex of $X\cap V(G)$ to $V(G_3) \setminus X$, we get $e_{G_3}(X, V(G_3) \setminus X) \ge \Delta+2-2|X|+\delta-|X| >\Delta$, a contradiction. Thus $X$ is an independent set in $G_0$.  
Consequently, $G_3[X] \subseteq G_3-E(G_0)$. 
Since $L$ is bipartite with one bipartition as $\{v_1, \ldots, v_p\}$, by the construction of $G_3$,  
we know that the vertex set of the only possible odd cycle  of $G_3-E(G_0)$ contains $\{v_{p-1},  v_p,  v_{p+1}\}$. 
Thus $\{v_{p-1}, v_p, v_{p+1}\} \subseteq X$.   By Lemma~\ref{lem:graphical-biparite}(d)-(f), 
 the underlying simple graph of $L$  is a forest such that either $v_{p-1}$ and $v_{p+1}$ are in different components of the forest or $v_{p-1}v_{p+1}$ is an edge in the forest. It then follows that 
 the underlying simple graph of $G_3-E(G_0)-v_p$  is a forest by the construction of $G_3$. 
 As $X$ is an independent set in $G_0$ and $v_p$ is only adjacent to $v_{p-1}$ and $v_{p+1}$  in $G_3-E(G_0)$, we have $e_{G_3}(v_p, X\setminus \{v_{p-1}, v_p, v_{p+1}\}) =0$.  
Hence 
the multigraph  $R^*$ obtained from $G_3[X]$ by identifying $v_{p-1}$, $v_p$, and $v_{p+1}$ as a single vertex (deleting resulting loops) is bipartite. 
We claim that  $\Delta(R^*) \le \Delta$. It suffices to verity that  
$e_{G_3}(v_{p-1}, X\setminus \{v_{p-1}, v_p, v_{p+1}\}) +e_{G_3}(v_p, X\setminus \{v_{p-1}, v_p, v_{p+1}\}) +e_{G_3}(v_{p+1}, X\setminus \{v_{p-1}, v_p, v_{p+1}\}) \le \Delta$.  Recall that  $e_{G_3}(v_p, X\setminus \{v_{p-1}, v_p, v_{p+1}\}) =0$, and  by the assumption that $X$ is an independent set in $G_0$ and the evolution process from $G_0$ to $G_3$, we have 
\begin{eqnarray*}
	e_{G_3}(v_{p-1}, X\setminus \{v_{p-1}, v_p, v_{p+1}\}) & \le &d_L(v_{p-1}) =\df_{G_0}(v_{p-1}) -\df_{G_0}(v_{p}),  \quad \text{and} \\
	e_{G_3}(v_{p+1}, X\setminus \{v_{p-1}, v_p, v_{p+1}\}) & \le & d_L(v_{p+1}) \le \df_{G_0}(v_{p+1}) \le \df_{G_0}(v_p). 
\end{eqnarray*}
Thus   
\begin{eqnarray*}
	 && e_{G_3}(v_{p-1}, X\setminus \{v_{p-1}, v_p, v_{p+1}\}) +e_{G_3}(v_p, X\setminus \{v_{p-1}, v_p, v_{p+1}\})  \\
	 &&+e_{G_3}(v_{p+1}, X\setminus \{v_{p-1}, v_p, v_{p+1}\}) \\
	 &=& e_{G_3}(v_{p-1}, X\setminus \{v_{p-1}, v_p, v_{p+1}\})  +	e_{G_3}(v_{p+1}, X\setminus \{v_{p-1}, v_p, v_{p+1}\}) \\ 
	 &\le & d_L(v_{p-1})+d_L(v_{p+1})   
	 \le  \df_{G_0}(v_{p-1}) -\df_{G_0}(v_{p})+\df_{G_0}(v_{p+1})  \\
	 &\le &\df_{G_0}(v_{p-1}) -\df_{G_0}(v_{p}) + \df_{G_0}(v_{p}) = \df_{G_0}(v_{p-1})  \le \Delta, 
\end{eqnarray*}
and so $\Delta(R^*) \le \Delta$. 
Therefore by ~\eqref{claim:degree-of-v-p+1}, 
we get 
\begin{eqnarray*}
	e(G_3[X])  &\le&   e(R^*)+e(G_3[\{v_{p-1}, v_p, v_{p+1}\}])  \le \Delta \frac{|X|-3}{2} + e(G_3[\{v_{p-1}, v_p, v_{p+1}\}])  \\
	 &\le & \Delta \frac{|X|-3}{2}  +  \frac{1}{2}(3\Delta -e_{G_2}(\{v_{p-1}, v_p, v_{p+1}\}, V(G_2) \setminus \{v_{p-1}, v_p, v_{p+1}\})) \\
	 &\le & \Delta \frac{|X|-3}{2}  + \Delta \le \Delta \frac{|X|-1}{2}, 
\end{eqnarray*}
contradicting the assumption that $G_3[X]$ is $\Delta$-overfull. 
\qed 

We now have a $\Delta$-regular multigraph $G_3$ with no $\Delta$-overfull subgraph. Disregarding some multiple edges of $G_3$ joining $v_p$ and $v_{p+1}$, we will 
decompose the rest edges of $G_3$ into $\Delta$ edge-disjoint matchings in the rest of the proof. The rest edges still induce
a supergraph of $G_2$, and so we 
get a contradiction to $\chi'(G_2) \ge \Delta+1$, 
where we assume $\chi'(G_2) \ge \Delta+1$ by Claim~\ref{claim:G2-properties}(i).

\subsection{Decompose the edges of a subgraph of $G_3$}

Let $U^*=U\cup \{v_1\}$.  We summarize  some properties of $G_3$ as follows before proceed to decompose its edges. 
\begin{enumerate}[(i)]
	\item The multigraph  $G_3$ is  constructed based on $G_2$ by Procedure~\ref{Pro2:Edge-Addition}, $G_2$ is  obtained from $G_1$
	by Procedure~\ref{opr:vertex-contraction},  $G_1$ is the union of $G_0$ and the multigraph $L$, and $G_0=G$ if $n$ is even  and  $G_0$ is the union of $G$
	and an isolated vertex if $n$ is odd. 
	\item  We have $|V(G_3)|\le n+1$ but $|V(G_3)|> n-\ve n-2\eta n =(1-\ve-2\eta )n  \ge  (1-\ve -2\eta)n_0$. As 
	$\frac{1}{(1-\ve -2\eta)n_0} \le \frac{2}{n_0}$, we have $\frac{1}{(1-\ve -2\eta )n_0}  \le \min \left \{\frac{1}{n_1^*}, \frac{1}{n_2}, \frac{1}{n_3},   \frac{1}{n^*_4} \right \}  \ll  \eta  \ll  \tau^*, \ve$
	by our choice of  $n_0$  in~\eqref{eqn:parameters}. 
	\item By Claim~\ref{claim:G2-properties}(iii)-(iv), we have 
	\begin{eqnarray}
	d_{G_3}^s(v) &\ge&  \delta \quad \text{for $v\in U$};  \label{eqn:U-vertex-degree-in-G3}\\
	d_{G_3}^s(v) &\ge&  \Delta -\eta n-2\eta n=\Delta-3\eta n \quad \text{for $v\in V(G_3)\setminus U^*$.}\label{eqn:non-U-vertex-degree-in-G3}
	\end{eqnarray}
	\item For any vertex  $v\in V(G_3)\setminus U^*$,   we have 
	\begin{eqnarray}
		\mu_{G_3}(v)&<&3\eta n.   \label{eqn: vertex-multiplicity-not-from-U*}
	\end{eqnarray}
	This is because  $v\in V(G_3)\setminus U^*$, and  if any, we only identified vertices from $U^*$ in Procedure~\ref{opr:vertex-contraction}. Thus   $d_{G_3}^s(v) \ge d_G(v, V(G) \setminus U) \ge (\Delta -\eta n -|U|) > \Delta -3\eta n$. 
	Thus $\mu_{G_3}(v)<3\eta n$.

	\item For notation simplicity, we will assume $V(G_3)=V(G_1)$, but take care of the possible multiple 
	edges resulted  in the  Vertex Identification Procedure (note that $|V(G_3)| =|V(G_2)|$ and $|V(G_2)|$ is even as in the Vertex Identification Procedure, 
	we identified  three vertices into one at a time). 
	So we  have that  $V(G_3) =\{v_1, \ldots, v_m\}$, $m$ is even, and $d_{G_0}(v_1) \le d_{G_0}(v_2) \le \ldots \le d_{G_0}(v_m) =\Delta$. 
	Note that $G$ is not $\Delta$-regular and $U\ne \emptyset$ by $\delta<\frac{1}{2}(1+\ve)n$ from~\eqref{eqn:upper-bound-delta}. 
	Let 
	\begin{eqnarray*}
		g \in [1, m-1]  &&\text{be the index such that $d_{G_0}(v_g)<\Delta$ but $d_{G_0}(v_{g+1}) =\Delta$, and} \\
		h\in [1, g] && \text{ be the  index such that $v_h\in U^*$ but $v_{h+1} \not\in U^*$.}
	\end{eqnarray*}
	\item To streamline the proof, we assume without loss of generality that $v_{p+1} \in  U^*$. 
	This implies $h\ge p+1$ and so  $\{v_1, \ldots, v_{p+1}\} \subseteq U^*$.  
	This assumption is valid because the proof for $v_{p+1} \in  U^*$
	encompasses the case for  $v_{p+1}\not\in U^*$.  Here's why:
	\begin{itemize}
		\item Vertices in $\{v_1, \ldots, v_{p+1}\}$: when $h\ge p+1$, arguments concerning these vertices can be directly applied to just $\{v_1, \ldots, v_h\}$, because they are all within $U^*$. 
		\item Vertices in $\{v_{p+1}, \ldots, v_h\}$: these vertices exist only when $h\ge p+1$. 
		\item Vertices in $V(G_3)\setminus \{v_1, \ldots, v_{h}\}$: arguments for these vertices hold true regardless of whether $v_{p+1}$ is in $U^*$ or not. However,  when $p\ge h$ we just deal  vertices of $\{v_1, \ldots, v_{p+1}\}\setminus \{v_1, \ldots v_h\}$ the same way as we do for vertices of $V(G_3)\setminus \{v_1, \ldots, v_{h}\}$ when $h\ge p+1$. 
	\end{itemize}
	
		\item  Let $L^*=G_3-E(G_3^s)$, where $G_3^s$ is the underlying simple graph of $G_3$. 
	Under the assumption that $v_{p+1} \in  U^*$, we know that the  partite set  $\{v_1, \ldots, v_p\}$ of the bipartite 
	multigraph $L$ is a  proper subset of $U^*$.  Since vertices of $\{v_{2i}: i\in [1,p/2]\}$ have degree 0 in $L$, and  every vertex  $v\in V(G_3)\setminus\{v_1, \ldots, v_p\}$
	has at most one neighbor from $\{v_{2i}: i\in [1,p/2]\}$ in $L^*$, 
	we know that 
	\begin{equation}\label{eqn:neighbors-in-L*}
		|N_{L^*}(v)| \le \frac{1}{2}( |U^*|-1 )+1<\eta n+1 \quad \text{for any $v\in V(G_3)\setminus\{v_1, \ldots, v_p\}$}. 
	\end{equation}
\end{enumerate}

 We now fill in the details of the five steps outlined in Section 1. 
 
\begin{center}
	{\bf 	Step 1: Partition  $V(G_3)$ into two desired subsets $A$ and $B$}
\end{center}

 If $h \ge p+1$, then  we let  $u_{p+1}, \ldots, u_h$ be $h-p$ distinct vertices from $V_\Delta \setminus \{v_{g+1}\}$. 
Note that $u_{p+1}, \ldots, u_h$ exist as  $h-p<|U|<2\eta n$  but  $|V_\Delta|> (\frac{1}{2}-\frac{3\ve}{2} ) n$ by~\eqref{eqn:number-maximum-degree-vertex}. 
Let 
\begin{numcases}{ \ell^*=}
	\lceil \frac{1}{2} g\rceil &  \text{if $h \le p$;}  \nonumber  \\
	\frac{p}{2}+ h-p+  \lceil \frac{1}{2}(g -h) \rceil  & \text{if $h\ge p+1$.} \nonumber 
\end{numcases}
We then pair up some vertices of $G_3$ in order to get a desired partition of $V(G_3)$ into two subsets. 
If $h\le p$, for each $i\in [1,	\lceil \frac{1}{2} g\rceil ]$, let  $x_i=v_{2i-1}$ and $y_i=v_{2i}$.
If $h\ge p+1$, for each $i\in [1,\frac{p}{2}]$, let  $x_i=v_{2i-1}$ and $y_i=v_{2i}$; for each $i\in [1, h-p]$, let $x_{\frac{p}{2}+i}=v_{p+i}$ and $y_{\frac{p}{2}+i}=u_{p+i}$;  
for each $i\in [1, \lceil \frac{1}{2}(g-h) \rceil ]$, let $x_{\frac{p}{2}+h-p+i}=v_{h+2i-1}$ and $y_{\frac{p}{2}+h-p+i}= v_{h+2i}$. 
We let 
$$
N=\{x_i, y_i: i\in [1,\ell^*]\}.  
$$
We say that  $x_i$ and $y_i$ from the set $N$ are \emph{partners} of each other. 
Note that under our assumption that $h\ge p+1$, we know  $v_{2i-1}$
and $v_{2i}$ are  partners of each other for each $i\in [1,\frac{p}{2}]$, $u_{j}$ and $v_j$
are partners of each other for each $j\in [p+1,h]$.

Applying Lemma~\ref{lem:partition} on  $G^s_3$ and $N$,
we obtain a partition $\{A, B\}$ of $V(G_3)$ satisfying the following properties in  $G^s_3$:  $|A|=|B|$;  $|A\cap \{x_i,y_i\}|=1$ for each   $i\in [1,\ell^*]$; and 
$| d^s_{G_3}(v, A)-d^s_{G_3}(v, B)| \le n^{\frac{2}{3}}$ for each $v\in V(G_3)$. 
We  refine the partition through the procedure below.

\begin{Pro}[Partition Modification Procedure]\label{procedure:the partition modification}
	By moving  a vertex from its own part to its partner's part and vice versa, we  make the following assumptions. 
	\begin{enumerate}[(i)]
		\item $v_{2i-1}\in A$ and $v_{2i}\in B$ for each $i\in [1,\frac{p}{2}]$;
		\item 	$v_{p+1}, v_{p+2},  \ldots, v_h\in B$. 
	\end{enumerate}
	See Figure~\ref{f0} for an illustration of  the partners   in $A$ and $B$ for the case that $g-h$ is even. 
\end{Pro}

\begin{figure}[!htb]
	\begin{center}
		
		\begin{tikzpicture}[scale=0.8]

			\begin{scope}[shift={(-8,-8)}]
				\draw[rounded corners, fill=white!90!gray] (6, 0) rectangle (22, 2) {};
				
				\draw[rounded corners, fill=white!90!gray] (6, -4) rectangle (22, -2) {};
				
				{\tikzstyle{every node}=[draw ,circle,fill=white, minimum size=0.3cm,
					inner sep=0pt]
					\draw[black,thick](9,1) node (c)  {};
					
					\draw[black,thick](10.5,1) node (c1)  {};
					\draw[black,thick](12.5,1) node (c2)  {};
					\draw[black,thick](14,1) node (vp-1)  {};
					\draw[black,thick](9,-3) node (d)  {};
					\draw[black,thick](10.5,-3) node (d1)  {};
					\draw[black,thick](12.5,-3) node (d2)  {};
					\draw[black,thick](14,-3) node (vp)  {};
					\draw[black,thick](15.5,-3) node (vp+1)  {};
					\draw[black,thick](17.5,-3) node (vh)  {};
					\draw[black,thick](15.5,1) node (up+1)  {};
					\draw[black,thick](17.5,1) node (uh)  {};
					\draw[black,thick](19,1) node (vp+2)  {};
					\draw[black,thick](19,-3) node (vp+2)  {};
					\draw[black,thick](21,1) node (vl-1)  {};
					\draw[black,thick](21,-3) node (vl)  {};
				}

				{\tikzstyle{every node}=[draw ,circle,fill=black, minimum size=0.1cm,
					inner sep=0pt]
					\draw[black,thick](11,1) node (x1)  {};
					\draw[black,thick](11.5,1) node (x2)  {};
					\draw[black,thick](12,1) node (x3)  {};

				}
				
				\begin{scope}[shift={(0,-4)}]
					{\tikzstyle{every node}=[draw ,circle,fill=black, minimum size=0.1cm,
						inner sep=0pt]
						\draw[black,thick](11,1) node (x1)  {};
						\draw[black,thick](11.5,1) node (x2)  {};
						\draw[black,thick](12,1) node (x3)  {};

					}
				\end{scope}
				
				\begin{scope}[shift={(5,-4)}]
					{\tikzstyle{every node}=[draw ,circle,fill=black, minimum size=0.1cm,
						inner sep=0pt]
						\draw[black,thick](11,1) node (x1)  {};
						\draw[black,thick](11.5,1) node (x2)  {};
						\draw[black,thick](12,1) node (x3)  {};

					}
					
				\end{scope}
				
				\begin{scope}[shift={(5,0)}]
					{\tikzstyle{every node}=[draw ,circle,fill=black, minimum size=0.1cm,
						inner sep=0pt]
						\draw[black,thick](11,1) node (x1)  {};
						\draw[black,thick](11.5,1) node (x2)  {};
						\draw[black,thick](12,1) node (x3)  {};

					}
					
				\end{scope}

				\begin{scope}[shift={(8.5,0)}]
					{\tikzstyle{every node}=[draw ,circle,fill=black, minimum size=0.1cm,
						inner sep=0pt]
						\draw[black,thick](11,1) node (x1)  {};
						\draw[black,thick](11.5,1) node (x2)  {};
						\draw[black,thick](12,1) node (x3)  {};

					}
					
				\end{scope}
				
				\begin{scope}[shift={(8.5,-4)}]
					{\tikzstyle{every node}=[draw ,circle,fill=black, minimum size=0.1cm,
						inner sep=0pt]
						\draw[black,thick](11,1) node (x1)  {};
						\draw[black,thick](11.5,1) node (x2)  {};
						\draw[black,thick](12,1) node (x3)  {};

					}
					
				\end{scope}
				
				\node at (9,1.4) {$v_1$}; 
				\node at (9,-3.4) {$v_2$}; 
				\node at (10.5, 1.4) {$v_3$}; 
				\node at (10.5, -3.4) {$v_4$}; 
				\node at (12.5, 1.4) {$v_{p-3}$}; 
				\node at (12.5, -3.4) {$v_{p-2}$}; 
				\node at (14,1.4){$v_{p-1}$}; 
				\node at (14,-3.4) {$v_p$}; 
				\node at (15.5,-3.4) {$v_{p+1}$}; 
				\node at (17.5,-3.4) {$v_{h}$}; 
				\node at (15.5,1.4) {$u_{p+1}$}; 
				\node at (17.5,1.4) {$u_{h}$}; 
				\node at (19,1.4) {$a$}; 
				\node at (19,-3.4) {$b$}; 
				\node at (21,1.4) {$c$}; 
				\node at (21,-3.4) {$d$}; 
				
				\node at (7.2,1) {$A$};
				\node at (7.2,-3) {$B$};
				
			\end{scope}	
			
		\end{tikzpicture}
	\end{center}
	\caption{The location of the partner vertices  from $N$ in the partition $\{A,B\}$ when  $g-h$ is even, where $\{a,b\} =\{v_{h+1}, v_{h+2}\}$ and $\{c,d\} =\{ v_{g-1}, v_g\}$.}
	\label{f0}
\end{figure}

Still denote by $A$ and $B$ the 
resulting partition of $V(G_3)$.  
We  exchanged the locations of at most  $|U|+1 <2\eta n+1$ pairs of  partner vertices. 
As a result, the  partition $\{A,B\}$ in the multigraph $G_3$ 
has the properties below: 
	 \begin{eqnarray}
		|A| &=&|B|;  \nonumber \\
		|A\cap \{x_i,y_i\}| &=&1   \quad \text{for each   $i\in [1,\ell^*]$};  \nonumber\\
		|d^s_{G_3}(v, A)-d^s_{G_3}(v, B)| &\le& n^{\frac{2}{3}}+4\eta n +2\quad \text{for each $v\in V(G_3)$.  } \label{eqn:degree-diff}
	\end{eqnarray}
%where the term $2\eta n+2$ in the last line above is from the modification  of the original partition $\{A,B\}$.  

\begin{center}
	{\bf 	Step 2: Form multigraph $G_{A,B}$ and edge color it}
\end{center}

Let 
$$ G_A=G_3[A], \quad G_B=G_3[B], \quad \text{and } \quad H=G_3[A,B]. $$

  Define 
  \begin{eqnarray}
  	\Delta'&=& \left\lceil \frac{\Delta}{2}+5.3\eta n  \right\rceil \quad \text{and}\quad k=\Delta'+ \lceil\sqrt{\Delta'} \rceil.   \label{eqn:D'-and-k}
  	  \end{eqnarray}  	  
  	  If $\lceil \frac{1}{2}e_{G_3}(v_p, v_{p+1}) \rceil< 6\eta^{1/2} n$, then we  let $e_p=0$. 
  	  Otherwise, let 
  \begin{eqnarray}	  
  	e_p&=&\min\{k-\lfloor \frac{1}{2} e_{G_3}(v_{p-1}, v_{p}) \rfloor -\lfloor \frac{1}{2} e_{G_3}(v_{p-1}, v_{p+1}) \rfloor  -\lfloor \frac{1}{2} e_{G_3}(v_{p}, v_{p+1}) \rfloor, \lfloor 6\eta n\rfloor\}. \label{eqn:ep-definition}
  \end{eqnarray}
  	Since $e(G_3[\{v_{p-1}, v_p, v_{p+1}\}]) \le \Delta$ by Claim~\ref{claim:G1-no-overfullness}, it follows 
  	that $$k-\lfloor \frac{1}{2} e_{G_3}(v_{p-1}, v_{p}) \rfloor -\lfloor \frac{1}{2} e_{G_3}(v_{p-1}, v_{p+1}) \rfloor  -\lfloor \frac{1}{2} e_{G_3}(v_{p}, v_{p+1}) \rfloor>5.3 \eta n.$$ 
  	Thus $e_p>5.3\eta n$ if $e_p>0$.

Since  $u_{p+1} \in V_\Delta$, we have 
$d^s_{G_2}(u_{p+1}) \ge \Delta-2\eta n$ by Claim~\ref{claim:G2-properties}(iii).  
Thus   \eqref{eqn:degree-diff} implies that there exist two disjoint edge sets $E_1$ and $E_2$ such that 
\begin{itemize}
	\item $E_1 \subseteq E_{G_A}(u_{p+1}, A\setminus \{u_{p+1}\})$ and $E_2 \subseteq E_{H}(u_{p+1}, B \setminus U\})$, and 
	\item $|E_1|=|E_2|=\max\{e_{G_3}(v_{p}, v_{p+1})-\left \lfloor\frac{1}{2} e_{G_3}(v_{p}, v_{p+1}) \right \rfloor-e_p, 0\}$. 
\end{itemize}
Furthermore, we partition the edges in $E_{G_3}(v_p, v_{p+1})$ into two subsets as follows:
\begin{itemize}
	\item $F_1$ is a subset of edges from $E_{G_3}(v_p, v_{p+1})$ of size the same as $|E_1|$; 
	\item $F_2=E_{G_3}(v_p, v_{p+1})\setminus F_1$. 
\end{itemize}
Since 
$e_{G_3}(v_p, v_{p+1})\le \frac{1}{2}\df_{G_2}(v_{p+1}) \le \frac{\Delta}{3}$ by Claim~\ref{claim:G2-properties}(v), 
it follows that 
\begin{equation}\label{eqn:E2-size}
|E_1| <\frac{\Delta}{6}-5.3\eta n. 
\end{equation}
As $d^s_{G_3}(v_p), d^s_{G_3}(v_{p+1}) \ge \delta$ and $e_p \le 6\eta n$, by~\eqref{eqn:degree-diff}, there exist vertex-disjoint 
sets $F_{21}$ and $F_{22}$ 
such that 
\begin{itemize}
	\item $F_{21} \subseteq E_{G_B}(v_p, B\setminus U^*)$ and $F_{22} \subseteq E_{G_B}(v_{p+1}, B\setminus U^*)$,  
	\item $|F_{21}|=|F_{22}| =e_p$,  and 
	\item both $G_B[F_{21}]$ and $G_B[F_{22}]$ are simple. 
\end{itemize}

We can ``ignore''  the edges in $F_1$ for now, as they were all added to $G_2$ when forming $G_3$. This allows us to maintain the inequality $\chi'(G)\le \chi'(G_2) \le \chi'(G_3-F_1)$. Ignoring these edges also simplifies $G_3$ (reduces the multiplicity of $G_B$), which is helpful for the 1-factor extension in Step 4.
We will temporarily disregard the edges in $E_1$ until Step 5, where we will reintroduce them to form a nearly-bipartite multigraph and apply Theorem~\ref{thm:chromatic-index-nearly-bipartite-graph}. 
The other reason for ignoring edges  in  $F_1$ during Steps 2-4 is to ensure that sets $A$ and $B$ have the same number of remaining edges.  Edges in $F_{21}\cup F_{22}$ will be removed in Steps 2-3  as to guarantee that 
the degrees of $v_p$ and $v_{p+1}$ in the multigraph $G_{A,B}$ formed below do not 
exceed $\Delta'$. 
These edges  of  $F_{21}\cup F_{22}$ will be colored  in Step 4. 
 We can now construct  $G_{A,B}$ by the following procedure.

\begin{Pro}[$G_{A,B}$ Formation Procedure]\label{pro:G-A-B}
	We  define $G_{A,B}$ initially  as an empty  spanning subgraph of $G_3$, and then  adding the following edges.
\begin{enumerate}[(i)]
	\item $A$-edges:  $(E(G_A)\setminus E_1)$;
	\item $B$-edges:  $\left(E(G_B) \setminus (E_{G_3}(v_{p-1}, v_{p}) \cup F_{21} \cup F_{22})  \right)$ and  all edges of $F_2$;
	\item $AB$-edges: 
	\begin{itemize}
			\item Edges of  $E_2$;
		\item For each    vertex  $ v\in \{v_{2i-1}: i\in [1, \frac{p}{2}]\}$, 
		we add    edges    from $E_{L^*}(v,  B\cap U^*)$ to $G_{A,B}$. 
		Precisely, 
		for each $u\in N_{L^*}(v)\cap (B\cap U^*)$, we add  either $ \lfloor  \frac{1}{2}(e_{L^*}(u, v) )\rfloor$ or $ \lceil  \frac{1}{2}(e_{L^*}(u, v)) \rceil$ 
		edges from  $E_{L^*}(u, v)$ to $G_{A,B}$ such that 
		the total number of edges from $E_{L^*}(v,  B\cap U^*)$ added to $G_{A,B}$ is $\lfloor \frac{1}{2} e_{L^*}(v, B\cap U^*) \rfloor$.
	\end{itemize}
\end{enumerate}
\end{Pro}

\begin{OBV}\label{obv:G-AB}
We have the following observations. \begin{enumerate}[(i)]
	\item Since $L^*=G_3-E(G_3^s)$, for  any $ v\in \{v_{2i-1}: i\in [1, \frac{p}{2}]\}$ and any $u\in N_{G_3}(v, B)$, we have 
	$e_{L^*}(v,u) =e_{G_3}(v,u)-1$. Thus,  by (iii) of Procedure~\ref{pro:G-A-B}, we have $e_{G_3-E(G_{A,B})}(v,u) \ge \frac{1}{2}e_{L^*}(v,u)$. 
	\item As $G_3$ is regular and $|A|=|B|$, we have $e(G_A)=e(G_B)$.   
	Thus by the construction of $G_{A,B}$ and the fact that $|E_1|=|F_1|$,   we have $e(G_{A,B}[A])=e(G_{A,B}[B])+|F_{21}|+|F_{22}|$. 
\end{enumerate}
\end{OBV}

Next, we study the degrees of vertices in $G_{A,B}$. 
\begin{CLA}\label{claim: degrees-in-G-A-B}
	The following statements hold. 
	\begin{enumerate}[(i)]
		\item For any $v\in V(G_3)\setminus U^*$,  we have $$ \frac{1}{2}\Delta-3.6\eta n \le d_{G_{A,B}} (v)\le \frac{1}{2}\Delta+5.1\eta n.$$
		\item  For any $v\in  \left(\{v_{2i}:  i\in [1, \frac{p}{2}-1]\} \cup \{v_i: i\in [p+2, h]\} \right)$, we have $$ \frac{1}{2}\Delta-3.1\eta n \le d_{G_{A,B}} (v)\le \frac{1}{2}\Delta+2.1\eta n.$$  
		\item For any $v\in \{v_1, v_3,  \ldots,  v_{p-3}, v_{p-1}\}$, we have $$ \frac{1}{2}\Delta-3.1\eta n \le d_{G_{A,B}} (v)\le \frac{1}{2}\Delta+3.1\eta n.$$  
		\item  For any   $v\in \{v_p, v_{p+1}\}$,  we have 
	 $$ \frac{1}{2}\Delta-3.1\eta n \le d_{G_{A,B}} (v)\le \frac{1}{2}\Delta+2.1\eta n.$$ 
	\end{enumerate}
	
\end{CLA}

\pf    Recall that $L^*=G_3-E(G_3^s)$ where $G_3^s$ is the underlying simple graph of $G_3$. 
Since  $G_3[E_1\cup E_2]$ is a simple graph by the choice of $E_1$ and $E_2$ (the multiple edges between vertices of $V_\Delta$ and $V(G_3)\setminus V_\Delta$ are incident with $v_{p+1}$ by the Vertex Identification Procedure), and 
$F_{21}$ and $F_{22}$ are vertex-disjoint such that $G_3[F_{21}\cup F_{22}]$ is simple, 
by the  property  of  the partition $\{A,B\}$ and  Procedure~\ref{pro:G-A-B}, the $G_{A,B}$ Formation Procedure,  for any $v\in V(G_{A,B}) \setminus \{v_p, v_{p+1}\}$, we have 
\begin{equation}\label{eqn:degree-in-G-A-B}
	\frac{1}{2}(d^s_{G_3}(v)-(n^{\frac{2}{3}}+4\eta n +2))-1  \le  d^s_{G_{A,B}} (v)  \le  \frac{1}{2}(d^s_{G_3}(v)+(n^{\frac{2}{3}}+4\eta n +2))+1. 
\end{equation}

For (i),  let $v\in V(G_3)\setminus U^*$.  Then  $d_{G_3}^s(v) \ge 
\Delta -3\eta n$ by~\eqref{eqn:non-U-vertex-degree-in-G3}.   
Thus by~\eqref{eqn:degree-in-G-A-B} and the fact that   $d^s_{G_{A,B}} (v) \le d_{G_{A,B}} (v) \le d^s_{G_{A,B}} (v) +d_{L^*}(v) \le d^s_{G_{A,B}} (v) +3\eta n$ ($d_{L^*}(v)\le 3\eta n$ by~\eqref{eqn:non-U-vertex-degree-in-G3}), we get 
\begin{eqnarray*}
\frac{1}{2}\Delta-3.6\eta n=\frac{1}{2}(\Delta -7.2\eta n)  \le & d_{G_{A,B}} (v) &\le  \frac{1}{2}( \Delta+(n^{\frac{2}{3}}+4\eta n +4))+3\eta n = \frac{1}{2}\Delta+5.1\eta n. 
\end{eqnarray*}
The same conclusion as above holds even if $v=u_{p+1}$, as $|E_1|=|E_2|$.

For (ii), let $ v\in \left(\{v_{2i}:  i\in [1, \frac{p}{2}-1]\} \cup \{v_i: i\in [p+2, h]\} \right)$.   By (i) and (ii) of the Partition Modification Procedure,  we have  $v\in B$
and $N_{L^*}(v)\cap U^* \subseteq A$.    By    (iii) of the $G_{A,B}$ Formation Procedure, we have 
$$\frac{1}{2} d_{L^*}(v) -\frac{1}{2}|U^*| \le  \frac{1}{2} d_{L^*}(v) -|N_{L^*}(v)\cap U^*|   \ \le  d_{G_{A,B}}(v, A)  \le  \frac{1}{2} d_{L^*}(v).$$ 
Since $d_{G_{A,B}} (v) =d_{G_3}(v,B)+d_{G_{A,B}}(v, A) =d^s_{G_3}(v,B)+d_{G_{A,B}}(v, A) $, it follows that 
\begin{eqnarray*}
	 d_{G_{A,B}} (v) &\le &  \frac{1}{2}(d^s_{G_3}(v)+(n^{\frac{2}{3}}+4\eta n +4)) +\frac{1}{2} d_{L^*}(v)\\
	  &\le & \frac{1}{2}(d^s_{G_3}(v)+ d_{L^*}(v))+2.1\eta n  \\
	  &= &  \frac{1}{2} \Delta+2.1\eta n, 
\end{eqnarray*} 
and that 
\begin{eqnarray*}
	d_{G_{A,B}} (v) & \ge &  \frac{1}{2}(d^s_{G_3}(v)-(n^{\frac{2}{3}}+4\eta n +4)) +\frac{1}{2} d_{L^*}(v)-\frac{1}{2}|U^*|\\
	&\ge  & \frac{1}{2}(d^s_{G_3}(v)+ d_{L^*}(v))-3.1\eta n  \\
	&=&  \frac{1}{2} \Delta-3.1\eta n. 
\end{eqnarray*}

For (iii), let $v\in \{v_1, v_3,  \ldots,  v_{p-3}, v_{p-1}\} $.  By Procedure~\ref{procedure:the partition modification}, the Partition Modification Procedure, we have $v\in A$ and $N_{L^*}(v)\cap U^* \subseteq B$.   We first claim that 
$$\frac{1}{2} d_{L^*}(v) -\eta n-1  \ \le  d_{G_{A,B}}(v, B)  \le  \frac{1}{2} d_{L^*}(v)+\eta n.$$ 
Suppose,   $v_{s_1}, \ldots, v_{s_t}$  are all the neighbors of $v$ in $L^*$
from $V(L^*)\setminus U^*$, where  $t \ge 1$ is an integer, 
 $s_1 < s_2 <\ldots < s_t$,  and $s_1,  s_2, \ldots, s_t\in [h+1, g]$.  By Lemma~\ref{lem:graphical-biparite}(e) and the definition of $h$, we know that  
 $s_1, s_2, \ldots, s_t$ are consecutive integers satisfying 
 $e_L(v,v_{s_t}) \le \ldots \le e_L(v,v_{s_2})$ and  $e_L(v,v_{s_1}), e_L(v,v_{s_2})<\eta n$. By the definition of $N$ and the construction of $\{A,B\}$, we know that 
$$
 \left ||A\cap \{v_{s_1}, \ldots, v_{s_t}\}| -|B\cap \{v_{s_1}, \ldots, v_{s_t}\}|  \right| \le 2. 
$$
Thus,  using $e_L(v,v_{s_t}) \le \ldots \le e_L(v,v_{s_2})<\eta n$ and  $e_L(v,v_{s_1})<\eta n$, we get 
$$
|d_{L}(v, A\cap \{v_{s_1}, \ldots, v_{s_t}\}) - d_L(v, B\cap \{v_{s_1}, \ldots, v_{s_t}\}) | <2\eta n. 
$$
Therefore,  by (i) of  the $G_{A,B}$ Formation Procedure,  and the fact that two consecutive vertices in $v_{s_1}, \ldots, v_{s_t}$ are separated in different parts  of $\{A,B\}$, we get 
$$
d_{G_{A,B}}(v, A\cap \{v_{s_1}, \ldots, v_{s_t}\})  
 \le \frac{1}{2} \left(d_{L}(v, \{v_{s_1}, \ldots, v_{s_t}\}) + 2 \eta n \right) \quad \text{and}$$
$$
 d_{G_{A,B}}(v, A\cap \{v_{s_1}, \ldots, v_{s_t}\})  
 \ge \frac{1}{2}\left(d_{L}(v, \{v_{s_1}, \ldots, v_{s_t}\}) -2 \eta n\right).$$
 Then by (iii) of the $G_{A,B}$ Formation Procedure, we get 
 \begin{eqnarray*}
 &&d_{G_{A,B}} (v, B\cup \{v_{s_1}, \ldots, v_{s_t} \})  \\
 &=& d_{G_{A,B}} (v, B\cap U^*) +d_{G_{A,B}}(v, A\cap \{v_{s_1}, \ldots, v_{s_t} \}) \\  
   & \le & d_{G_{A,B}} (v, B\cap U^*) +\frac{1}{2} \left(d_{L}(v, \{v_{s_1}, \ldots, v_{s_t}\}) + 2 \eta n \right) \\ 
   & \le & \left \lfloor\frac{1}{2} e_{L^*}(v, B\cap U^*)  \right \rfloor+\frac{1}{2} \left(d_{L^*}(v, \{v_{s_1}, \ldots, v_{s_t}\}) + |N_{G_3}(v)\cap \{v_{s_1}, \ldots, v_{s_t}\}|+2 \eta n \right)  \\
   &\le & \frac{1}{2} d_{L^*}(v)+ \frac{1}{2}|N_{G_3}(v)\cap \{v_{s_1}, \ldots, v_{s_t}\}|+\eta n, 
 \end{eqnarray*}
 and 
  \begin{eqnarray*}
 	&&d_{G_{A,B}} (v, B\cup \{v_{s_1}, \ldots, v_{s_t} \})  \\
 	&= & d_{G_{A,B}} (v, B\cap U^*) +d_{G_{A,B}}(v, A\cap \{v_{s_1}, \ldots, v_{s_t} \}) \\  
 	& \ge & d_{G_{A,B}} (v, B\cap U^*) +\frac{1}{2} \left(d_{L}(v, \{v_{s_1}, \ldots, v_{s_t}\}) -2 \eta n \right) \\ 
 	& \ge & \left \lfloor \frac{1}{2} e_{L^*}(v, B\cap U^*)  \right \rfloor +\frac{1}{2} \left(d_{L^*}(v, \{v_{s_1}, \ldots, v_{s_t}\}) + |N_{G_3}(v)\cap \{v_{s_1}, \ldots, v_{s_t}\}|- 2 \eta n \right)  \\
 	& \ge & \frac{1}{2} d_{L^*}(v)+ \frac{1}{2}|N_{G_3}(v)\cap \{v_{s_1}, \ldots, v_{s_t}\}|-\eta n-1. 
 \end{eqnarray*}
 
Since  $d_{G_{A,B}}(v)=d_{G^s_{A,B}}(v, A\setminus  \{v_{s_1}, \ldots, v_{s_t}\} )+d_{G_{A,B}}(v,  B\cup \{v_{s_1}, \ldots, v_{s_t}\} )$, 
$d_{G^s_{A,B}}(v, A)-1 \le d_{G^s_{A,B}}(v, A\setminus  \{v_{s_1}, \ldots, v_{s_t}\} )+\frac{1}{2}|N_{G_3}(v)\cap \{v_{s_1}, \ldots, v_{s_t}\}| \le d_{G^s_{A,B}}(v, A)+1$ (two consecutive vertices in $v_{s_1}, \ldots, v_{s_t}$ are separated in different parts  of $\{A,B\}$), 
and by~\eqref{eqn:degree-diff},  we get 
\begin{eqnarray*}
d_{G_{A,B}}(v) & \ge& \frac{1}{2}(d^s_{G_3}(v)-(n^{\frac{2}{3}}+4\eta n +4)) -1 +\frac{1}{2}d_{L^*}(v)-\eta n -1\\
 &=& \frac{1}{2}(d^s_{G_3}(v)+d_{L^*}(v))-3.1\eta n=\frac{1}{2}\Delta-3.1\eta n,  
\end{eqnarray*}
and 
\begin{eqnarray*}
	d_{G_{A,B}}(v) & \ \le & \frac{1}{2}(d^s_{G_3}(v)+(n^{\frac{2}{3}}+4\eta n +4)) +1 +\frac{1}{2}d_{L^*}(v)+\eta n\\
	&<& \frac{1}{2}\Delta+3.1\eta n. 
\end{eqnarray*}

For (iv), let $v\in \{v_{p}, v_{p+1}\}$.   
Note that $d_{G_{A,B}} (v) =d^s_{G_{A,B}} (v, B)-e_p+d_{ G_{A,B}\cap L^*}(v, A)+d_{G_{A,B} \cap L^*}(v,B)$, $d^s_{G_{A,B}} (v, B)=d^s_{G_3}(v, B)$, 
and $  \frac{1}{2}(d^s_{G_3}(v)+(n^{\frac{2}{3}}+4\eta n +2)) \ge  d^s_{G_3}(v, B) \ge \frac{1}{2}(d^s_{G_3}(v)-(n^{\frac{2}{3}}+4\eta n +2)) $ by~\eqref{eqn:degree-diff}. 
 By  the Partition Modification Procedure and the $G_{A,B}$ Formation Procedure, we have $v\in B$, 
$N_{L^*}(v) \subseteq A \cup \{v_p, v_{p+1}\}$, and 
 $$d_{G_{A,B}\cap L^*}(v, A) \le  \frac{1}{2} d_{L^*}(v) -\left\lfloor \frac{1}{2} e_{G_3}(v_{p}, v_{p+1}) \right\rfloor$$ and 
$$ d_{G_{A,B}\cap L^*}(v, A) \ge \sum_{u\in N_{L^*}(v)\setminus \{v_p, v_{p+1}\}} \left  \lfloor \frac{1}{2} e_{L^*}(v,u)  \right \rfloor  \ge  \frac{1}{2} d_{L^*}(v) -\left \lceil \frac{1}{2} e_{G_3}(v_{p}, v_{p+1}) \right \rceil  -|N_{L^*}(v)|.$$  
By (ii) of 
the  $G_{A,B}$ Formation Procedure, we have 
$ d_{G_{A,B} \cap L^*}(v,B)= |F_2|$.  As $|F_2|=\left \lfloor\frac{1}{2} e_{G_3}(v_{p}, v_{p+1}) \right \rfloor+e_p$, we get  
\begin{eqnarray*}
	d_{G_{A,B}} (v) & \ge &  \frac{1}{2}(d^s_{G_3}(v)-(n^{\frac{2}{3}}+4\eta n +4))-e_p +\frac{1}{2} d_{L^*}(v)-\left \lceil \frac{1}{2} e_{G_3}(v_{p}, v_{p+1}) \right \rceil \\
	&&-|N_{L^*}(v)|+|F_2|\\
	& \ge &  \frac{1}{2} ( d^s_{G_3}(v)+d_{L^*}(v)) -3.1\eta n  \quad (\text{by~\eqref{eqn:neighbors-in-L*} for $v_{p+1}$ and the fact that $|N_{L^*}(v_p)| \le 2$})\\
	&=&\frac{1}{2} \Delta-3.1\eta n. 
\end{eqnarray*} 
On the other hand, 
\begin{eqnarray*}
	d_{G_{A,B}} (v) &  \le  &  \frac{1}{2}(d^s_{G_3}(v)+(n^{\frac{2}{3}}+4\eta n +4))-e_p+ \frac{1}{2} d_{L^*}(v)-\left\lfloor \frac{1}{2} e_{G_3}(v_{p}, v_{p+1})\right\rfloor+ |F_2| \\
	& \le &  \frac{1}{2} ( d^s_{G_3}(v)+d_{L^*}(v)) +2.1\eta n =\frac{1}{2} \Delta+2.1\eta n. 
\end{eqnarray*} 
\qed

Next we show that $G_{A,B}$ contains no $k$-overfull subgraph, which will be used to bound up 
the chromatic index of $G_{A,B}$. 

\begin{CLA}\label{claim:G-AB-no-overfull-subgrap}
The  multigraph $G_{A,B}$ contains no $k$-overfull subgraph. 
\end{CLA}

\pf  By  Claim~\ref{claim: degrees-in-G-A-B}, we have $\Delta(G_{A,B}) \le \Delta'<k=\Delta'+\lceil \sqrt{\Delta'} \rceil$.
We suppose to the contrary that there exists $X \subseteq V(G_{A,B})$ such that 
$G_{A,B}[X]$ is $k$-overfull or $ 2e_{G_{A,B}}(X) \ge k (|X|-1)+2$.  As $\Delta(G_{A,B}) \le \Delta'$ and so 
$ 2e(G_{A,B}[X])\le \Delta'|X|$, 
it follows that $|X| \le \sqrt{\Delta'}+1$.  
Since $d^s_{G_{A,B}}(v) \ge  \frac{1}{2} \left(\delta-(n^{\frac{2}{3}}+4\eta n +2) \right)\ge \frac{1}{2}(\frac{1}{7}(1-\frac{4}{14})n-(n^{\frac{2}{3}}+4\eta n +2) )> \frac{4.5}{98} n> \frac{1}{22}n$ 
for any $v\in V(G_{A,B}) \setminus \{v_1\}$ and $e_{G_{A,B}}(X, V(G_{A,B})\setminus X) <k<n$, 
we get 
$$n>e_{G_{A,B}}(X, V(G_{A,B})\setminus X) \ge \frac{1}{22}n \left (|X|-1 \right)-{|X| \choose 2}. $$
Since $|X| \le \sqrt{\Delta'}+1$ and so the function  $f(|X|)=\frac{1}{22}n (|X|-1)-{|X| \choose 2} -n$  is increasing in $|X|$, and  $f(24) >0$,  we conclude that $\frac{1}{22}n (|X|-1)-{|X| \choose 2} >n$ when $|X| \ge 24$. 
Thus  $|X|  \le 23$. 

We first show that $ e(G_{A,B}[\{v_{p-1}, v_p, v_{p+1}\}])\le k$. 
By the definition of $G_{A,B}$ and the definition of $e_p$ from~\eqref{eqn:ep-definition}, we have 
\begin{eqnarray*}
	&& e(G_{A,B}[\{v_{p-1}, v_p, v_{p+1}\}])  = \left\lfloor  \frac{1}{2}e_{G_3}(v_{p-1}, v_p)\right \rfloor+\left\lfloor \frac{1}{2}e_{G_3}(v_{p-1}, v_{p+1})\right\rfloor +\left \lfloor \frac{1}{2} e_{G_3}(v_{p}, v_{p+1}) \right \rfloor+e_p \\
	&\le &  \max\left\{k, \frac{1}{2}e_{G_3}(v_{p-1}, v_p)+\frac{1}{2}e_{G_3}(v_{p-1}, v_{p+1}) +\frac{1}{2}e_{G_3}(v_{p}, v_{p+1})\right\}\\
	&=& \max\left\{k, \frac{1}{2} \left( \frac{1}{2} \left(3\Delta -e_{G_2}(\{v_{p-1}, v_p, v_{p+1}\}, V(G_2) \setminus \{v_{p-1}, v_p, v_{p+1}\})\right)\right)\right\} \\
	&\le & \max\left\{k,\frac{1}{4} (3\Delta -\Delta)\right\}  \quad \text{(by ~\eqref{claim:degree-of-v-p+1})} \\ 
	&=& k. 
\end{eqnarray*}
Assume first that $|X|=3$. If $X\ne \{v_{p-1}, v_p, v_{p+1}\}$, then by the construction of $G_3$ and $G_{A,B}$, we know that the number of 
edges among any three vertices of $G_{A,B}$ is at most $\frac{1}{2} (\Delta-(\delta-1))+2 <k$ (at most $\Delta-(\delta-1)$ edges joining two of the vertices  and 2 edges joining the third vertex to those two vertices)  or $\frac{1}{2}d_{G_3}(v_{p+1})+2<k$ (when two of the vertices are $v_p$ and $v_{p+1}$ and the third vertex is not $v_{p-1}$). 
If
$X=\{v_{p-1}, v_p, v_{p+1}\}$, then by the definition of $G_{A,B}$ and the calculations above, we have 
$
e(G_{A,B}[X])   = e(G_{A,B}[\{v_{p-1}, v_p, v_{p+1}\}])  \le k, 
$
a contradiction.

Thus we assume $|X| \ge 5$.
If  $G_{A,B}[X] \cap L^*$ is bipartite, then we have  
$\chi'(G_{A,B}[X] ) \le \chi'(G_{A,B}[X] \cap L^*)+{|X|\choose 2} \le \Delta'+ {23 \choose 2}<k$, a contradiction. Thus $G_{A,B}[X] \cap L^*$ contains an odd cycle.  As the only odd cycle in $L^*$ is on $\{v_{p-1}, v_p, v_{p+1}\}$, we have $\{v_{p-1}, v_p, v_{p+1}\} \subseteq X$. As the only  possible odd cycle  of $G_3-E(G_0)$ has its vertex set as $\{v_{p-1}, v_p, v_{p+1}\}$, 
and   the underlying simple graph of $L$  is a forest by Lemma~\ref{lem:graphical-biparite}(d),  it follows that 
 the underlying simple graph of $L^*-v_p$  is a forest by the construction of $G_3$. 
By the construction of $G_3$, we have $e_{L^*}(v_p, X\setminus \{v_{p-1}, v_p, v_{p+1}\}) =0$.  
Hence 
the multigraph  $R^*$ obtained from $G_{A,B}[X]\cap L^*$ by identifying $v_{p-1}$, $v_p$, and $v_{p+1}$ as a single vertex is bipartite. 
We claim that  $\Delta(R^*) \le \frac{1}{2}\Delta$. It suffices to verity that  
$e_{G_{A,B} \cap L^*}(v_{p-1}, X\setminus \{v_{p-1}, v_p, v_{p+1}\}) +e_{G_{A,B} \cap L^*}(v_p, X\setminus \{v_{p-1}, v_p, v_{p+1}\}) +e_{G_{A,B} \cap L^*}(v_{p+1}, X\setminus \{v_{p-1}, v_p, v_{p+1}\}) \le \frac{1}{2} \Delta$.  Recall that  $e_{G_{A,B} \cap L^*}(v_p, X\setminus \{v_{p-1}, v_p, v_{p+1}\}) =0$, and  by the evolution process from $G_0$ to $G_3$ and Procedure~\ref{pro:G-A-B}(ii)-(iii), we have 
\begin{eqnarray*}
	e_{G_{A,B} \cap L^*}(v_{p-1}, X\setminus \{v_{p-1}, v_p, v_{p+1}\}) & \le & \frac{1}{2}d_L(v_{p-1}) = \frac{1}{2}(\df_{G_0}(v_{p-1}) -\df_{G_0}(v_{p})),  \quad \text{and} \\
	e_{G_{A,B} \cap L^*}(v_{p+1}, X\setminus \{v_{p-1}, v_p, v_{p+1}\}) & \le & \frac{1}{2}d_L(v_{p+1}) \le \frac{1}{2}\df_{G_0}(v_{p+1}) \le \frac{1}{2}\df_{G_0}(v_p). 
\end{eqnarray*}
Thus   
\begin{eqnarray*}
	&& e_{G_{A,B} \cap L^*}(v_{p-1}, X\setminus \{v_{p-1}, v_p, v_{p+1}\}) +e_{G_{A,B} \cap L^*}(v_p, X\setminus \{v_{p-1}, v_p, v_{p+1}\})  \\
	&&+e_{G_{A,B} \cap L^*}(v_{p+1}, X\setminus \{v_{p-1}, v_p, v_{p+1}\}) \\
	&=& e_{G_{A,B} \cap L^*}(v_{p-1}, X\setminus \{v_{p-1}, v_p, v_{p+1}\})  +	e_{G_{A,B} \cap L^*}(v_{p+1}, X\setminus \{v_{p-1}, v_p, v_{p+1}\}) \\ 
	& \le & \frac{1}{2}d_L(v_{p-1})+\frac{1}{2}d_L(v_{p+1})   =\frac{1}{2}(\df_{G_0}(v_{p-1}) -\df_{G_0}(v_{p}))+\frac{1}{2}d_L(v_{p+1})\\
	&\le & \frac{1}{2}(\df_{G_0}(v_{p-1}) -\df_{G_0}(v_{p}))+\frac{1}{2}\df_{G_0}(v_{p+1})  \\
	&\le & \frac{1}{2}  \left (\df_{G_0}(v_{p-1}) -\df_{G_0}(v_{p}) + \df_{G_0}(v_{p})\right) = \frac{1}{2}\df_{G_0}(v_{p-1})  \le \frac{1}{2}\Delta, 
\end{eqnarray*}
and so $\Delta(R^*) \le \frac{1}{2}\Delta$. 
Therefore, 
\begin{eqnarray*}
	e(G_{A,B}[X])  &\le&   e(R^*)+e(G_{A,B}[\{v_{p-1}, v_p, v_{p+1}\}]) \\
	& \le& \frac{1}{2} \Delta\frac{|X|-3}{2} + e(G_{A,B}[\{v_{p-1}, v_p, v_{p+1}\}]) +{|X| \choose 2}  \\
	&\le & \frac{1}{2} \Delta\frac{|X|-3}{2}  +  k+{23 \choose 2}
	<   k\frac{|X|-1}{2} , 
\end{eqnarray*}
contradicting the assumption that $G_{A,B}[X]$ is $k$-overfull. 
\qed

As $\Delta(G_{A,B} )\le \Delta'$ by Claim~\ref{claim: degrees-in-G-A-B} and $G_{A,B}$
contains no $k$-overfull subgraph by Claim~\ref{claim:G-AB-no-overfull-subgrap}, 
Theorem~\ref{thm:chromatic-index-bound} implies that 
$G_{A,B}$ has an edge coloring using $k$ colors. By Theorem~\ref{lem:equa-edge-coloring}, 
we find an equalized edge $k$-coloring  $\varphi_0$ of  $G_{A,B}$.  In particular, $\varphi_0$
satisfies the following properties: 
each color is missed by less than  $9\eta n-2$ vertices of $G_{A,B}$ or 
\begin{eqnarray}\label{s1.1}
	|\pbar_0^{-1}(i)| &<& 9\eta n-2 \quad \text{for each $i\in [1,k]$}. 
\end{eqnarray}
The reason that~\eqref{s1.1} holds is from the following argument: 
by Claim~\ref{claim: degrees-in-G-A-B}, we have $k-d_{G_{A,B}}(v) <8.92\eta n <9 \eta n-3$.   Thus 
\begin{eqnarray} \label{eqn:total-missing-color-a-d}
	\sum\limits_{v\in G_{A,B}}|\pbar_0(v)| &< &  (9 \eta n-3) \times n=9 \eta n^2-3n. 
\end{eqnarray}
As $\varphi_0$
is an equalized edge  $k$-coloring of $G_{A,B}$, we have $ \left||\pbar_0^{-1}(i)| -|\pbar_0^{-1}(j)| \right| \le 2$ for any $i,j\in [1,k]$. Therefore, \eqref{s1.1} holds.  
This  partial edge coloring  $\varphi_0$  of $G_3$ will be  extended  in the subsequent  steps.

\begin{center}
	{\bf Step 3: Extending  the $k$ color classes from Step 2   into  1-factors}
\end{center}
We modify the partial edge coloring of $G_3$ obtained in Step 2 by exchanging  alternating paths (swapping ``uncolor'' and a given color ``$i$'' on the edges of a path $P$). Upon the completion of Step 3, each of the $k$ color class will be a 1-factor of $G_3$. In the process of Step 3,  a few edges of $H-E(G_{A,B})$ will be colored and  a few edges of $G_A \cup G_B$  will be uncolored. Denote by  $R_A$  the subgraph of $G_A-E_1$ consisting of the uncolored edges, which will be   empty initially. 
Denote by $R_B$ the union of the multigraph formed by edges of $F_{21}\cup F_{22}$ 
and the subgraph of $G_B-(F_1 \cup F_{21}\cup F_{22})$ consisting of the uncolored edges.  
One up to  four edges will be added to at least one of  $R_A$ and $R_B$
when each time we exchange colors on an alternating path. 
The conditions below will be satisfied at the completion of this step.
\begin{enumerate}[S3.1]
	\item  The total number of uncolored edges produced in Step 3 in each of $R_A$
	and $R_B$ is less than $18 \eta n^2$. 
	Furthermore,  $R_A$
	and $R_B$ have the same number of uncolored edges. 	\label{s31}
	\item  $\Delta(R_A)$ and $\Delta(R_B)$ are less than $5\eta^{\frac{1}{2}} n$.  
	\label{s32}
	
		\item	We  require that $U^*\cap V(R_A) =\emptyset$, and 
	all the edges  of $R_B$ that are incident with a vertex  of $U^*$ either  incident with 
	one vertex of $\{v_2, v_p, v_{p+1}\}$ or have degree one in $R_B$; 
	edges of $R_B$ that are incident with a vertex of $\{v_2, v_p, v_{p+1}\}$ are 
	contained in $F_{21}\cup F_{22} \cup E_{R_B}(v_2, V(R_B)\setminus U^*)\cup E_{R_B}(v_{p+1}, V(R_B)\setminus U^*)$,
	and $e_{R_B}(v_2, V(R_B)\setminus U^*),  e_{R_B}(v_{p+1}, V(R_B)\setminus U^*) <9\eta n$. 
	If $p\ne 2$, then 
	$d_{R_B}(v_2) \le |\pbar_0(v_1)|<9 \eta n$; if $p=2$, then $d_{R_B}(v_2), d_{R_B}(v_3)\le |\pbar_0(v_1)|+6\eta n<15 \eta n$. 
	\label{s34}
	
	\item 	We also require that 
	every vertex  $v\in V(G_3)\setminus U^*$ is incident in $G_3$ with fewer than $|\pbar_0(v)|+1+5\eta^{\frac{1}{2}} n$ colored edges of $H-E(G_{A,B})$; 
	the vertex $v_1$ is incident in $G_3$ with  $|\pbar_0(v_1)|<9\eta n$ colored edges of $H-E(G_{A,B})$ and all these edges are from $E_{H}(v_1, v_2)$ if $v_1 \in V(G_0)\setminus V(G)$. Furthermore, in $G_3$,  
for  $v\in  U^*$ with $v\ne v_1$, if $p\ne 2$ or $v\ne v_3$, $v$  is incident in $G_3$ to  at most  $|\pbar_0(v)| +1<9\eta n$ colored edges of $H-E(G_{A,B})$; when $p=2$ and $v=v_3$,  the vertex $v_3$ is incident  in $G_3$ to  at most $|\pbar_0(v_1)|+|\pbar_0(v_3)|<18\eta n$ colored edges of $H-E(G_{A,B})$. 
	\label{s33}
\end{enumerate}

To ensure  Condition S3.\ref{s32} is satisfied,  we say that  an edge  $e=uv\in E(G_{A,B})$ is \emph{good} if $e\not \in E(R_A \cup R_B)$  and  the degree of $u$ and $v$ 
in both $R_A$ and $R_B$ is less than $5\eta^{\frac{1}{2}} n$ (actually, note that when  $uv\in E(G_A)$, then the degree of $u$
and $v$ is zero in $R_B$  and vice versa). 
Thus  a good edge can be added to $R_A$ or $R_B$ without violating S3.\ref{s32}. 
 
 We call a pair of distinct vertices $(a,b)$  a \emph{missing-common-color pair} or \emph{MCC-pair} in short with respect to a color $i$ if  $i$
 is missing at both $a$ and $b$ with respect to the current coloring. 
As  $|\pbar_0^{-1}(i)\cap A|-|\pbar^{-1}_0(i)\cap B|$  is even by the Parity Lemma, for each color 
 $i\in [1,k]$, we can pair up a vertex from $\pbar^{-1}_0(i)\cap A$ 
 and a vertex from $\pbar^{-1}_0(i)\cap B$,  and then pair up 
 the remaining unpaired vertices from  $\pbar^{-1}_0(i)\cap A$ or $\pbar^{-1}_0(i)\cap B$.  
 Thus we can form in total $|\pbar^{-1}_0(i)|/2$ MCC-pairs with respect to $i$.

 For every MCC-pair $(a,b)$ with respect to a  color $i\in[1,k]$, 
 we will exchange colors on an alternating path $P$ from $a$ to $b$ with at most 13
 edges,  where the path starts with an uncolored edge  of $H$ and alternating between uncolored 
 edges and good edges colored by $i$. 
 After $P$ is exchanged, $a$ and $b$ will be incident with edges  colored by $i$, and at most four good edges will be added to each of $R_A$ and $R_B$. With this information at hand, 
 before demonstrating the existence of such paths,  we  show that Conditions S3.\ref{s31}, S3.\ref{s32},  S3.\ref{s34}, and S3.\ref{s33}  can be guaranteed  at the end of Step 3.

 For~S3.\ref{s31}: We will add at most $9\eta n$ additional MCC-pairs in the initial stage of Step 3 when we deal with missing colors 
 at $v_1$ when $v_1\in V(G_0)\setminus V(G)$.  Thus 
 after the completion of Step 3, by~\eqref{eqn:total-missing-color-a-d}, the total number 
 of missing colors from vertices in $A \cup B$
is  less than   $9\eta n^2$. 
 Thus there are less than  $4.5 \eta n^2$ MCC-pairs. 
 For each MCC-pair $(a,b)$ with $a,b\in V(G_3)$, 
 at most four edges will be added to each of $R_A$ and $R_B$
 when we exchange an alternating path from $a$ to $b$.  Thus 
 there will always be less than 
 $18 \eta n^2$ newly uncolored 
 edges in each of $R_A$ and $R_B$.  
 At the completion of 
 Step 3, each of the $k$ color classes is a 1-factor of $G_3$
 so each of $G_A$ and $G_B$ have the same number of 
 colored edges.   As edges of $F_{21}\cup F_{22}$ are contained in $R_B$, by Observation~\ref{obv:G-AB}(ii), 
 we know that $e(R_A)=e(R_B)$. 
As $R_A$ is initially empty and less than $18 \eta n^2$ edges will be added to $R_A$
 at the end of Step 3, we have $e(R_A) <18 \eta n^2$ and so  $e(R_B) <18 \eta n^2$. 
  Thus  Condition S3.\ref{s31} will be satisfied at the end of Step 3. 
  
 For~S3.\ref{s32}: As we only ever add good edges to $R_A$ and $R_B$, Condition S3.\ref{s32} will  hold automatically for $R_A$
 and $R_B-\{v_p, v_{p+1}\}$.  As  there will be less than $9\eta n$ edges added to $R_B$ in the entire process of Step 3
 such that these edges 
 are incident with $v_p$
 or $v_{p+1}$,  and $|F_{21}|, |F_{22}| \le 6 \eta n$, 
 it follows 
  that the degree of 
 both $v_{p}$ and $v_{p+1}$ will be less than $15\eta n $ in $R_B$.  Thus $\Delta(R_B)<5 \eta^{\frac{1}{2}} n$ as well.

 For~S3.\ref{s34}:   Other than edges of $F_{21}\cup F_{22}$ which are contained in $R_B$ in
 the beginning of Step 3, all other edges of $G_{A,B}$ incident with vertices of $U^*$
 are added only to $R_B$ in the initial stage of Step 3 in the entire process of this step, where in the initial state,  we deal with missing colors
 at $v_1$ if $v_1\in V(G_0)\setminus V(G_1)$. Thus  $U^*\cap V(R_A) =\emptyset$.  
 When we deal
 with  missing colors at $v_1$,  we will only use edges 
 from $E_{H}(v_1, v_2)$ of $E(H)$ and add at most  $2|\pbar_0(v_1)|$  edges of $G_B$ to $R_B$ when we exchange an alternating path starting from $v_1$.  
 When we exchange 
 an alternating path starting at $v_1$ with respect to a color $i$,  we either add an edge $v_2 u $ from $G_B$ to $R_B$ if $\varphi_0(v_2u)=i$ and $u\not\in U^*$, 
 or when $\varphi_0(v_2u)=i$ and $u\in U^*$, 
 an edge, say  $wv_2$,   incident 
 with $v_2$  will be added to $R_B$, and another edge $uv$ of $G_B$,  
 will also be added to $R_B$.   We choose $w,v$ such that $w,v\in B\setminus U^*$ and $w,v$ have not yet  included in $R_B$.  Under 
 the second case that $\varphi_0(v_2u)=i$ and $u\in U^*$, if $u\ne v_{p+1}$, then by the construction of $G_3$, we know that $e_{G_3}(v_2,u)=1$. 
 As colors in $\pbar_0(v_1)$ are distinct, it follows that all the edges  of $R_B$ that are incident with a vertex  of $U^*$ either  incident with 
 one vertex of $\{v_2, v_p, v_{p+1}\}$ or have degree one in $R_B$. 
 Since there will be at most $ |\pbar_0(v_1)|$ edges  incident with $v_2$ added to $R_B$ in
 the initial stage and also in the entire Step 3,  we know that $d_{R_B}(v_2) \le |\pbar_0(v_1)|<9 \eta n$   if $p\ne 2$, and $d_{R_B}(v_2), d_{R_B}(v_3) \le |\pbar_0(v_1)|+6\eta n<15 \eta n$ 
 if $p=2$ ($v_3$ could be the vertex $u$ with $u\in U^*$ in the entire initial stage).   
All edges of $R_B$ that are incident with a vertex of $\{v_2, v_p, v_{p+1}\}$ are 
contained in $F_{21}\cup F_{22}$,  or   from $E_{R_B}(v_2, V(R_B)\setminus U^*) \cup E_{R_B}(v_{p+1}, V(R_B)\setminus U^*)$ by our choice of vertices like $w, v$ above.    
 Thus Condition S3.\ref{s34}   will be satisfied.

 For~S3.\ref{s33}:
 In the process of Step 3,   the number of newly colored edges of $H$ that are incident with a vertex   $u\in V(G_3)\setminus U^*$ will equal  the number of alternating paths containing $u$ that have been exchanged. As $u$ can have degree at most one after the initial stage of Step 3, 
 the number of such alternating paths of which $u$ is the first vertex  at the end of the initial stage of Step 3 is at most $1+|\pbar_0(u)|$.  
 The number of alternating paths in which $u$ is not the first  vertex will  equal  the degree of $u$ in $R_A\cup R_B$, and so will be less than $5\eta^{\frac{1}{2}} n$. Hence the number of  edges of $H$ that are colored in Step 3 and are incident with $u$ will be less than $|\pbar_0(u)|+1+5\eta^{\frac{1}{2}} n$.  The number of edges of $H$ that are colored in Step 3 and are incident with  $v_1$ is $|\pbar_0(v_1)|$,  as $v_1$ will only be used in alternating paths starting at $v_1$. By~S3.\ref{s34}, when $v_1\in V(G_0)\setminus V(G)$, all these newly colored 
 edges incident with $v_1$ are from $E_H(v_1, v_2)$.  
The number of newly colored edges of $H$ that are incident with a vertex $u\in U^*\setminus \{v_1\}$, if $p\ne 2$ or $u\ne v_3$,  will be at most   $|\pbar_0(u)|+1<9 \eta n$  at the end of Step 3,  as the  number of edges of $H$ that are colored in Step 3 and are incident with  $u$   will only be used in alternating paths  starting at $u$ after the initial stage of Step 3, and after the initial stage,  the total number 
of colors missing at $u$ is at most $|\pbar_0(u)|+1<9\eta n$ by S3.\ref{s34} (the set of colors missing at $v_2$ after the initial stage is the same as $\pbar_0(v_2)$).   If $p=2$ and $u=v_3$, then the number of newly colored edges of $H$ that are incident with $v_3$ will be at most $|\pbar_0(u)|+|\pbar_0(v_1)|<18\eta n$ by S3.\ref{s34}, as again, $v_3$   will only be used in alternating paths  starting at $v_3$ after the initial stage of Step 3. 
  Thus Condition S3.\ref{s33}   will be satisfied. 
  
  \medskip

\hl{Initial Stage}. \quad 
For  the vertex $v_1$, if it was added to $G$ to form $G_0$, it may have a quite small simple degree in $G_3$,  so we deal with 
colors in $\pbar_0(v_1)$ first  if $v_1\in V(G_0)\setminus V(G_1)$.  As $d_{G_2}(v_1)=\delta$, 
 we know that $e_{G_3}(v_1,v_2) = (\Delta -\delta)$ if $p\ge 4$
and $e_{G_3}(v_1,v_2) \ge (\Delta -\delta)-\frac{1}{2} \df_{G_2}(v_3) \ge \frac{1}{2}(\Delta -\delta)$ if $p=2$. 
As $|\pbar_0(v_1)| <9\eta n$ by~\eqref{s1.1},  by (iii) of the  $G_{A,B}$ formation Procedure, the number of uncolored edges  of $H$ from $E_{G_3}(v_1, v_2)$ 
is at least 
\begin{eqnarray*}
	   \frac{1}{2} (\Delta-\delta ) -|\pbar_0(v_1)| &\ge& \frac{1}{2} (\Delta-\delta) -9\eta n  \ge \frac{1}{4}(1-3\ve )n-9\eta n>0 \quad \text{if $p\ge 4$},
\end{eqnarray*}
 and at least 
\begin{eqnarray*}
	\frac{1}{4} (\Delta-\delta ) -|\pbar(v_1)| &\ge& \frac{1}{4} (\Delta-\delta) -9\eta n   \ge \frac{1}{8}(1-3\ve )n-9\eta n>0 \quad \text{if $p=2$},
\end{eqnarray*}
where we used $\delta <\frac{1}{2}(1+\ve) n$ from~\eqref{eqn:upper-bound-delta} in the calculations above.

 For each $i\in \pbar_0(v_1)$, if  $i \in \pbar_0(v_2)$, we take an uncolored edge $e$ from $E_{G_3}(v_1, v_2)$  
 and color $e$ by $i$. Thus we assume that $i \notin \pbar_0(v_2)$.  Since only edges incident with $v_1$  in $H-E(G_{A,B})$ have been colored in Step 3 so far and $(G_{A,B}\cap H)-v_1$ contains no edge incident with $v_2$ by the  $G_{A,B}$ Formation Procedure ($v_2$ is only adjacent to $v_1$ in $L^*\cap H$), 
 it follows that there exists $w\in N_{G_{B}}(v_2)$ such that $v_2 w$ is colored by $i$. 
 If $w\in B\setminus U^*$, 
 we take an uncolored edge $e$ from $E_{G_3}(v_1, v_2)$, 
 color $e$ by $i$, and uncolor $v_2w$.  Now the color $i$ presents at the vertex $v_1$ 
 but is missing at $w$. 
 (We say that $w$ is \emph{originated} from $v_1$ with respect to the color $i$). 
 If $w\in U^*$, then  as 
 $d_G(w) \ge d_G(v_2) \ge \delta$, 
 $|U| <2\eta n$, $|\pbar_0(v_2)|, |\pbar_0(w)|<9\eta n$ by~\eqref{s1.1}, there 
 exist distinct  $x\in N_{G_B}(v_2)\setminus U$ and $y\in N_{G_B}(w)\setminus U$ such that 
 $\varphi_0(v_2x)=\varphi_0(wy)$ and that $x$ and $y$ are incident in $G_B$ with no uncolored edge so far.  
 Now we take an uncolored edge $e$ from $E_{G_3}(v_1, v_2)$, 
 color $e$ by $i$,  recolor $v_2w$ by the color used on $v_2x$ and $wy$ under $\varphi_0$, and 
 uncolor both $v_2x$ and $wy$.  Now the color $i$ presents at $v_1$ and is missing at both $x$ and $y$,
 and we also say that $x$ and $y$  are \emph{originated} from $v_1$ with respect to the color $i$.

  The  procedure above 
 guarantees  that all edges of $R_B$ to be good, as  $|\pbar_0(v_1)| < 9\eta n$ and the maximum degree of $R_B$ is at most $6\eta n$ before we start the procedure. 
 After this procedure, all the $k$ colors present at $v_1$. 
 Thus $v_1$ is not contained in any MCC-pairs with respect to the current edge coloring if  $v_1\in V(G_0)\setminus V(G)$. 
 We denote by $\varphi_1$ the current edge coloring.  
 As every time when we uncolor an edge at $v_2$ in the initial stage, we also color 
 an edge from $E_H(v_1, v_2)$ using the color we just took off, thus we have $\pbar_1(v_2)=\pbar_0(v_2)$.  For all the uncolored edges like $v_2w$, $v_2x$ or $wy$ in the stage,
 the process guarantees that $w, x,y\in B\setminus U^*$ and that  $w, x,y$ have  degree one in $R_B$.
 Furthermore,  for any  $v\in B$ with $v\ne v_2$,  if $p\ne 2$ or $v\ne v_3$, we have $|\pbar_1(v)| \le |\pbar_0(v)| +1$;
 and if $p=2$ and $v=v_3$, we have $|\pbar_1(v_3)| \le |\pbar_0(v_3)| +|\pbar_0(v_1)| <18\eta n$.

 Let $U_A=(U^*\cap A) \cup \{u_{p+1}\}$ and  $U_B=U^*\cap B$. 
 We now show below the existence of alternating paths for the current MCC-pairs.
    For a given color $i\in [1,k]$,  and vertices $a\in A$ and $b\in B$, let $N_B(a)$ be the set
 of vertices in $B$ that are joined with $a$ by an uncolored edge and are incident with a
 good edge colored $i$ such that the good edge is not incident with any vertex of $U_B$, and let $N_A(b)$ be the set of vertices in $A$ that are joined with
 $b$ by an uncolored edge and are incident with a good edge colored $i$ such that the good edge is not incident with any vertex of $U_A$. 
 In order to estimate the sizes of $N_A(b)$ and $N_B(b)$, we 
 show that $A$ and $B$ contain only a few vertices 
 that either miss the color $i$, 
 or are incident  with a non-good edge colored $i$, or are incident with a good edge colored $i$
 but the good edge is incident with a vertex from $U^*$.  
 By S3.\ref{s31}, there are less than  $18 \eta n^2$ edges in $R_B$, so there are fewer than $\frac{36 \eta n^2}{5 \eta^{\frac{1}{2}}n}  \le 7.2 \eta^{\frac{1}{2}} n$ vertices of
 degree at least $5 \eta^{\frac{1}{2}}n$ in $R_B$. 
 Each non-good edge is incident with one or two 
 vertices of $R_B$ through the color $i$, 
 so there are fewer than
 $
 2\times 7.2\eta^{\frac{1}{2}} n =14.4 \eta^{\frac{1}{2}}n
 $ 
 vertices in $B$ that are
 incident with a non-good edge colored $i$. 
 Furthermore, there are at most $2|U_B| <2|U| \le 4\eta n$ vertices in $B$
that are  either contained in $U_B$ 
 or adjacent to a vertex from $U_B$ through an edge with color $i$. 
Finally, there are fewer than $9 \eta n$
 vertices in $B$ that are missing the color $i$ by~\eqref{s1.1} and the procedure in the Initial Stage. So the number of vertices in $B$ that are
 not incident with a good edge colored $i$  such that the good edge is not 
 incident with any vertex from $U_B$ is less than
 $$
 14.4\eta^{\frac{1}{2}}n+ 4\eta n+9\eta n<15 \eta^{\frac{1}{2}}n. 
 $$
 Similarly, 
 the number of vertices  in $A$ that are
 not incident with a good edge colored $i$ such that the good edge is not 
 incident with any vertex from $U_A$ is less than
 $15\eta^{\frac{1}{2}}n$.
 
We have  $|\pbar_0(v)| <9 \eta n-2$ for any $v\in V(G_3)$ by~\eqref{s1.1}. 
If $v_1\in V(G_0)\setminus V(G_1)$, 
after the initial operation on the vertex $v_1$, we have $|\pbar_1(v)| \le |\pbar_0(v)|+1<9 \eta n$ if $p\ne 2$ or $v\ne v_3$,
and $|\pbar_1(v_3)| \le |\pbar_0(v_3)|+|\pbar_0(v_1)|<18\eta n$ if $p=2$ and $v=v_2$ by S3.\ref{s33}. 
Note that   $\pbar_1(v_2) =\pbar_0(v_2)$ by the  initial operation. 

Using  the inequality in~\eqref{eqn:degree-diff}, 
we have the following lower bounds on $|N_B(a)|$ and $|N_A(b)|$. Let $c\in \{a,b\}$. 
\begin{enumerate}[]
	\item When $ c\not\in  U^* $ and $c\ne u_{p+1}$, as $\ve \le \frac{1}{14}$, 
	\begin{eqnarray}\label{eqn: not-in-S-AB}
		|N_A(c)|, |N_B(c)| &\ge& 
		\frac{1}{2}\left((\Delta -3\eta)n-(n^{\frac{2}{3}}+4\eta n+2)\right)-(9\eta n+5\eta^{\frac{1}{2}}n)-15\eta^{\frac{1}{2}}n  
		>\frac{3}{7} n. 
	\end{eqnarray} 
where $9\eta n+5\eta^{\frac{1}{2}}n$ is the  upper bound of the number of edges of $H-E(G_{A,B})$  incident to $c$ that were colored  in Step 3. 
\item When $c= u_{p+1}$ (recall that $u_{p+1}$ is the partner of $v_{p+1}$), 
\begin{eqnarray}\label{eqn: not-in-S-AB3}
	|N_A(c)|, |N_B(c)| &\ge& 
	\frac{1}{2}\left((\Delta -3\eta)n-(n^{\frac{2}{3}}+4\eta n+2)\right)-|E_2|-(9\eta n+5\eta^{\frac{1}{2}}n)-15\eta^{\frac{1}{2}}n  \nonumber \\
	&> &\frac{\Delta}{3}-21\eta^{\frac{1}{2}}n>\frac{2}{7} n, 
\end{eqnarray} 
where we used $|E_2|=|E_1| <\frac{\Delta}{6}$ by~\eqref{eqn:E2-size}. 

\item When $c\in  U$, 
 we have 
	\begin{eqnarray}\label{eqn: not-in-S-AB2}
	|N_A(c)|, |N_B(c)| &\ge& 
	\frac{1}{2}\left(\delta-(n^{\frac{2}{3}}+4\eta n+2)\right)-18\eta n-15\eta^{\frac{1}{2}}n  
	>\frac{2}{49} n,
\end{eqnarray} 
where $\delta>\frac{1}{7}(1-4\ve)n$ by~\eqref{eqn:upper-bound-delta}. 
\end{enumerate}

 Let $M_B(a)$ be the set of vertices in $B$ that are joined with a vertex in $N_B(a)$ by an edge
 of color $i$, and let $M_A(b)$ be the set of vertices in $A$ that are joined with a vertex in $N_A(b)$
 by an edge of color $i$.  Note that $(U_A \cup U_B)\cap ( M_A(b) \cup M_B(a))=\emptyset$ by the choice of $N_A(b)$ and $N_B(a)$. 
 Note also that $|M_B(a)|=|N_B(a)|$     
 but some vertices
 may be in both. Similarly
 $|M_A(b)|=|N_A(b)|$. 
 
 For an MCC-pair $(a,b)$, in order to have a unified discussion as in the case that $\{a,b\}\cap (U \cup \{u_{p+1}\})= \emptyset$,  if necessary, by exchanging an alternating path of length 2 from $a$ to another vertex $a^*$,
 and exchanging an alternating path from $b$ to another vertex $b^*$, 
 we will replace the pair $(a,b)$ by $(a^*,b^*)$ such that $\{a^*,b^*\}\cap (U \cup \{u_{p+1}\})=\emptyset$.  
 Precisely, we will implement the following 
 operations to vertices in $U \cup \{u_{p+1}\}$. 
 For any vertex $a\in U_A$, and for each color $i\in \pbar_1(a)$, we take an edge $b_1b_2$ with $b_1\in N_B(a)$ and $b_2\in M_B(a)$ such that $b_1b_2$ is colored by $i$, where the edge $b_1b_2$ exists by~\eqref{eqn: not-in-S-AB3} and~\eqref{eqn: not-in-S-AB2} and the fact that $|M_B(a)|=|N_B(a)|$.   Then we exchange the path $ab_1b_2$ by coloring $ab_1$
 with $i$ and uncoloring the edge $b_1b_2$ (See Figure~\ref{f1}(a)).
  After this, the edge $ab_1$  of $H$ is now colored by $i$, and the uncolored edge $b_1b_2$ is added to $R_B$.
 %We then update the original MCC-pair that contains $a$ with respect to the color $i$ by replacing the vertex $a$ with $b_2$. 
 We do this at the vertex $a$ for every color $i\in \pbar_1(a)$ and then repeat the same process for every vertex in $U_A$. 
 Similarly, 
 for any vertex $b\in U_B$, and for each color $i\in \pbar_1(b)$, we take an edge $a_1a_2$ with $a_1\in N_A(b)$ and $a_2\in M_A(b)$ such that $a_1a_2$ is colored by $i$, where the edge $a_1a_2$ exists by~\eqref{eqn: not-in-S-AB3} and~\eqref{eqn: not-in-S-AB2} and the fact that $|M_B(b)|=|N_B(b)|$.   Then we exchange the path $ba_1a_2$ by coloring $ba_1$
 with $i$ and uncoloring the edge $a_1a_2$. 
 %The same, we update the original MCC-pair that contains $b$ with respect to the color $i$ by replacing the vertex $b$ with $a_2$.  
 (Although in this process, we increased the number of missing colors at vertices like $a_2$ and $b_2$, but the increase is calculated as the degrees of $a_2$ and $b_2$ in $R_A\cup R_B$. Thus 
 all the calculations in~\eqref{eqn: not-in-S-AB} to \eqref{eqn: not-in-S-AB2} are still valid after the implementations above.)

After the procedure above,    we have now three types MCC-pair $(u,v)$ with respect to a given color $i$: $u,v\in A$, $u,v\in B$, and $A$
contains exactly one of $u$ and $v$ and $B$ contains the other. However, in either case, $\{u,v\}\cap (U \cup \{u_{p+1}\})=\emptyset$. We will exchange alternating path for each of such pairs. 

We deal with each of the colors from $[1,k]$ in turn. 
Let $i\in [1,k]$ be a color.  We consider first an MCC-pair $(a,b)$
with respect to  $i$ such that $a\in A$ and $b\in B$. 
By~\eqref{eqn: not-in-S-AB}, we have 
$
|M_B(a)|,  |M_A(b)|>\frac{3}{7} n.  
$
We choose $a_1a_2$ with color $i$ such that $a_1\in N_A(b)$
and $a_2\in M_A(b)$. Now as $|M_B(a)|, |N_B(a_2)| >\frac{3}{7} n$ by~\eqref{eqn: not-in-S-AB}, we know that $N_B(a_2)\cap M_B(a) \ne \emptyset$. We choose $b_2\in N_B(a_2)\cap M_B(a)$
and let $b_1\in N_B(a)$ such that $b_1b_2$ is colored by $i$. 
Then $P=ab_1b_2a_2a_1 b$ is an alternating path from $a$ to $b$ (See Figure~\ref{f1}(b)). We exchange $P$ by coloring $ab_1, b_2a_2$ and $a_1b$
with color $i$ and uncoloring the edges $a_1a_2$ and $b_1b_2$. 
After the exchange,  the color $i$ appears on edges incident with $a$ and $b$,
the edge $a_1a_2$ is added to $R_A$
and the edge $b_1b_2$ is added to $R_B$. 
As $a$ and $b$ could be originated from some vertices $b'\in U_B$ and $a'\in U_A$, respectively, 
and  $b'$  could be originated from  $v_1$
with respect to the color $i$, 
we added at most 
two edges to each of $R_A$ and $R_B$ prior to having $(a,b)$ as an MCC-pair with respect to $i$ such that 
$\{a,b\}\cap (U^*\cup \{u_{p+1}\}) =\emptyset$. 
  Thus  we added at most three edges to each of $R_A$ and   $R_B$ for each MCC-pair with respect to $i$ under $\varphi_0$.

We consider then  an MCC-pair $(a,a^*)$
with respect to  $i$ such that $a,a^* \in A$. 
By~\eqref{eqn: not-in-S-AB}, we have  $|M_B(a^*)|>\frac{3}{7} n$.
We take an edge $b_1^*b_2^*$ colored by $i$ with $b_1^* \in N_B(a^*)$ and  $b_2^*\in M_B(a^*)$.  
Then again, by~\eqref{eqn: not-in-S-AB}, we have 
$
 |M_B(a)|,  |M_A(b_2^*)|>\frac{3}{7} n.  
$
Therefore, as each vertex  $c\in M_A(b_2^*)$ 
satisfies $|N_B(c)|>\frac{3}{7} n$,  we have $|N_B(c)\cap M_B(a)| >\frac{2}{7} n$. We take $a_2a_2^*$ colored by $i$ with $ a_2^*\in N_A(b_2^*)$ and  $a_2\in   M_A(b_2^*)$. 
Then we let $b_2\in N_B(a_2)\cap M_B(a)$, and let $b_1$
be the vertex in $N_B(a)$ such that $b_1b_2$ is colored by $i$. 
Now we get the alternating path $P=ab_1b_2 a_2 a_2^* b_2^* b_1^* a^*$ (See Figure~\ref{f1}(c)).  
We exchange $P$ by coloring $ab_1, b_2a_2, a_2^*b_2^*$ and $b_1^*a^*$
with color $i$ and uncoloring the edges $b_1b_2, b_1^*b_2^*$ and $a_2a_2^*$. 
After the exchange,  the color $i$ appears on edges incident with $a$ and $a^*$,
the edges $b_1b_2$ and $b_1^*b_2^*$ are added to $R_B$
and the edge $a_2a_2^*$ is added to $R_A$. 
As $a$ and $a^*$ could be originated from some vertices $b_1', b_2'\in U_B$, and  
$b_1'$ and $b_2'$ could be originated from  $v_1$
with respect to the color $i$, 
we added at most 
two edges to each of $R_A$ and $R_B$ prior to having $(a,a^*)$ as an MCC-pair with respect to $i$ such that 
$\{a,a^*\}\cap (U^*\cup \{u_{p+1}\}) =\emptyset$. 
Thus  we added at most four  edges to $R_A$
and at most four edges to $R_B$ for each MCC-pair with respect to $i$ under $\varphi_0$.  The maximum length of an alternating path combined from the three procedures (dealing missing colors at $v_1$, dealing missing colors at vertices from $U\cup \{u_{p+1}\}$, and dealing missing colors at vertices from $V(G_3)\setminus (U^*\cup \{u_{p+1}\})$) together is at most $2+2+2+7=13$. 

By symmetry, we can  deal with an MCC-pair $(b,b^*)$
with respect to  $i$ such that $b,b^* \in B$ similarly as above.  
 By finding
 such paths  for all MCC-pairs with respect to the color $i$,  we can increase the number of
 edges colored  by $i$ until the color class is a 1-factor of $G_3$. By doing this for all colors,
 we can make each of the $k$ color classes  into a 1-factor of $G_3$.
We denote by $\varphi_2$ the resulting coloring at the end of Step 3. 

\begin{figure}[!htb]
	\begin{center}
		
		\begin{tikzpicture}[scale=0.8]
		
			\begin{scope}[shift={(-10,0)}]
		\draw[rounded corners, fill=white!90!gray] (8, 0) rectangle (12, 2) {};
		
		\draw[rounded corners, fill=white!90!gray] (8, -4) rectangle (12, -2) {};
		
		{\tikzstyle{every node}=[draw ,circle,fill=white, minimum size=0.5cm,
			inner sep=0pt]
			\draw[black,thick](9,1) node (c)  {$a$};
			\draw[black,thick](9,-3) node (d)  {$b_{1}$};
			\draw[black,thick](11,-3) node (d1)  {$b_{2}$};
		}
		\path[draw,thick,black!60!white,dashed]
		(c) edge node[name=la,pos=0.7, above] {\color{blue} } (d)
		;
		
		\path[draw,thick,black]
		(d) edge node[name=la,pos=0.7, above] {\color{blue} } (d1)
		;
		\node at (7.2,1) {$A$};
		\node at (7.2,-3) {$B$};
		\node at (10,-4.8) {$(a)$};	
		\end{scope}	
		
		\begin{scope}[shift={(5,0)}]
		\draw[rounded corners, fill=white!90!gray] (0, 0) rectangle (6, 2) {};
		
		\draw[rounded corners, fill=white!90!gray] (0, -4) rectangle (6, -2) {};
		
		{\tikzstyle{every node}=[draw ,circle,fill=white, minimum size=0.5cm,
			inner sep=0pt]
			\draw[black,thick](1,1) node (a)  {$a$};
			\draw[black,thick](3,1) node (a1)  {$a_1$};
			\draw[black,thick](5,1) node (a2)  {$a_2$};
			\draw[black,thick](1,-3) node (b)  {$b$};
			\draw[black,thick](3,-3) node (b1)  {$b_1$};
			\draw[black,thick](5,-3) node (b2)  {$b_2$};
		}

		\path[draw,thick,black!60!white,dashed]
		(a) edge node[name=la,pos=0.7, above] {\color{blue} } (b1)
		(a2) edge node[name=la,pos=0.7, above] {\color{blue} } (b2)
		(b) edge node[name=la,pos=0.6,above] {\color{blue}  } (a1)
		;
		
		\path[draw,thick,black]
		(a1) edge node[name=la,pos=0.7, above] {\color{blue} } (a2)
		(b1) edge node[name=la,pos=0.7, above] {\color{blue} } (b2)
		;
		
		\node at (-0.5,1) {$A$};
		\node at (-0.5,-3) {$B$};
		\node at (3,-4.8) {$(b)$};

		\end{scope}

		\begin{scope}[shift={(-8,-8)}]
		\draw[rounded corners, fill=white!90!gray] (8, 0) rectangle (16, 2) {};
		
		\draw[rounded corners, fill=white!90!gray] (8, -4) rectangle (16, -2) {};
		
		{\tikzstyle{every node}=[draw ,circle,fill=white, minimum size=0.5cm,
			inner sep=0pt]
			\draw[black,thick](9,1) node (c)  {$a$};
			\draw[black,thick](11,1) node (c1)  {$a_{2}$};
			\draw[black,thick](13,1) node (c2)  {$a^*_{2}$};
			\draw[black,thick](15,1) node (c3)  {$a^*$};
			\draw[black,thick](9,-3) node (d)  {$b_{1}$};
			\draw[black,thick](11,-3) node (d1)  {$b_{2}$};
			\draw[black,thick](13,-3) node (d2)  {$b^*_{2}$};
			\draw[black,thick](15,-3) node (d3)  {$b^*_{1}$};
		}
		\path[draw,thick,black!60!white,dashed]
		(c) edge node[name=la,pos=0.7, above] {\color{blue} } (d)
		(c1) edge node[name=la,pos=0.7, above] {\color{blue} } (d1)	
		(c3) edge node[name=la,pos=0.7, above] {\color{blue} } (d3)	
		(c2) edge node[name=la,pos=0.7, above] {\color{blue} } (d2)	
		;
		
		\path[draw,thick,black]
		(d) edge node[name=la,pos=0.7, above] {\color{blue} } (d1)
		(c2) edge node[name=la,pos=0.7, above] {\color{blue} } (c1)
		(d3) edge node[name=la,pos=0.7, above] {\color{blue} } (d2)
		;
		\node at (7.2,1) {$A$};
		\node at (7.2,-3) {$B$};
		\node at (12.5,-4.8) {$(c)$};	
		\end{scope}	
		
		\end{tikzpicture}
			\end{center}
	\caption{The alternating path $P$. Dashed lines indicate uncoloured edges, and solid
		lines indicate edges with color $i$.}
	\label{f1}
\end{figure}

\begin{center}
{\bf 	Step 4: Coloring $R_A$ and $R_B$ and extending the new color classes in 1-factors }
\end{center}

Each of the color classes for the colors from $[1,k]$ is now a 1-factor of $G_3$. We 
now consider the multigraphs $R_A$ and $R_B$ that consist of the uncolored edges of $G_A-E_1$ and $G_B-F_1$. 
By Conditions S3.\ref{s31} and~S3.\ref{s32}, $R_A$ and $R_B$ each has fewer than $18\eta n^2$ edges, and $\Delta(R_A ), \Delta(R_B)< 5\eta^{\frac{1}{2}} n$.  
Let $R$ be the subgraph of $G_3$ consisting of the remaining uncolored edges at the completion of Step 4. 
We need to make sure that 
\begin{equation}\label{eqn:v_p,p-1,p=1-edge-in-R}
	e(R[\{v_{p-1}, v_p, v_{p+1}\}]) \le \Delta(R). 
\end{equation}
For this purpose, we color edges in this step in two stages.

Let $F=E_{R_B}(v_2, V(R_B)\setminus \{v_2\}) \cup E_{R_B}(v_p, V(R_B)\setminus \{v_p\})$ be the set of edges in $R_B$ incident with $v_2$ or $v_p$,  $\ell_1=|F|$, and $F^*$
be an arbitrary set of $\ell_1$ edges of $R_A$.   Note that by Condition S3.\ref{s34}, $\ell_1<9 \eta n$ if $p\ne 2$, and $\ell_1 < 15\eta n$ if $p=2$. 
We will edge color each of $R_A$
and $R_B$ and color a few uncolored edges of $H$ using another 
$\ell$ colors, where $$\ell =\ell_1+\ell_2 \quad \text{and} \quad \ell_2:=\Delta(R_A\cup R_B)+\mu(R_A\cup R_B) \le \lceil  5\eta^{\frac{1}{2}} n+3\eta n\rceil.$$  

\medskip 
\hl{Stage 1:  Color edges in $F\cup F^*$ using $\ell_1$ colors.}

\medskip

We color each edge in $F$ using a distinct color from $[k+1, k+\ell_1]$, 
and also color each edge in $F^*$ using a distinct color from $[k+1, k+\ell_1]$. 
Given a
color $i$ with $i\in [k+1, k+\ell_1]$, we let $A_i$ and $B_i$ be the sets of vertices in $A$ and
$B$ respectively that are incident with edges colored by $i$. 
Note that $|A_i| =  |B_i|   =2$.
Let $H_i$ be the subgraph of $H$ obtained by
deleting the vertex sets $A_i $ and $B_i$ and removing all colored edges. We will show next that $H_i$ has a perfect matching and we will color 
the edges in the matching by the color $i$ to extend the color class $i$ into a 1-factor of $G_3$. 

As $d^s_{H_i}(v_1)$ can be very small, and we 
want to use edges of $E_{H_i}(v_{p-1}, v_p) \cup E_{H_i}(v_{p-1}, v_{p+1})$ as many as possible 
so that~\eqref{eqn:v_p,p-1,p=1-edge-in-R} holds at the end of Step 4, we match $v_1$
and $v_{p-1}$ to other vertices first. We consider the following cases.

{\bf Case 1: $p=2$ and $e_p=0$. } In this case, $v_1=v_{p-1}$. 
If $v_{1}\in V(G)$, then 
$d_{G_3}^s(v_{1})\ge \delta$.  
As $v_1$ is incident with at most $9\eta n$ colored edges of $H-G_{A,B}$
by Condition~S3.\ref{s34}, and at most $\ell_1$ edges of $H$ incident with $v_1$ have been colored in 
Step 4, it follows that 
\begin{eqnarray*}
	d_{H_i}^s(v_{1})&> & \frac{1}{2}\left(\delta-(n^{\frac{2}{3}}+4\eta n +2)\right)-9\eta n -\ell_1-2>0. 
\end{eqnarray*}
Thus we can match $v_{1}$ to one of its neighbors, say  $w$,  in $H_i$. 
If  $v_1\not\in V(G)$,  
assume that $v_2u$ is the edge with color $i$.  
By the formation of $R_B$, we know that $u\ne v_3$. 
Then 
by Condition S3.\ref{s34}, we know that $u\not\in U^*$. 
Thus $e_{G_3}(v_1,u)<3\eta n$ by~\eqref{eqn: vertex-multiplicity-not-from-U*}. 
As $d_{G_3-v_2}(v_1)=\delta$ and so $d_{G_{A,B}-v_2}(v_1) \ge \frac{1}{2}\delta $ by 
(iii) of $G_{A,B}$ Formation Procedure, we have 
\begin{eqnarray*}
	d_{H_i}(v_{1}) &\ge&  \frac{1}{2}\delta -9\eta n -\ell_1-e_{L}(v_1, B_i)
	>\frac{1}{2}\delta-(9+15+3)\eta n>0. 
\end{eqnarray*}
Again, we can match $v_{1}$ to one of its neighbors, say  $w$, in $H_i$. 

{\bf Case 2: $p=2$ and $e_p>0$. } 
This implies that $\lceil \frac{1}{2}e_{G_3}(v_p, v_{p+1}) \rceil \ge 6\eta^{1/2} n$. 
Since $e_{G_3}(v_{p-1}, v_{p+1}) \ge e_{G_3}(v_p, v_{p+1}) -1$ by the construction of $G_3$, 
it follows that $\lfloor  \frac{1}{2}e_{G_3}(v_{p-1}, v_{p+1}) \rfloor\ge 6\eta^{1/2} n-2$. 
Then we get 
\begin{eqnarray*}
	e_{H_i}(v_{1},v_{3})&> &6\eta^{1/2} n-2-9\eta n -\ell_1>0,  
\end{eqnarray*}
and so we can match $v_{1}$ to $w:=v_{3}$  in $H_i$.

{\bf Case 3: $p\ge 4$ and $e_p=0$. }
If $v_{1}\in V(G)$, then 
$d_{G_3}^s(v_{1})\ge \delta$.  Then the same as in Case 1, 
we can again match $v_{1}$ to one of its neighbors, say  $w$,  in $H_i$. 
If $v_1\not\in V(G)$, 
then we have 
$e_{H_i}(v_1,v_2) \ge \frac{1}{2}(\Delta-\delta)-9\eta n -\ell_1.$
Thus we can match $v_1$ to $w:=v_2$ in $H_i$. 
Since $d_{H_i}^s(v_{p-1})> \frac{1}{2}\left(\delta-(n^{\frac{2}{3}}+4\eta n +2)\right)-9\eta n -\ell_1-2>0$, 
we can match $v_{p-1}$ to one of its neighbors, say  $w_1$,  in $H_i-\{v_1,w\}$.

{\bf Case 4: $p\ge 4$ and $e_p>0$. }
We match $v_1$ to a vertex $w$ of $H_i$ the same way as in Case 3. 
Since $e_{G_3}(v_{p-1}, v_{p+1}) \ge e_{G_3}(v_p, v_{p+1}) -1$ by the construction of $G_3$, 
it follows that $\lfloor  \frac{1}{2}e_{G_3}(v_{p-1}, v_{p+1}) \rfloor\ge 6\eta^{1/2} n-2$. 
Then we get 
$e_{H_i}(v_{p-1},v_{p+1})> 6\eta^{1/2} n-2-9\eta n -\ell_1>0,$
and so we can match $v_{p-1}$ to $w_1:=v_{p+1}$  in $H_i$.

Let $H_i^*=H_i-\{v_{1}, w\}$ if $p=2$ and $H_i^*=H_i-\{v_{1}, w, v_{p-1}, w_1\}$ if $p\ge 4$.   
If $u\in  U^*\cap V(H_i^*)$,  then we have 
\begin{eqnarray*}
	d_{H^*_i}^s(u)&> & \frac{1}{2}\left(\delta-(n^{\frac{2}{3}}+4\eta n +2)\right)-9\eta n -\ell_1-2>2\eta n>|V(H^*_i)\cap U^*|. 
\end{eqnarray*}
For each $u\in V(H^*_i)\setminus U^*$,   we have 
\begin{eqnarray*}
	d_{H_i^*}^s(u)&> & \frac{1}{2}\left((\Delta-3\eta n)-(n^{\frac{2}{3}}+4\eta n +2)\right)-(9\eta n+5\eta^{\frac{1}{2}}n)   -\ell_1-2 \\
	&>&  \frac{1}{2}\Delta-6 \eta^{\frac{1}{2}} n>\frac{n}{4}.  
\end{eqnarray*}
Now applying Lemma~\ref{lem:matching-in-bipartite}, $H^*_i$ has a perfect matching $M$. 
Let $M^*=M\cup \{v_1w\}$ if $p=2$, and $M^*=M\cup \{v_1w, v_{p-1}w_1\}$ if $p\ge 4$. 
Then $M^*$ is a perfect matching of $H_i$. 
We color all edges of $M^*$ by the color $i$. 
This extends the color class $i$ into a 1-factor of $G_3$. 
We repeat this procedure for each of the colors from $[k+1,k+\ell_1]$.  After this has been
done, each of these $\ell_1$ colors  in Stage 1 of Step 4 have  presented at all vertices of $G_3$.

\medskip 
\hl{Stage 2:  Color the rest edges of $R_A\cup R_B$ using $\ell_2$ colors.}

\medskip 

Since all edges of $L^*$ are incident with a  vertex of  $\{v_{2i+1}: i\in [0, \frac{p}{2}] \}$ and  $\{v_{2i-1}: i\in [1, \frac{p}{2}] \}\cap (V(R_A\cup R_B))=\emptyset$, 
we know that the maximum number of edges joining two vertices in $R_A \cup R_B$ have one of their endvertices from $V(G_3)\setminus U^*$. 
Therefore, by~\eqref{eqn: vertex-multiplicity-not-from-U*} and~S3.\ref{s34}, we get 
$\mu(R_A), \mu(R_B) <  3\eta n$.  
Since $\Delta(R_A ), \Delta(R_B)< 5\eta^{\frac{1}{2}} n$, and $\mu(R_A), \mu(R_B) <  3\eta n$, 
by Theorem~\ref{thm:chromatic-index}  and Lemma~\ref{lem:equa-edge-coloring}, 
there is an equalized edge $\ell_2$-coloring   of the rest uncolored edges of  $R_A\cup R_B$.  
Since $e(R_A)=e(R_B)$,  by renaming some color classes of $R_A$ if necessary,  we can assume that in the edge colorings of 
the  uncolored edges of  $R_A$ and $R_B$ after Stage 1 of Step 4, each color
appears on the same number of edges in $R_A$ as it does in $R_B$. 
Each color will appear on at most $(18 \eta n^2-\ell_1)/ 5\eta^{\frac{1}{2}} n <4 \eta^{\frac{1}{2}} n$ edges. 
Given a
color $i$ with $i\in [k+\ell_1+1, k+\ell]$, we let $A_i$ and $B_i$ be the sets of vertices in $A$ and
$B$ respectively that are incident with edges colored $i$. 
Note that $|A_i| =  |B_i|   \le 8 \eta^{\frac{1}{2}} n$.
Let $H_i$ be the subgraph of $H$ obtained by
deleting the vertex sets $A_i $ and $B_i$ and removing all colored edges. We will show next that $H_i$ has a perfect matching and we will color 
the edges in the matching by the color $i$.  By Condition~S3.\ref{s33}, we have $U^*\cap A \subseteq V(H_i)$;
and by Stage 1 of Step 4, we have $v_2, v_{p}\in V(H_i)$.  We consider two cases for matching $v_1$ and $v_{p-1}$ to other vertices first.

{\bf Case 1:  $e_p=0$. }
If $v_{1}\in V(G)$, then 
$d_{G_3}^s(v_{1})\ge \delta$.  Then  
we can again match $v_{1}$ to one of its neighbors, say  $w$,  in $H_i$. 
If $v_1\not\in V(G)$, 
then we have 
$e_{H_i}(v_1,v_2) \ge \frac{1}{2}(\Delta-\delta)-9\eta n -\ell>0.$
Thus we can match $v_1$ to $w:=v_2$ in $H_i$. If $p\ge 4$, 
since $d_{H_i}^s(v_{p-1})> \frac{1}{2}\left(\delta-(n^{\frac{2}{3}}+4\eta n +2)\right)-9\eta n -\ell-8 \eta^{\frac{1}{2}} n>0$, 
we can match $v_{p-1}$ to one of its neighbors, say  $w_1$,  in $H_i-\{v_1,w\}$.

{\bf Case 2:  $e_p>0$. } Since  
 $e_{G_3}(v_{1}, v_{2}) \ge e_{G_3}(v_p, v_{p+1}) -1$ by the construction of $G_3$, 
it follows that $\lfloor \frac{1}{2}e_{G_3}(v_{1}, v_{2}) \rfloor, \lfloor \frac{1}{2}e_{G_3}(v_{p-1}, v_{p}) \rfloor\ge 6\eta^{1/2} n-2$. 
Thus 
\begin{eqnarray*}
e_{H_i}(v_{1},v_{2}), 	e_{H_i}(v_{p-1},v_{p})&> &6\eta^{1/2} n-2-9\eta n -\ell>0,  
\end{eqnarray*}
and so we can match  $v_1$ to $w:=v_2$ and $v_{p-1}$ to $w_1:= v_{p}$  in $H_i$.

Let $H_i^*=H_i-\{v_{1}, w\}$ if $p=2$ and $H_i^*=H_i-\{v_{1}, w, v_{p-1}, w_1\}$ if $p\ge 4$.   
If $u\in  U^*\cap V(H_i^*)$,  then we have 
\begin{eqnarray*}
	d_{H^*_i}^s(u)&> & \frac{1}{2}\left(\delta-(n^{\frac{2}{3}}+4\eta n +2)\right)-9\eta n -\ell-8 \eta^{\frac{1}{2}} n-4>2\eta n>|V(H^2_i)\cap U^*|. 
\end{eqnarray*}
For each $u\in V(H^*_i)\setminus U^*$,   we have 
\begin{eqnarray*}
	d_{H_i^*}^s(u)&> & \frac{1}{2}\left((\Delta-3\eta n)-(n^{\frac{2}{3}}+4\eta n +2)\right)-(9\eta n+5\eta^{\frac{1}{2}}n)   -\ell-8 \eta^{\frac{1}{2}} n-4 \\
	&>&  \frac{1}{2}\Delta-14 \eta^{\frac{1}{2}} n>\frac{n}{4}.  
\end{eqnarray*}
Now applying Lemma~\ref{lem:matching-in-bipartite}, $H^*_i$ has a perfect matching $M$. 
Let $M^*=M\cup \{v_1w\}$ if $p=2$, and $M^*=M\cup \{v_1w, v_{p-1}w_1\}$ if $p\ge 4$. 
Then $M^*$ is a perfect matching of $H_i$. 
We color all edges of $M^*$ by the color $i$. 
This extends the color class $i$ into a 1-factor of $G_3$. 
We repeat this procedure for each of the colors from $[k+\ell_1+1,k+\ell]$.  After this has been
done, each of these $\ell_2$ colors  in Stage 2 of Step 4 have  presented at all vertices of $G_3$.

Let $R$ be the subgraph of $G_3$ consisting of the remaining uncolored edges. 
By the coloring process in Step 4, we know that when $e_p>0$, 
only edges from $E_H(v_{p-1}, v_p) \cup E_H(v_{p-1}, v_{p+1})$ were used 
to extend each of the $\ell$ color classes at the vertex $v_{p-1}$. 
As a consequence and (iii) of $G_{A,B}$ Formation Procedure, when $e_p>0$, 
we have 
$
e_R(v_{p-1}, v_p)+e_R(v_{p-1}, v_{p+1}) \le  \lceil \frac{1}{2}e_{G_3}(v_{p-1}, v_p) \rceil + \lceil \frac{1}{2}e_{G_3}(v_{p-1}, v_{p+1}) \rceil -\ell. 
$ 
Let $\beta=\lceil \frac{1}{2}e_{G_3}(v_{p-1}, v_p) \rceil + \lceil \frac{1}{2}e_{G_3}(v_{p-1}, v_{p+1}) \rceil+\lceil \frac{1}{2}e_{G_3}(v_{p}, v_{p+1}) \rceil$. 
Now by the definition of $e_p$, we get 
\begin{eqnarray*}
e(R[\{v_{p-1}, v_p, v_{p+1}\}]) &\le& \lceil \frac{1}{2}e_{G_3}(v_{p-1}, v_p) \rceil + \lceil \frac{1}{2}e_{G_3}(v_{p-1}, v_{p+1}) \rceil -\ell +|F_1|\\
&=& \lceil \frac{1}{2}e_{G_3}(v_{p-1}, v_p) \rceil + \lceil \frac{1}{2}e_{G_3}(v_{p-1}, v_{p+1}) \rceil+\lceil \frac{1}{2}e_{G_3}(v_{p}, v_{p+1}) \rceil-\ell -e_p\\
&=& \max\left\{e(G_3[\{v_{p-1}, v_p, v_{p+1}\}])-k,  \beta-\lfloor 6\eta n \rfloor\right\} -\ell \\
&\le & \max\left\{ \Delta-k-\ell,  \Delta/2+3-\lfloor 6\eta n \rfloor-\ell\right\}\quad \text{(by~\eqref{claim:degree-of-v-p+1})}\\
&=& \Delta-k-\ell=\Delta(R). 
\end{eqnarray*}
Thus~\eqref{eqn:v_p,p-1,p=1-edge-in-R} holds. 

\begin{center}
	{\bf Step 5: Coloring the multigraph $R-F_1$}
\end{center}

Since all edges of $G_A-E_1$ and all edges of $G_B-F_1$ are already colored, we know that 
$ R-(E_1 \cup F_1)$ is bipartite. 
Let  $R^*=R-F_1$, where note  that $F_1= E_R(v_{p}, v_{p+1})$. 
As all edges of $E_1$ are incident with $u_{p+1}$ in $R$,   it follows that $R^*$ is nearly bipartite. 
Since $R$ is regular and $F_1= E_R(v_{p}, v_{p+1})$, we know that $\Delta(R^*)=\Delta(R)$. 

We may assume $|E_1|>0$, as otherwise $R^*=R$ is bipartite
regular with degree $\Delta(R)=\Delta-k-\ell$. By Theorem~\ref{konig}, we can then color the
edges of $R$ with $\Delta(R)$ colors from $[k+\ell+1, \Delta]$. 
Thus $\chi'(G_3) \le k+\ell +(\Delta-k-\ell)=\Delta$.  As $G_2\subseteq G_3$, we get $\chi'(G_2)=\Delta$
and so $\chi'(G)=\Delta$ by Claim~\ref{claim:G2-properties}(i). This gives a contradiction to the assumption that  $\chi'(G)=\Delta+1$.

Thus $|E_1|=|F_1|>0$.
We show that $R^*$ contains no $\Delta(R^*)$-overfull subgraph, and thus getting $\chi'(R^*)=\Delta(R^*)$ by
Theorem~\ref{thm:chromatic-index-nearly-bipartite-graph}.  
  To show that $R^*$ contains no $\Delta(R^*)$-overfull subgraph, 
 it suffices to show that $R$ contains no $\Delta(R)$-overfull subgraph. Suppose to the contrary that $R$ contains a $\Delta(R)$-overfull subgraph. 
That is, there exists $X\subseteq V(R)$ such that both $R[X]$ and $R-X$
are $\Delta(R)$-overfull (as $R$ is regular). Since $R-E_1-F_1$ is bipartite and has no $\Delta(R)$-overfull subgraph, 
it follows that $u_{p+1}$ is contained in exactly one of $R[X]$ and $R-X$ and $v_{p}, v_{p+1}$ 
are contained in the other one.  Without loss of generality, we assume that $u_{p+1} \in X$
and  $ v_p, v_{p+1}\in  V(R)\setminus X$.

We first notice that 
for any $u\in V(R)$, since   at most $9\eta n+5\eta^{\frac{1}{2}}n$  edges of $H-G_{A,B}$ incident with 
$u$ in $G_3$ were colored in Step 3 
by Condition~S3.\ref{s33}, and   $\ell$ edges of $H-G_{A,B}$ incident with $u$ in $G_3$ were colored in 
Step 4,  we have 
\begin{eqnarray}\label{eqn:R-degree}
	d_{R}^s(u)&\ge & \frac{1}{2}\left(d^s_{G_3}(u)-(n^{\frac{2}{3}}+4\eta n +2)\right)-(9\eta n+5\eta^{\frac{1}{2}}n)   -\ell \nonumber\\
	&>&  \frac{1}{2}d^s_{G_3}(u)-10.5 \eta^{\frac{1}{2}} n. 
\end{eqnarray}
Thus 
 \begin{numcases}{ d_{R}^s(u)\ge }
	\frac{1}{2}\Delta -11 \eta^{\frac{1}{2}} n  &  \text{if $u\in V(R)\setminus U^*$;}  \label{eqn:non-U-vertex-degree-inR}\\
	\frac{5}{98} n-10.5 \eta^{\frac{1}{2}} n> \frac{1}{22}n &  \text{if $u\in U$}. \label{eqn:U-vertex-degree-inR}
\end{numcases}

{\bf \noindent Claim 1: $|X\cap A|=|X\cap B|+1$ and $|B\setminus X|=|A\setminus X|+1$. }
If $|X\cap A| \le \frac{|X|-1}{2}$, then we get $e(R[X]) \le \Delta(R)\frac{|X|-1}{2}$, 
a contradiction to $R[X]$ being $\Delta(R)$-overfull. 
We also must have $|X\cap B| >\frac{|X|-3}{2}$.  For otherwise,  we have $e(R[X]) \le \Delta(R)\frac{|X|-3}{2}+|E_1| \le  \Delta(R)\frac{|X|-1}{2}$, where the last inequality is obtained by $|E_1| <\frac{\Delta}{6}<\Delta(R)$,
as $\Delta(R)=\Delta-k-\ell > \frac{\Delta}{2}-6\eta^{\frac{1}{2}} n$.  
Thus  $|X\cap A|=|X\cap B|+1$. This also gives $|B\setminus X|=|A\setminus X|+1$.

{\bf \noindent Claim 2: Let  $Z\in \{X, V(R)\setminus X\}$ such that  $|Z| \le m/2$. Then  $|Z|  \le 11$. }
We first show that 
$|Z| < 2\eta n+6$. 
For otherwise, we have $|Z|  \ge  2\eta n+6$ and so $|Z\setminus U^*| \ge 5$. Since $v_p, v_{p+1}\in U^*$, 
it follows that $|Z\setminus (U^*\cup \{u_{p+1}, v_p, v_{p+1}\})|  \ge 4$. 
As $|Z| \le m/2$, $|X\cap A|=|X\cap B|+1$, and $|B\setminus X|=|A\setminus X|+1$, 
it follows that $|Z\cap A|, |Z\cap B| \le m/4+1$, and so every vertex $u\in Z\setminus (U^*\cup \{u_{p+1}, v_p, v_{p+1}\})$
can have in $R$ at most $m/4+1$ neighbors from $Z$. Thus 
\begin{eqnarray*}
	e_{R}(Z, V(R)\setminus Z)& \ge&  4\left(\frac{1}{2}\Delta-11 \eta^{\frac{1}{2}} n -(m/4+1)\right)\\ 
	&> & n >\Delta(R), 
\end{eqnarray*}
a contradiction to $R[Z]$ being $\Delta(R)$-overfull.
By~\eqref{eqn:U-vertex-degree-inR}, if $|Z| \ge 13$, then 
we get 
$e_{R}(Z, V(R)\setminus Z) \ge  12 (\frac{1}{22}n -(|Z|+1)/2) >n/2>\Delta(R)$, a contradiction to $R[Z]$ being $\Delta(R)$-overfull. 
Thus $|Z| \le 12$.  Since $|Z|$ is odd, we have $|Z| \le 11$.

\bigskip 

If $|X| \le 11$, then as $u_{p+1}\in V(R)\setminus U^*$ and $X\cap B$ contains a vertex of $V(G)$,  we have 
\begin{eqnarray*}
	e(R[X])& \le&  \Delta(R)|X\cap B|+6-\left(\frac{1}{22}n -(|X|-1)/2\right)< \Delta(R)\frac{|X|-1}{2}, 
\end{eqnarray*}
a contradiction to $R[X]$ being $\Delta(R)$-overfull.
Thus $|X| \ge m/2+1$. 
Since $|X| \ge m/2+1$, we have 
  $|V(R)\setminus X| \le 11$ by Claim 2.

  {\bf \noindent Claim 3: We have  $e(R[\{v_p, v_{p+1}, u\}]) \le \Delta(R)$ for any $u\in V(R)\cap A$.} Consider 
  first that  $u=v_{p-1}$. If $e_p>0$, then $e(R[\{v_p, v_{p+1}, u\}]) \le \Delta(R)$  by~\eqref{eqn:v_p,p-1,p=1-edge-in-R}. 
  If $e_p=0$, then we have $e(R[\{v_p, v_{p+1}, u\}]) \le 3\times 6\eta^{1/2} n +3<\Delta(R)$  by the definition of $e_p$ and (iii) 
  of  $G_{A,B}$ Formation Procedure. 
  Thus we assume $u\ne v_{p-1}$. Then by the construction of $G_3$, we have $e_R(v_p, u)\le e_{G_3}(v_p,u) \le 1$, 
  and $e_R(v_{p+1}, \{v_p, u\}) \le \frac{1}{2}e_{G_3}(v_{p+1}, \{v_p, u\}) +1 \le \Delta/3+1$ (by Claim~\ref{claim:G2-properties}(v) and (ii) and (iii) of $G_{A,B}$ Formation Procedure). Thus $e(R[\{v_p, v_{p+1}, u\}]) \le \Delta/3+2 <\Delta(R)$.

   {\bf \noindent Claim 4: There are no two distinct vertices $u,v\in V(R)\setminus X$ such that 
   	$d^s_{G_3}(u), d^s_{G_3}(v) \ge \frac{\Delta}{2}$. } Suppose to the contrary that such $u$ and $v$ exist. By Claim 3, we have  $|V(R)\setminus X| \ge 5$.  
   Thus $|(V(R)\setminus X)\setminus \{v_1, ,u, v\}|\ge 2$. By~\eqref{eqn:R-degree} and~\eqref{eqn:U-vertex-degree-inR}, we have 
   	$e_{R}(V(R-X), X) \ge 2(\frac{\Delta}{4}-10.5\eta^{\frac{1}{2}} n -(|V(R)\setminus X|+1)/2)+2( \frac{1}{22}n-(|V(R)\setminus X|+1)/2)>\Delta(R)$, 
   	a contradiction to $R-X$ being $\Delta(R)$-overfull.

  By the construction of $L$ and $G_3$, we know that $N_{L^*}(v_p) \subseteq \{v_{p-1}, v_{p+1}\}$. 
  Thus $d_{R-X-\{v_{p-1},v_{p+1}\}}(v_p) \le 5$. 
  For the vertex $v_{p+1}$, if we have $\mu_{R-X-\{v_{p-1},v_{p}\}}(v_{p+1}) \le \frac{n}{23}$, then as  $B\setminus X \subseteq V(G)$ and $|A\setminus X|<|B\setminus X|$ by Claim 1, 
  we have 
  \begin{eqnarray*}
  	e(R-X) &\le&  e(R[\{v_p, v_{p+1}, u\}]) +d_{R-X-\{v_{p-1},v_{p+1}\}}(v_p)+d_{R-X-\{v_{p-1},v_{p}\}}(v_{p+1}) \\
  	&&+(|B\setminus X|-2)\Delta(R) -|B\setminus X|(\frac{n}{22}-|A\setminus X|) \\
  	& \le & \Delta(R)+5+|A\setminus X|\times \frac{n}{23}+(|B\setminus X|-2)\Delta(R) -|B\setminus X|(\frac{n}{22}-5)\\
  	&<& \Delta(R)\frac{|V(R-X)|-1}{2}, 
  \end{eqnarray*}
a contradiction to $R-X$ being $\Delta(R)$-overfull.

Thus we assume that $\mu_{R-X-\{v_{p-1},v_{p}\}}(v_{p+1}) >  \frac{n}{23}$. Since $\mu_{G_3}(u) \le 3\eta n$ for any $u\in V(G_3)\setminus U^*$ by~\eqref{eqn: vertex-multiplicity-not-from-U*} and $N_{L^*}(v_{p+1}) \cap (U^*\setminus \{v_{p-1}, v_p\}) \subseteq \{v_{2i-1}: i\in [1,p/2-1]\}$, 
it follows that $(A\setminus  X)\cap  \{v_{2i-1}: i\in [1,p/2-1]\} \ne \emptyset$.  Let $v_{i_1}, \ldots, v_{i_s}$
be all the elements of  $(A\setminus  X)\cap  \{v_{2i-1}: i\in [1,p/2-1]\}$, where  $i_1, \ldots i_s \in \{2i-1: i\in [1,p/2-1]\}$ and  $s \le 4$ by $|V(R)\setminus X| \le 11$ and Claim 1. 
If there exists $j\in [1,s]$ such that $v_{i_j+1} \not\in B\setminus X$,  then as 
$d_{G}(v_{p-1})<\frac{\Delta}{2}$ by Claim 4 (for otherwise, $d^s_{G_3}(v_p) \ge d_G(v_p) \ge d_G(v_{p-1}) \ge \frac{\Delta}{2}$
and $d^s_{G_3}(v_{p-1}) \ge d_{G}(v_{p-1}) \ge \frac{\Delta}{2}$),   and  $d_G(v_1) \le \ldots \le d_G(v_{p+1})$, 
we have  $d_G(v_{i_j+1})<\frac{\Delta}{2}$ 
and so $e_R(v_{i_j}, v_{i_j+1}) \ge \frac{1}{2}e_{G_3}(v_{i_j}, v_{i_j+1})-(9\eta n+5\eta^{\frac{1}{2}}n)   -\ell>\frac{\Delta}{4}-11\eta^{\frac{1}{2}} n$. 
As a consequence, we have 
 \begin{eqnarray*}
	e(R-X) &\le&  |A\setminus X|\Delta(R)+|F_1|-e_R(v_{i_j}, v_{i_j+1}) <  \Delta(R)\frac{|V(R-X)|-1}{2},
\end{eqnarray*}
as $|F_1|<\frac{\Delta}{6}<e_R(v_{i_j}, v_{i_j+1}) $. We again get a contradiction to $R-X$ being $\Delta(R)$-overfull.
Therefore $v_{i_1+1} , \ldots, v_{i_s+1}\in B\setminus X$.  Then by the construction of $L$, we know that $d_G(v_{i_j+1}) \ge d_G(v_{i_j})+e_{G_1}(v_{i_j}, v_{p+1})$. 
Since $e_{G_1}(v_{i_j}, v_{p+1}) \ge e_{G_3}(v_{i_j}, v_{p+1})-2\eta n$ (because of the Vertex Identification Procedure),  we get 
$e_R(v_{i_j+1}, X) >\frac{1}{2}e_{G_3}(v_{i_j}, v_{p+1}) -11\eta^{\frac{1}{2}}n$ by~\eqref{eqn:R-degree}.   By (iii) of  $G_{A,B}$ Formation Procedure, we  have $d_{R-X-\{v_{p-1},v_{p}\}}(v_{p+1}) \le \sum_{j=1}^s(\frac{1}{2}e_{G_3}(v_{i_j}, v_{p+1}) +1)+3\eta n (5-s)$. 
As $v_p, v_{p+1} \in U$ and $s\le 4$, we have 
 \begin{eqnarray*}
	e(R-X) &\le&  e(R[\{v_p, v_{p+1}, u\}]) +d_{R-X-\{v_{p-1},v_{p+1}\}}(v_p)+d_{R-X-\{v_{p-1},v_{p}\}}(v_{p+1}) \\
	&&+(|B\setminus X|-2)\Delta(R) -2(\frac{n}{22}-5) -\sum_{j=1}^s(\frac{1}{2}e_{G_3}(v_{i_j}, v_{p+1}) -11\eta^{\frac{1}{2}}n-5)\\
	& \le & \Delta(R)+5+s+3\eta n(5-s) +(|B\setminus X|-2)\Delta(R) -2(\frac{n}{22}-5)+s(11\eta^{\frac{1}{2}}n+5)\\
	&<& \Delta(R)\frac{|V(R-X)|-1}{2}, 
\end{eqnarray*}
a contradiction to $R-X$ being $\Delta(R)$-overfull.

Therefore, $R^*$ contains no $\Delta(R^*)$-overfull subgraph. Hence,  $\chi'(R^*)=\Delta(R^*)$ by
Theorem~\ref{thm:chromatic-index-nearly-bipartite-graph}.  
This implies that $\chi'(G_3-F_1) \le k+\ell +\Delta(R^*)=\Delta$. 
Since $G_2$ is a subgraph of $G_3-F_1$, we then get  $\chi'(G_2)=\Delta$. 
This further implies that $\chi'(G)=\Delta$ by Claim~\ref{claim:G2-properties}(i), showing
a contradiction to the assumption that $\chi'(G)=\Delta+1$. The proof of Theorem~\ref{thm:1}
is now complete. 
\qed

\bibliographystyle{abbrv}
\bibliography{SSL-BIB_08-19}
\end{document}